\documentclass[preprint,3p,11pt,sort&compress]{elsarticle}

\usepackage{hyperref}
\hypersetup{backref, colorlinks=true, linkcolor=blue}

\usepackage{stfloats}
\usepackage{subcaption}

\usepackage{pict2e}
\usepackage{xcolor}

\usepackage{graphicx} 
\usepackage{booktabs} 
\usepackage{microtype} 
\usepackage{amssymb,amsmath}
\usepackage[ruled,vlined]{algorithm2e}

\usepackage{bm}
\usepackage{amsthm}

\theoremstyle{plain}
\newtheorem{theorem}{Theorem}[section]

\usepackage{placeins}
\usepackage{longtable}

\begin{document}

\begin{center}
{\bf \Large Multi-Level Hybrid Monte Carlo / Deterministic  Methods \\
\vspace{0.1cm}
	for  Particle Transport Problems}
\end{center}

\author[ncsu]{Vincent N. Novellino}
\ead{vnnovell@ncsu.edu}
\author[ncsu]{Dmitriy Y. Anistratov}
\ead{anistratov@ncsu.edu}
\address[ncsu]{Department of Nuclear Engineering,
North Carolina State University Raleigh, NC}

\begin{frontmatter}
	\begin{abstract}
This paper presents multilevel hybrid transport (MLHT) methods for solving the neutral-particle Boltzmann transport equation. 
The proposed MLHT methods are formulated on a sequence of spatial grids using a multilevel Monte Carlo (MLMC) approach. 
The general MLMC algorithm is defined by recursively estimating the expected value of the correction to a solution functional on a neighboring grid. 
MLMC theory optimizes the total computational cost for estimating a functional to within a target accuracy.
The proposed MLHT algorithms are based on the quasidiffusion (variable Eddington factor) and second-moment methods. 
For these methods, the low-order equations for the angular moments of the angular flux are discretized in space. 
Monte Carlo techniques compute the closures for the low-order equations; then the equations are solved, yielding a single realization of the global flux solution. 
The ensemble average of the realizations yields the level solution.
The results for 1-D slab transport problems demonstrate weak convergence of the functionals.
We observe that the variance of the correction factors decreases faster than the computational cost of generating an MLMC sample increases. 
In the problems considered, the variance and cost of the MLMC solution are driven by the coarse-grid calculations.
\end{abstract}

\end{frontmatter}

\section{Introduction}

Solving the particle transport problem using a Monte Carlo method involves simulating many particles and collecting information about their interaction histories in tally bins split by space, time, energy, and angle.
Statistical uncertainty in the computed quantities depends on the number of particles simulated.
Estimators for the mean value of these quantities will converge by the Central Limit Theorem $\frac{\big<{X}_N\big> - \mu}{\sigma\sqrt{N}} \xrightarrow{d} \mathcal{N}(0,1)$, where $\big<{X}_N\big> = \frac{1}{N}\sum_{n=1}^N X(\omega_n)$, $X$ is a random variable of interest, and $\omega_n \in \Omega$ is a random walk of a particle sampled from the collection of all random walks $\Omega$. 
Reducing the uncertainty in high-resolution simulations can be expensive, since it necessitates the generation of many particle histories to reduce the variance, and the overhead of tallying results in many small tally volumes cannot be ignored \citep{hoogenboom-martin-mc-2011}.
Suppose the resolution of the simulation is decreased; then the cost of an individual particle history decreases since fewer boundary crossing events occur, and larger elements may contain more tally events, meaning there is less statistical uncertainty in the solution.
However, the solution obtained has a lower spatial resolution than we desire.
What if the lower-variance, cheaper-to-generate solution can be used to reduce the computational effort required for refined computational models while maintaining the desired model resolution?
That is the goal of Multilevel Monte Carlo (MLMC) methods \citep{giles-AN-2015}.

The fundamental idea of MLMC is to use a computationally inexpensive, low-fidelity solution that can be solved with low statistical uncertainty, and then compute correction terms to remove the discretization error by applying a telescopic sum of solution differences between two grids.
The MLMC algorithm optimizes the simulation by changing the number of simulations requested on a sequence of computational grids to reduce the total variance estimate to within a threshold in the most computationally efficient manner.
The application of the basic MLMC idea is problem-dependent and requires developing methods that account for the problem's specific features.

MLMC was first applied to parametric integration problems using the method of dependent tests \citep{HEINRICH1998151, Heinrich-2000, heinrich-LSSC2001, frolov-chentsov-1963,sobol-1963}.
Then the MLMC method was applied to ordinary and partial differential equations with stochastic coefficients \citep{giles-OP-2008, cliffe-2011, Barth-mlmc-sPDE-2011,giles-AN-2015}.
MLMC methods for uncertainty quantification problems have been developed \citep{MFmethods-siam-review-2018, zhang-mlcl-uq-review-2020, Geraci-2017}.
The MLMC approach was applied to optimize a single functional during simulation, such as hydraulic conductivity for groundwater flow simulations \citep{cliffe-2011}.
A recent extension of this work used MLMC in conjunction with an aggregation-based algebraic multigrid (AMG) coarsening strategy to solve the Darcy equation with a stochastic permeability field \citep{fairbanks-2024}.
In addition, one could also optimize a vector of output functionals in a simulation by checking the convergence of each component \citep{giles-AN-2015}.
Some conditions under which MLMC can be applied include a decrease in variance as grid fidelity increases for the correction terms and a method that converges under grid refinement, i.e., the correction factors shrink as more computational levels are added due to the decreasing amplitude of differences in the numerical solution on neighboring grids.
The distribution of computational work should minimize the total computational cost compared to an MC simulation on the fine target grid.
For cases where the variance decreases faster than the computational cost grows, the cost will be minimized by allocating most of the effort to the coarsest grids.

Another efficient approach to improving the MC solution is to apply hybrid MC/deterministic (HMCD) methods. 
There exists a family of HMCD methods for solving particle transport problems that were developed for fission source convergence \citep{larsen-yang-nse-2008, Lee-Joo-Lee-Smith, wolters-nse-2013, Willert-Knoll-Kelley-Park} and to remove effective scattering events in Implicit MC calculations \citep{pozulp-mc2023}.
Hybrid-MC-$S_2$ and Hybrid-MC-$S_2$X methods for solving the fixed-source problem have been derived by formulating equations for the partial reaction rates to avoid the approximation errors introduced by energy and angle discretization \citep{wolters-mc2009, Wolters_2011}.
This family of methods uses MC to compute non-linear functionals that depend weakly on the high-order transport solution; thus, the closures can have lower variance than the scalar flux solution.
The Coarse Mesh Finite Difference (CMFD) and Quasidiffusion (QD) equations have been analyzed to demonstrate the variance reduction of the non-linear methods for eigenvalue problems \citep{Lee-Joo-Lee-Smith,larsen-yang-nse-2008, Wolters_2011}.

In this paper, we present a novel multilevel HMCD method for solving particle transport problems by combining a geometric multigrid algorithm in space with the MLMC approach \citep{multigrid,mg-tutorial}.
The HMCD schemes are based on the low-order equations of the  Quasidiffusion (QD)/Variable Eddington Factor (VEF) and Second Moment methods \citep{gol'din-cmmp-1964,auer-mihalas-1970,sm-1976}.
We consider problems in 1-D slab geometry with isotropic sources and scattering.
The low-order equations for the angular moments of the angular flux are discretized with a finite volume (FV)  scheme that is of second-order accuracy \citep{dya-vyag-vant-1986,dya-vyag-ttsp}.
The discretized QD and SM low-order equations are defined for the cell-average scalar flux. 
The functionals of the transport solution, which define the exact closures of the low-order equations, are cell-average quantities as well.
The closures for the low-order equations are computed using an MC algorithm.
A preliminary analysis of these single-level FV-based HMCD transport schemes demonstrated promising results in filtering statistical noise effects \citep{vnn-dya-ans-annual-2024}.

The proposed multilevel HMCD method is defined on multiple spatial grids.
At the $0^{th}$ level, the algorithm obtains the hybrid transport solution on the coarsest grid.
At each level, the algorithm computes the difference between the hybrid solutions on two neighboring grids. 
The solutions from different levels are prolongated on the given target grid of the problem. 
The final solution on the target grid is defined using a telescopic sum, inspired by the MLMC approach.
MLMC optimization is applied to calculate a set of grid functionals of the transport solution.
Based on the variance of the functionals and runtime estimates from an initial set of samples, the MLMC algorithm determines the optimal number of samples at each level to achieve a user-specified accuracy tolerance.

The remainder of this paper is organized as follows.
We formulate HMCD transport methods in  Section \ref{sec:hmcd}.
Section \ref{sec:mlmc} describes the essential elements of the MLMC approach.
The Multilevel Hybrid Transport (MLHT) methods are formulated in  Section \ref{sec:mlht}.
Section \ref{sec:mlmc-opt}  reviews the MLMC optimization algorithm and relevant theory.
Section \ref{sec:mlmc-mlht}  describes MLHT algorithms with the MLMC optimization procedure.
Numerical results are presented in Section \ref{sec:res}.
We conclude with a discussion in Section~\ref{sec:con}.

\section{Hybrid Transport Methods Based on Low-Order Equations for Moments \label{sec:hmcd}}

We consider the 1-D slab geometry, steady-state particle transport equation with isotropic scattering and source:
\begin{equation} \label{t-eq}
	\mu \frac{\partial \psi}{\partial x}(x,\mu) + \Sigma_t(x) \psi(x,\mu) = \frac{\Sigma_s(x)}{2} \int_{-1}^1 \psi(x,\mu') d\mu' + \frac{q(x)}{2} \, ,
\end{equation}
\[x\in D, \quad D = [0,X], \quad \mu \in [-1,1] \, ,
\]
\begin{equation*}
	\psi(0,\mu) = \psi_{in}^+ , \ \mu > 0 \, , \quad 
	\psi(X,\mu) = \psi_{in}^-, \ \mu < 0 \, .
\end{equation*}
$x$ is the location in the slab, $X$ is the length of the slab, $\mu$ is the cosine of the angle between the direction of particle motion and the $x$-axis, $\Sigma_t$ is the total cross-section, $\Sigma_s$ is the scattering cross-section, and $q$ is the external source.
$\psi$ is the angular flux, $\psi_{in}^{\pm}$ are the angular fluxes of incoming particles.
The neutron scalar flux and current are defined by {the angular moments of $\psi$ given by
\begin{equation}
	\phi(x) = \int_{-1}^1 \psi(x,\mu)d\mu, \quad J(x) = \int_{-1}^1 \mu \psi(x,\mu)d\mu \, ,
\end{equation}
respectively.

\subsection{Hybrid Quasidiffusion/VEF Method  \label{sec:LOQD}}
To formulate an HMCD transport method, we apply the QD/VEF  method \citep{gol'din-cmmp-1964,auer-mihalas-1970}.
The	 low-order QD (LOQD) equations  for the scalar flux and current  are derived 
by taking the zeroth and first angular moments of the transport equation (Eq. \eqref{t-eq})
and formulating an exact nonlinear closure defined by means of the high-order transport solution.
The  LOQD equations are given by:
\begin{equation}
   \frac{dJ}{dx}(x) +   \big(\Sigma_t(x) - \Sigma_s(x) \big) \phi(x) = q(x) \; ,
\label{loqd1}
\end{equation}
\begin{equation}
  \frac{d}{dx}\big( E(x)\phi(x) \big) + \Sigma_t(x) J(x)  = 0 \; ,
\label{loqd2}
\end{equation}
where the  closure for the second moment $\int_{-1}^{1} \mu^2 \psi d \mu$ in the first moment equation (Eq. \eqref{loqd2}) 
	is defined by the QD (Eddington) factor
\begin{equation} \label{QDf}
   E(x) = \frac{ \displaystyle{\int_{-1}^{1} \mu^2 \psi(x,\mu) d \mu} }{ \displaystyle{\int_{-1}^{1} \psi(x,\mu) d \mu} }\, .
\end{equation}
The boundary conditions are given by  \citep{gol'din-cmmp-1964,gol'din-chetverushkin-1972}:
\begin{equation}
J(0)=  B_L (\phi(0) - \phi_{in}^+) + J_{in}^+\, ,
\quad
J(X)= B_R (\phi(X) - \phi_{in}^-) + J_{in}^- \; ,
\label{loqd-bcs}
\end{equation}
where 
\begin{equation}
     B_L  = \frac{ \displaystyle{\int_{-1}^{0}\mu \psi(0,\mu) d\mu}}
    {\displaystyle{\int_{-1}^{0} \psi(0,\mu) d\mu}}\; , \quad
     B_R = \frac{ \displaystyle{\int_{0}^{1} \mu \psi(X,\mu) d\mu}}
    {\displaystyle{\int_{0}^{1} \psi(X,\mu) d\mu}}\; ,
 \label{eqn:qd_bf}
\end{equation}
are the boundary QD factors.
 The partial fluxes and currents at  boundaries  are defined by the incoming angular flux distribution:
\begin{equation}
\phi_{in}^{\pm} = \pm \int_0^{\pm 1} \psi_{in}^{\pm }(\mu) d \mu \, , \quad
J_{in}^{\pm} = \pm \int_0^{\pm 1} \mu \psi_{in}^{\pm }(\mu) d \mu \, \, .
\end{equation}

We discretize the LOQD equations (Eqs. \eqref{loqd1} and \eqref{loqd2}) by a second-order finite volume (FV) method \citep{dya-vyag-vant-1986,dya-vyag-ttsp}.
We define the spatial grid $\{x_i \}_{i=0}^I$ and assume that cross sections and source are piece-wise functions over the set of spatial cells $\{\tau_i \}_{i=1}^I$, where $\tau_i=[x_{i-1},x_i]$.
 The balance equation (Eq. \eqref{loqd1}) is integrated over the $i^{th}$ spatial cell  to obtain
\begin{equation}\label{loqd1-d}
    J_i- J_{i-1} + \big(\Sigma_{t,i} - \Sigma_{s,i} \big)\Delta x_i \phi_i = q_i \Delta x_i  \, , \quad i \in \mathbb{N}(I) \, ,
\end{equation}
where $\Sigma_{t,i}$, $\Sigma_{s,i}$, and  $q_i$ are cross sections and the source in $\tau_i$,
$\Delta x_i =x_i - x_{i-1}$ is the cell width, $J_i = J(x_i)$ is the cell-edge current,
\begin{equation}
	 \phi_i  = \frac{1}{\Delta x_i} \int_{x_{i-1}}^{x_i} \phi dx \, 
\end{equation}
is the cell-average scalar flux. 
The first moment equation (Eq. \eqref{loqd2}) is integrated over $[\bar x_{i-1}, \bar x_i]$, $ i \in \mathbb{N}(I+1)$, where $\bar x_i=0.5(x_i + x_{i-1})$ for $ i \in \mathbb{N}(I)$, $\bar x_0 = x_0$ and $\bar x_{I+1}= x_I$. The FV discretization of Eq. \eqref{loqd2} is given by
\begin{equation}\label{loqd2-d}
    E_i \phi_i - E_{i-1} \phi_{i-1}+ \hat \Sigma_{t,i}  \Delta \hat x_{i} J_i  = 0,   \,  \quad i \in \mathbb{N}(I+1) \,  .
\end{equation}
Here $\phi_0 = \phi(x_0)$, $ \phi_{I+1} = \phi(x_I) $,
\begin{equation}  
   E_{i}(x) = \frac{\displaystyle{ \int_{-1}^{1}  \int_{x_{i-1}}^{x_i}  \mu^2 \psi(x,\mu) dx d \mu} }{\displaystyle{ \int_{-1}^{1}   \int_{x_{i-1}}^{x_i}\psi(x,\mu) dx d \mu }} \, ,  \quad i \in \mathbb{N}(I)
\end{equation}
is the cell-average QD (Eddington) factor,
$E_0 = E(x_0)$, $E_{I+1} = E(x_I)$ are the factors at the boundaries,
\begin{equation}
    \hat \Sigma_{t,i} = \frac{\Sigma_{t,i}\Delta x_{i} + \Sigma_{t,i-1}\Delta x_{i-1}}{\Delta x_i + \Delta x_{i-1}}, \quad \Delta \hat x_i  = \frac{1}{2}(\Delta x_{i} + \Delta x_{i-1}) \, .
\end{equation}
The boundary conditions have the form:
\begin{equation} \label{loqd-bcs-d}
    J_0=  B_L (\phi_0 - \phi_{in}^+) + J_{in}^+\, ,
	\quad
	J_I= B_R  (\phi_{I+1} - \phi_{in}^-) + J_{in}^- \, .
\end{equation}
The hybrid QD (HQD) method is defined by the LOQD system of equations discretized by the FV scheme \eqref{loqd1-d}, \eqref{loqd2-d}, and \eqref{loqd-bcs-d} with the QD and boundary factors computed by MC.

\subsection{Hybrid Second Moment Method \label{sec:LOSM}}
 
Another HMCD method is formulated based on the Second Moment (SM) method \citep{sm-1976,mla-ewl-pne-2002}.
The low-order SM (LOSM) equations are  derived from the zeroth and first moments of the transport equation with a linear closure and are given by
\begin{equation} \label{losm1}
 	\frac{dJ}{dx}(x) +   \big(\Sigma_t(x) - \Sigma_s(x) \big) \phi(x) = q(x) \; ,
\end{equation}
\begin{equation}\label{losm2}
    \frac{1}{3} \frac{d \phi}{dx} (x) + \Sigma_t(x) J(x)  =    \frac{dH}{dx}(x),
\end{equation}
where the exact closure is defined with
\begin{equation}
 	H(x) = \frac{1}{3}\int_{-1}^{1}  (1- 3 \mu^2) \psi(x,\mu) d \mu \, .
\end{equation}
The boundary conditions are given by
\begin{equation} \label{losm-bcs}
    J(0)= -\frac{1}{2}\phi(0)   +2 J_{in}^+ + W_L\, ,
    \quad
    J(X)= \frac{1}{2}\phi(X)   + 2 J_{in}^- -  W_R\,  ,
\end{equation}
where the boundary functionals are defined as follows:
\begin{equation}\label{sm-boundary}
    W_L =  \frac{1}{2}\int_{-1}^{1}  (  1- 2 |\mu|  )  \psi(0,\mu) d\mu  \, , \quad
    W_R =  \frac{1}{2} \int_{-1}^{1} (1-  2 |\mu|   )  \psi(X,\mu) d\mu  \, .
\end{equation}
To discretize the LOSM equations, we apply a similar FV scheme as described in Section \ref{sec:LOQD}
 to obtain
\begin{equation}\label{losm1-d}
    J_i- J_{i-1} + \big(\Sigma_{t,i} - \Sigma_{s,i} \big)\Delta x_i \phi_i = q_i \Delta x_i  \, , \quad i \in \mathbb{N}(I) \, ,
\end{equation}
\begin{equation}\label{losm2-d}
    \frac{1}{3} \big( \phi_i -   \phi_{i-1} \big) + \hat \Sigma_{t,i}  \Delta \hat x_{i} J_i  =  H_i - H_{i-1}
     \quad i \in \mathbb{N}(I+1) \, ,
\end{equation}
\begin{equation} \label{losm-bcs-d}
    J_0= -\frac{1}{2}\phi_0   +2 J_{in}^+ + W_L\, ,
    \quad
    J_{I+1}= \frac{1}{2}\phi_{I+1}   + 2 J_{in}^- -  W_R\,  ,
\end{equation}
where
\begin{equation}
 		H_i = \frac{1}{3 \Delta x_i}\int_{-1}^{1} \int_{x_{i-1}}^{x_i} (1- 3 \mu^2) \psi(x,\mu) dx  d \mu \,,  \quad i \in \mathbb{N}(I) \, ,
\end{equation}
\begin{equation}
        H_0=H(x_0) \, , \quad   H_{I+1}=H(x_I) \, .
\end{equation}
The hybrid SM (HSM) method is formulated using the approximated LOSM equations \eqref{losm1-d}-\eqref{losm-bcs-d}, with the closure term $H$ and boundary functionals computed via MC.

\subsection{MC Estimators of Closure Functionals for Low-Order Equations} \label{sec:tallies}

To calculate functionals for the closures of the low-order equations that define HMCD methods, we collect scores in spatial cells on a computational grid to compute the corresponding tally quantities.
Track-length-based tallies are used for estimators of closure functionals $E_i$ and $H_i $ as well as of the scalar flux $\phi_i$.
The track-length tally estimators of the $r^{th}$ angular moment $\int_{-1}^1 \int_{\tau_{i}}\mu^r \psi(x,\mu)dx d\mu$ in the $\tau_i$ cell is given by
\begin{equation}\label{eqn:tracklength_estimator}
    T^{[r]}_i = \frac{1}{K} \sum_{k=1}^K \sum_{m=1}^{M_k} \mu_{k,m}^r  w_{k,m} \nu_{k,m} \, ,
\end{equation}
where  $k$ is the particle index, $\nu_{k,m}$ is the track-length of the $k^{th}$ particle in the cell $\tau_i$ traveling in the direction  $\mu_{k,m}$,  $w_{k,m}$ is the particle weight, $M_k$  number of $k^{th}$ particle tracks in the $i$-th cell, $K$ is the number of source particles.
As a result, we define estimators
\begin{equation}  \label{qdf-mc}
	\big< E \big>_i = \frac{T^{[2]}_i}{T^{[0]}_i} \, ,
\end{equation}
\begin{equation}\label{smf-mc}
	\big< H \big>_i = \frac{1}	{3\Delta x_i} (T_i^{[0]} - 3T_i^{[2]}) \, ,
\end{equation}
\begin{equation}
	\big<\phi\big>_i = \frac{1}{\Delta x_i} T_{i}^{[0]},
\end{equation}
for the cell-average values of $E$, $H$ and $\phi$ in the $i^{th}$ the cell $\tau_i$ ($i \in \mathbb{N}(I)$).

To compute functions at domain boundaries, we define the $r^{th}$  partial angular moment surface crossing tally at some plane located at $x^*$

\begin{equation}\label{eqn:surface_estimator}
    S_{x^*}^{[r]\pm} = \pm \int_{0}^{\pm 1} |\mu|^r \psi(x^{*},\mu) d\mu \approx \frac{1}{K} \sum_{k=1}^K \sum_{m=1}^{M_k}g(\pm\mu_{k,m}) |\mu_{k,m}|^{r-1}  w_{k,m} \bigg|_{x=x^*}
\end{equation}
where
\begin{equation}
  g(\mu) = 
    \begin{cases}
            1, \quad  \mu   > 0   \\
            0, \quad \mu \le 0 \,  
    \end{cases}
\end{equation}
  is a unit-step function.
The full range surface crossing tally is given by
\begin{equation}
    S_{x^*}^{[r]} = \int_{-1}^{ 1}  \mu^r \psi(x^{*},\mu) d\mu = S_{x^*}^{[r]+} + (-1)^r S_{x^*}^{[r]-}.
\end{equation}
The absolute surface crossing tally is defined by
\begin{equation}
    |S_{x^*}^{[r]}| =  \int_{-1}^{ 1}  |\mu|^r \psi(x^{*},\mu) d\mu = S_{x^*}^{[r]+} + S_{x^*}^{[r]-}.
\end{equation}
The QD factors at the boundaries and the boundary factors are estimated using face-crossing tallies of the form:
\begin{equation}\label{qdf-boundary-mc}
	\big<E_{0}\big>_0   = \frac{S_{0}^{[2]}}{S_{0}^{[0]}}\, , \quad
	\big<E \big>_{I+1} = \frac{S_{X}^{[2]}}{S_{X}^{[0]}}\, , \quad
	\big<B_{L}\big> = \frac{-S_{0}^{[1]-}}{S_{0}^{[0]-}}\, , \quad
	\big<B_{R}\big> = \frac{S_{X}^{[1]+}}{S_{X}^{[0]+}}\,.
\end{equation}
The second moment functionals on the boundaries and the boundary functionals are given by:
\begin{subequations}\label{smf-boundary-mc}
\begin{equation} 
	\big<H \big>_0    = \frac{1}{3} \Big(S_{0}^{[0]}  - 3S_{0}^{[2]} \Big) \, , \quad
	\big<H \big>_{I+1}  = \frac{1}{3} \Big(S_{X}^{[0]}  - 3S_{X}^{[2]}\Big) \, ,
\end{equation}
\begin{equation} 
	\big<W_{L}\big>  = \frac{1}{2} \Big(|S_{0}^{[0]}| - 2|S_{0}^{[1]}|\Big) \, , \quad
	\big<W_{R}\big>  = \frac{1}{3} \Big(|S_{X}^{[0]}| - 2|S_{X}^{[1]}|\Big) \, .
\end{equation}
\end{subequations}

\section{Basic Idea of MLMC \label{sec:mlmc}}

Consider a hierarchy  of spatial grids, $G_\ell$  for $\ell = 0, 1, \dots, L$ such that $G_0 \subset G_1 \subset G_2 ... \subset G_L$.
$G_0$ is the coarsest grid, $G_L$ is the finest grid.
Let $F$ be a functional of interest, and $F_\ell$ be an approximation of the functional on the grid $G_\ell$.
An estimator of the expected value $\mathbb{E}[ F_\ell]$ is defined  by
\begin{equation} \label{F-est}
	\big<F_\ell\big> = \frac{1}{N_{\ell}} \sum_{n=1}^{N_\ell} F_\ell(\omega^*_{n,\ell}),
\end{equation}
where $\omega^*_{n,\ell}$ is the $n^{th}$ random sample on $G_\ell$ coming from the probability space $(\Omega^*, \mathcal{F},P)$, $N_\ell$ is the number of samples for this grid.
$\Omega^*$ is the set of all Monte Carlo simulations comprised of $K$ particle samples.
The MLMC  approach is based on recursive estimation of the expected value of correction with respect to the neighboring grid  
$\mathbb{E}[F_{\ell} - F_{\ell-1}]$.
To compute $\mathbb{E}[ F_L ]$ on the finest grid $G_L$, the MLMC applies a telescoping sum \citep{giles-AN-2015}
\begin{equation} \label{E-tele-sum}
	\mathbb{E}[ F_L ]  = \mathbb{E}[ F_0 ]   + \sum_{\ell=1}^{L} \mathbb{E}[\Delta F_{\ell} ] \, ,
\end{equation}
where 
\begin{equation}
	\Delta F_{\ell}  =   F_{\ell} - F_{\ell-1} \, .
\end{equation}
At the $\ell^{th}$ level, the estimator of  $\mathbb{E}[\Delta F_{\ell}]$ is computed as follows:
\begin{equation} \label{delta-F-ell-est}
 	\big<\Delta F_{\ell} \big> =  	\frac{1}{N_{\ell}} \sum_{n=1}^{N_{\ell}}  \Delta F_{\ell}(\omega^*_{n,\ell}) =
 	\frac{1}{N_{\ell}} \sum_{n=1}^{N_{\ell}} \big(  F_{\ell}(\omega^*_{n,\ell}) - F_{\ell -1}(\omega^*_{n,\ell}) \big),
\end{equation}
\begin{equation}  
    \Delta F_{\ell}(\omega^*_{n,\ell}) =    F_{\ell}(\omega^*_{n,\ell}) - F_{\ell -1}(\omega^*_{n,\ell})  \, ,
\end{equation}
 where the same collection of particle histories $\omega^*_{n,\ell}$ is used to estimate the functional on both the  grid $G_\ell$ and its coarser neighboring grid  $G_{\ell-1}$.
This decreases the effects of statistical noise since the samples are correlated, yielding the following variance estimate for the correction terms \citep{frolov-chentsov-1963, Heinrich-2000} 
where the variance of $\big<\Delta F_{\ell} \big>$ is given by 

\begin{equation}
  \mathbb{V}(\Delta F_{\ell}(\omega^*_{n,\ell})) = \mathbb{V}(F_{\ell}(\omega^*_{n,\ell})) + \mathbb{V}(F_{\ell-1}(\omega^*_{n,\ell})) - 2 \cdot\text{COV}(F_{\ell}(\omega^*_{n,\ell}), F_{\ell-1}(\omega^*_{n,\ell})).
\end{equation}
The functional estimated on $G_{\ell-1}$ can be interpreted as a control variate for the given $\ell^{th}$ level.
The estimator of $\mathbb{E}[F_L]$ is given by
\begin{equation} \label{F-est-tele-sum}
	\big<{F}_L\big> = \big<F_0\big> + \sum_{\ell=1}^L \frac{1}{N_{\ell}}  \sum_{n=1}^{N_{\ell}} \Big(  F_{\ell}(\omega^*_{n,\ell}) - F_{\ell -1}(\omega^*_{n,\ell}) \Big),
\end{equation}
The estimations of $\mathbb{E}[\Delta F_{\ell} ]$   at all levels are performed independently. As a result,  the variance of  $\big<F_{L} \big>$ is given by
\begin{equation} \label{var-F-L}
	\mathbb{V}\big[ \big<F_{L} \big> \big] =  \mathbb{V}\big[\big<  F_0 \big> \big]  + \sum_{\ell=1}^{L} \mathbb{V}\big[\big<  \Delta F_{\ell}  \big>\big]  \, .
\end{equation}

The evaluation of each $F_0$ sample is computationally inexpensive on the coarse grid $G_{0}$.
Computational costs of $\Delta F_\ell$ increase from level to level due to grid refinement.
However, the magnitude of $\Delta F_\ell$ decreases with each level and the variance $\mathbb{V}[\Delta F_{\ell}]$ decreases.
Hence, fewer samples may be needed if the variance $\mathbb{V}[\Delta F_{\ell}]$ is far less than $\mathbb{V}[F_{\ell}]$.
MLMC enables one to minimize the computational cost for a given target accuracy of the functional $\big<{F}_L\big>$.
The number of independent random samples used in MLMC at different levels can be optimized, taking into account the rate of weak convergence, the decrease in variance, and the increase in computational cost.
The optimization algorithm of MLMC is described below in Section \ref{sec:mlmc-opt}.

\section{Multilevel Hybrid Transport Methods \label{sec:mlht}}

The multilevel hybrid transport (MLHT) methods are defined on a hierarchy of sequentially refined spatial grids $G_\ell=\{ x_{i,\ell}\}_{i=0}^{I_{\ell}}$ for $\ell = 0, 1, \dots, L$, where $I_{\ell} =  \rho I_{\ell-1}$ and $\rho>1 $ is the refining factor.
The spatial interval $\tau_{i,\ell}=[x_{i-1,\ell}, x_{i,\ell}]$ contains the corresponding intervals of the grid $G_{\ell+1}$.

Algorithm \ref{alg:mlht} presents the MLHT algorithm, which uses the HQD and HSM methods to formulate hybrid low-order problems. 
Hereafter, the MLHT algorithms based on the LOQD and LOSM equations are referred to as the MLHQD and MLHSM methods, respectively.
The MLHT algorithm consists of two stages:
\begin{itemize}
    \item {\it Stage 1:}  calculations at the level 0 on the coarsest grid $G_0$,
 	\item  {\it Stage 2:} calculations at  the  $\ell^{th}$ level ($\ell =1,\ldots,L$)  
 	solving a hybrid transport problem  on two neighboring grids $G_{\ell-1}$ and $G_{\ell}$.
\end{itemize}
An element of the algorithm at any $\ell^{th}$ level is solving the hybrid low-order problem to generate the vector of discrete scalar flux values $\boldsymbol{\phi}_{\ell}=\{ \phi_{i,\ell}\}_{i=0}^{I_{\ell}+1}$ on the grid $G_{\ell}$.
To compute the closure functionals for the low-order equations, the algorithm performs an MC simulation with $K_{\ell}$ particle histories and calculates the tallies necessary for the functionals. 
Each such low-order solve (i.e. one MC simulation followed by a hybrid solve) provides a single realization of the hybrid transport solution $\boldsymbol{\phi}_{n,\ell}=\{ \phi_{i,n,\ell}\}_{i=0}^{I_{\ell}+1}$ on $G_{\ell}$ obtained with $\{\mathcal{H}_{k,n,\ell}\}_{k=1}^{K_{\ell}}$ collection of particle histories, where $n$ is the realization index.
The number of realizations $N_{\ell}$ varies with levels and is defined according to the MLMC optimization procedure described below in Sections \ref{sec:mlmc-opt} and \ref{sec:mlmc-mlht}.
The final hybrid transport solution is the average over $N_{\ell}$ realizations used at that level.

In Stage 1, the hybrid solution at the level 0 is the ensemble average of realizations on $G_0$ and given by
\begin{equation}
		\big< \boldsymbol{\phi}_0 \big>  = \frac{1}{N_0} \sum_{n=1}^{N_0} \boldsymbol{\phi}_{n,0} \, .
\end{equation}
For Stage 2 and at the $\ell^{th}$ level ($\ell\ge 1$), the MLHT algorithm computes
\begin{equation}\label{delta-phi-ell}
	\big< \Delta \boldsymbol{\phi}_{\ell} \big> = 
	\frac{1}{N_{\ell}} \sum_{n=1}^{N_{\ell}} (\boldsymbol{\phi}_{n,\ell} -  \mathcal{I}_{\ell-1}^{\ell}\boldsymbol{\phi}_{n,\ell-1} ) \, ,
\end{equation}
where $\mathcal{I}_{\ell-1}^{\ell}$ is the prolongation operator  of the solution from $G_{\ell-1}$ to the finer grid $G_{\ell}$.
In this work, a constant prolongation operator is considered.
The $n^{th}$ realization of hybrid solutions on $G_{\ell}$ and $G_{\ell-1}$ are obtained from low-order equations on corresponding spatial grids with the closure functionals computed with the same $n^{th}$ ensemble of particle histories $\{\mathcal{H}_{k,n,\ell}\}_{k=1}^{K_{\ell}}$.
After moving through all levels ($\ell=0,...,L$), the hybrid solution on the spatial grid  $G_L$ is computed through   a telescopic summation   given by
\begin{equation} \label{tele-phi-L}
	\big< \boldsymbol{\phi}_L\big> =    \sum_{\ell=0}^L  \mathcal{I}_{\ell}^L \big< \Delta \boldsymbol{\phi}_\ell \big>  \,  , \quad
	\mbox{where} \quad 
\big< \Delta \boldsymbol{\phi}_0 \big> = 	\big< \boldsymbol{\phi}_0 \big> \, ,
\end{equation}
  which applies the prolongation of the numerical solutions from each level to $G_L$. 
The telescopic summation follows the basic idea of MLMC   (see Eq. \eqref{F-est-tele-sum})
and combines it with a geometric multigrid approach.
This formulates a multilevel estimator for the intermediate solution $\big< \Delta \boldsymbol{\phi}_{\ell} \big>$ and the solution on the finest grid $\big< \boldsymbol{\phi}_L\big>$.

The MLHT algorithm  requires calculations of tallies $\boldsymbol{T}_{n,\ell}$  for the closure functions $\big< \boldsymbol{\varGamma}_{n,\ell}  \big>$. 
The set of functionals for the MLHQD method includes $\big\{\big< E_i \big>\big\}_{i=0} ^{I_{\ell}}$,   $\big< B_L\big>$, $\big< B_R\big> \big\}$.
The tallies are defined by Eqs. \eqref{qdf-mc}-\eqref{qdf-boundary-mc}.
The set of functionals $\big< \boldsymbol{\varGamma}_{n,\ell}  \big>$ for the MLHSM method consists of $\big\{ \big< H_i \big>\big\}_{i=0} ^{I_{\ell}}$,   $\big< W_L\big>$,  $\big< W_R\big> \big\}$ given by Eqs. \eqref{smf-mc} and \eqref{smf-boundary-mc}.
At the $\ell^{th}$ level calculations of the MLHT algorithm, the ensemble of particle histories $\{\mathcal{H}_{k,n,\ell}\}_{k=1}^{K_{\ell}}$  for the $n^{th}$ realization is utilized to compute closure functions for low-order equations on both $G_{\ell}$ and $G_{\ell-1}$. 
The set of tallies $\boldsymbol{T}_{n,\ell}$ is combined to generate $\boldsymbol{T}_{n,\ell-1}$ at the $\ell^{th}$ level.
The coarse grid $G_{\ell-1}$ contains exactly $\rho$ fine cells, meaning scores collected on the fine grid can be added together to calculate the scores on the coarse grid for track-length-based tallies.
Boundary factors are the same for the coarse and fine grids.

 \medskip
\begin{algorithm}[h!]
	\caption{MLHT  Algorithm for  MLHQD and MLHSM Methods}\label{alg:mlht}
	\DontPrintSemicolon
	\SetAlgoVlined
    \texttt{$\bullet$ Stage 1 }\;
	\For{$n=1,\ldots, N_{0}$}{
 	\texttt{$\bullet$ $n^{th}$ realization at $0^{th}$ level}\;
		
		\For{$k=1,\ldots,K_0$}{
			simulate $k^{th}$  particle history  $\mathcal{H}_{k,n,0}$\;
	compute  tallies  $\boldsymbol{T}_{n,0}$ for closure functionals  on   $G_{0}$\;
		}
		compute closure functionals 	$\big< \boldsymbol{\varGamma}_{n,0}  \big>$  for $n^{th}$ realization on $G_0$\;
		solve the low-order equations defined with  $\big< \boldsymbol{\varGamma}_{n,0}  \big>$   for $\boldsymbol{\phi}_{n,0}$  on $G_0$\;
	}
	  $\big< \boldsymbol{\phi}_{0}\big>  =  N_0^{-1} \sum_{n=1}^{N_0}\boldsymbol{\phi}_{n,0}$\;

    \texttt{$\bullet$ Stage 2 }\;
	\For{$\ell=1,\ldots, L$}{
		
		\For{$n=1,\ldots, N_{\ell}$}{
 		\texttt{$\bullet$ $n^{th}$ realization at $\ell^{th}$ level}\;

		\For{$k=1,\ldots,K_{\ell}$}{
	simulate $k^{th}$  particle history  $\mathcal{H}_{k,n,\ell}$\;

	compute  tallies  $\boldsymbol{T}_{n,\ell}$ for closure functionals  on   $G_{\ell}$\;
}
			combine  scores on $G_{\ell}$  to compute  tallies $\boldsymbol{T}_{n,\ell-1}$ on   $G_{\ell-1}$\;

			compute closure functionals $\big< \boldsymbol{\varGamma}_{n,\ell}  \big>$ on $G_{\ell}$ 
			and $\big< \boldsymbol{\varGamma}_{n,\ell-1}  \big>$ on $G_{\ell-1}$\;
	
	     	solve  low-order equations  defined with  $\big< \boldsymbol{\varGamma}_{n,\ell-1}  \big>$ 
             for $\boldsymbol{\phi}_{n,\ell-1}$  on $G_{\ell-1}$\;	
     
         	solve  low-order equations  defined with  $\big< \boldsymbol{\varGamma}_{n,\ell}  \big>$ 
	     	  for $\boldsymbol{\phi}_{n,\ell}$  on $G_{\ell}$\;
		}
			 $\big< \Delta \boldsymbol{\phi}_{\ell} \big> = 
 N_{\ell}^{-1} \sum_{n=1}^{N_{\ell}} (\boldsymbol{\phi}_{n,\ell} -  \mathcal{I}_{\ell-1}^{\ell}\boldsymbol{\phi}_{n,\ell-1} )$\;
	}
	$\big< \boldsymbol{\phi}_L\big> =  \mathcal{I}_{0}^L \big< \boldsymbol{\phi}_0 \big>  +  \sum_{\ell=1}^L  \mathcal{I}_{\ell}^L \big< \Delta \boldsymbol{\phi}_\ell \big> $ \;
\end{algorithm}

\section{Multilevel Monte Carlo Optimization \label{sec:mlmc-opt}}

Let $F_{I}$ be a functional $F$ computed with a random vector, which is a numerical solution of a discretized PDE on a spatial grid with $I$ degrees of freedom (DF). 
In this study, we consider the discretized low-order particle transport equation with stochastic closure coefficients as the PDE of interest.
The spatial grid has $I$ intervals. 
In the case of a uniform grid, we have $I=X h^{-1}$, where $h$ is the width of spatial intervals.
The functional of interest is a moment of the angular flux, for example, the scalar flux.
We assume that
\begin{equation}
	\mathbb{E}[F_I]   \longrightarrow 	\mathbb{E}[F] \quad \mbox{as}  \quad I \to \infty \quad (h \to 0) \,  ,
\end{equation}
and 
\begin{equation}
	\mathbb{E}[F_I - F ]  =\mathcal{O}(h^{\alpha})  =  \mathcal{O}(I^{-\alpha}) \, .
\end{equation}
We note that the spatial discretization schemes for the LOQD and LOSM equations are second-order accurate and, hence, $\alpha=2$ provided that the approximation of the functional is at least of order 2.
The accuracy of an estimator $\big< F_I \big>=\frac{1}{N} \sum_{n=0}^N F_{I,n}$ is measured by the root mean square error (RMSE)
\begin{equation}
	RMSE\big( \big< F_I \big> \big) =  \sqrt{ \mathbb{E}\Big[ \big(\big< F_I \big> - F\big)^2 \Big] } \, .
\end{equation}
The mean square error (MSE) has the following form  \citep{giles-OP-2008, cliffe-2011, giles-AN-2015}:
\begin{equation} \label{mse-F_I}
	MSE\big( \big< F_I \big> \big) =  \mathbb{V}\big[ \big< F_I \big> \big] + 
\Big(\mathbb{E}\big[ \big< F_I \big>\big] - \mathbb{E}[F_I]\Big)^2 =
	 	\frac{1}{N}  \mathbb{V}\big[  F_I  \big]  +  \Big( \mathbb{E}\big[  F_I  - F] \Big)^2 \, ,
\end{equation}
where $N$ is the number of random samples. 
The equation \eqref{mse-F_I} shows contributions of stochastic and discretization errors to the estimator MSE.
Thus, to compute the functional with such accuracy that
\begin{equation}
	 RMSE\big( \big< F_I \big>   \big)  < \varepsilon,
\end{equation}
it is sufficient that
\begin{equation} \label{cond-1}
		\frac{1}{N}  \mathbb{V}\big[  F_I  \big]  < \frac{\varepsilon^2}{2} \, , \quad
		  \Big( \mathbb{E}\big[  F_I  - F] \Big)^2 < \frac{\varepsilon^2}{2}  \, .
\end{equation}
This leads to the following conditions on the number of samples and DF:
\begin{equation} \label{cond-2}
 N \gtrsim \varepsilon^{-2} \, , \quad 	I  \gtrsim \varepsilon^{-\frac{1}{\alpha}} \, , 
\end{equation}
provided that $\mathbb{V}\big[  F_I  \big] $ is constant and doesn't depend on $I$.
Here we use the notation $f  \gtrsim g$ for $f>0$ and $g>0$, indicating that the ratio $\frac{f}{g}$ is uniformly bounded and independent of the number of random samples and DF.

The cost $C_{I,n}$  of a single sample $F_{I,n}$ depends on DF. Assuming that
\begin{equation}
	C_{I,n}  \lesssim I^{\gamma}   \, , \quad  \mbox{and} 	\quad  C_{I,n}  \lesssim h^{-\gamma}  \, , \quad  \gamma> 0 \, ,
\end{equation}
then the cost of the estimator $\big< F_I \big>$  meets the condition
\begin{equation}
	C_{I}  \lesssim N I^{\gamma} \, . 
\end{equation}
Taking into account Eq. \eqref{cond-2},  the costs of computations for the given accuracy $\varepsilon$ satisfies
\begin{equation}
	C_{I}  \lesssim \varepsilon^{-2-\frac{\gamma}{\alpha}} \, .
\end{equation}

We now consider an MLMC algorithm on the sequence of grids $\{ G_{\ell}\}_{\ell=0}^L$ with $\{ I_{\ell}\}_{\ell=0}^L$ spatial intervals.  
In the case of a set uniform grids,   $I_{\ell}=X h_{\ell}^{-1}$ for the cell width $h_{\ell}$.
Let $F_{\ell}$ be an approximation of the functional $F$ by a hybrid solution of a PDE on $G_{\ell}$  estimated by Eq. \eqref{F-est}. 
The estimation of the functional on $G_L$ is  defined according to the MLMC method described above (see Section \ref{sec:mlmc}) and given by
\begin{equation} \label{F-est-tele-sum-1}
	\big<{F}_L\big> =  \sum_{\ell=0}^L \big< \Delta F_{\ell}\big> \, ,
\end{equation}
where   $\Delta F_{0} = F_{0}$.
The MSE of the MLMC estimator $\big< F_L \big>$ has the form  \citep{giles-OP-2008, cliffe-2011, giles-AN-2015}:
\begin{equation}
	MSE\big( \big< F_L \big> \big) = \mathbb{E}\Big[ \big( \big< F_L \big> - F \big)^2\Big]  =  
\sum_{\ell=0}^L	\frac{1}{N_{\ell}}  \mathbb{V}\big[  \Delta F_{\ell}  \big]  +  \Big( \mathbb{E}\big[  F_I  - F] \Big)^2 \, .
\end{equation}
Thus, the sufficient conditions for 
\begin{equation}
	RMSE\big( \big< F_L \big>   \big)  < \varepsilon
\end{equation}
are the following: 
\begin{equation} \label{v-part}
\sum_{\ell=0}^L	\frac{1}{N_{\ell}}  \mathbb{V}\big[  \Delta F_{\ell}  \big]  < \frac{\varepsilon^2}{2} \,  , 
\end{equation}
\begin{equation} \label{disc-part-err}
	\mathbb{E}\big[  F_L  - F] < \frac{\varepsilon}{2}  \, .
\end{equation}
To meet the condition on the approximation error (Eq. \eqref{disc-part-err}), it is sufficient that
\begin{equation}
	I_L \gtrsim \varepsilon ^{-\frac{1}{\alpha}} \quad \mbox{and} \quad
		h_L \lesssim \varepsilon ^{\frac{1}{\alpha}} \, .
\end{equation}

Let $C_{\ell}$ be the computational cost of one sample $\Delta F_{\ell} (\omega_{\ell,n})$ at the $\ell^{th}$ level, then the cost of the MLMC estimator is given by
\begin{equation}
	C_L = \sum_{\ell=0}^L N_{\ell} C_{\ell} \, .
\end{equation}
The variance of the MLMC estimator is minimized if   the number of samples at $\ell^{th}$ level  is defined by  \citep{giles-AN-2015}
\begin{equation} \label{estimation-N_ell}
	N_{\ell} = \displaystyle{\frac{2}{\varepsilon^2}}  
	 \sqrt{\frac{\mathbb{V}[ \Delta F_{\ell} ]}{  C_{\ell}} } \sum_{\ell=0}^{L} \sqrt{\mathbb{V}[\Delta F_{\ell}] C_{\ell} }  \,  .
\end{equation}

The performance of MLMC algorithms for discretized PDEs with stochastic coefficients is described by the following complexity theorem, which is based on conditions on the numerical solution approximation, properties of the multilevel estimator, and the computational cost of the algorithm's components \citep{cliffe-2011}. 
\begin{theorem}
\label{theorem}
Let   $\big< \Delta F_{\ell}\big>=	\frac{1}{N_{\ell}} \sum_{n=1}^{N_{\ell}} \big(  F_{\ell}(\omega^*_{n,\ell}) - F_{\ell -1}(\omega^*_{n,\ell}) \big)$  and assume that there are constants $\alpha>0$,  $\beta>0$,  $\gamma>0$ such that 
$\alpha \ge \frac{1}{2} min(\beta,\gamma)$,
and
\begin{equation}
\big| \mathbb{E} [ F_{\ell} - F ] \big| \lesssim  I_{\ell}^{-\alpha} \, ,
\end{equation} 
\begin{equation}
\mathbb{V}\big[ \Delta F_{\ell} \big] \lesssim  I_{\ell}^{-\beta} \, ,
\end{equation} 
\begin{equation}
C_{\ell} \lesssim  I_{\ell}^{\gamma} \, .
\end{equation} 
Then, $\forall \varepsilon < e^{-1}$, there exists a value $L$ (and corresponding $I_{L})$ and a sequence $\{ N_{\ell}\}_{\ell=0}^L$ such that
\begin{equation}
	MSE[\big<  F_L\big>] = \mathbb{E}\Big[ \big(  \big<  F_L\big> - \mathbb{E}[F]\big)^2\Big]  < \varepsilon^2 \, ,
\end{equation}
and
\begin{equation}
	C\big(\big<  F_L\big> \big) \lesssim
	\begin{cases}
	 \varepsilon^{-2}  & \mbox{if} \quad \beta> \gamma  \, , \\
	 \varepsilon^{-2} (\log \varepsilon) ^2 & \mbox{if} \quad \beta = \gamma  \, ,\\	 
	 	 \varepsilon^{-2-(\frac{\gamma - \beta}{\alpha})}  & \mbox{if} \quad \beta <\gamma \, .\\
	\end{cases}
\end{equation}
\end{theorem}
More general theorems related to MLMC can be found elsewhere  \citep{giles-OP-2008,giles-AN-2015}.

\section{The MLHT algorithms with MLMC Optimization \label{sec:mlmc-mlht}}

Algorithm  \ref{alg:hybrid_mlmc} presents the MLHT algorithms with MLMC optimization for computing a functional of transport solution  $F[\phi]$ with optimization of computational costs for the given error $\varepsilon$.
The algorithm consists of three stages for estimating the numbers of realizations $\{N_{\ell}\}_{\ell=0}^L$.
$s$ is the index of the algorithm stage. 
At the initial stage ($s=1$), the algorithm starts calculations with some initial given number of realizations $\bm{N}^{ini}=\{ N_{\ell}^{ini} \}_{\ell=0}^L$ prescribed for each level to evaluate variances $V_\ell = \mathbb{V}\big[\big< \Delta F_{\ell}\big> \big]$ for $\ell=0,\ldots, L$ and associated costs $C_\ell$.
These data provide the basis for the estimation of an extra number of realizations  $\tilde N_{\ell}$
Eq. \eqref{estimation-N_ell} needed to meet the tolerance $\varepsilon$ according to MLMC theory \citep{giles-AN-2015}.
At the second stage ($s=2$), the algorithm performs computations with an estimated number of realizations to obtain the functional with the given accuracy $\varepsilon$. 
The numerical solution at this stage is also used to improve estimations of $V_\ell = \mathbb{V}\big[\big< \Delta F_{\ell}\big> \big]$,  $C_\ell$ for $\ell=0,\ldots, L$ and run more realizations on the last stage ({$s=3$) if required.
The number of levels is chosen as an input parameter.
At the end of the second stage, the algorithm  checks  the weak convergence criterion
\begin{equation} \label{weak-conv}
W_{\hat{\ell}} = \frac{\big< \Delta F_{\hat{\ell}}\big>}{2^{\alpha} - 1} < \frac{\varepsilon}{\sqrt{2}}, \quad \mbox{where} \quad \hat{\ell} = \{L-2,L-1,L\}
\end{equation}
to ensure that the number of levels $L$ chosen as an input parameter is sufficient to converge to the desired mean squared error.
Hereafter, we refer to MLHQD and MLHSM  algorithms with MLMC optimization as MLMC-HQD and MLMC-HSM, respectively.
\begin{algorithm}[ht!]
	\DontPrintSemicolon
	\caption{The MLHT algorithm with MLMC Optimization for Calculation of $F[\phi]$
    }\label{alg:hybrid_mlmc}
	$\tilde N_\ell = N_\ell^{ini}\, , \ N_{\ell} = \tilde N_\ell \, , \ \ell=0,\ldots,L$ \;
	\For{$s=1, \ldots, 3$}{
        \texttt{$\bullet$ Stage 1}\;
			\For{$n = 1, \ldots, \tilde N_0$}{
				execute MLHT  Algorithm  \ref{alg:mlht} to compute $n^{th}$ realization at $0^{th}$ level\;
			}
			compute   $\big< \phi_{0}\big>^{(s)}$ using  $\tilde N_0$ realizations \;
			compute   $\big< F_{0}\big>^{(s)}$\;
			\eIf{$s =1$}{
				$\big< \phi_{0}\big>$ $\leftarrow$ $\big< \phi_{0}\big>^{(s)}$,
				 $\big< F_{0}\big>$ $\leftarrow$ $\big< F_{0}\big>^{(s)}$\;
			}{
				 $\big< \phi_{0}\big>$ $\leftarrow$ combine  $\big< \phi_{0}\big>^{(s)}$ and $\big< \phi_{0}\big>$\;
				$\big< F_{0}\big>$ $\leftarrow$ combine  $\big< F_{0}\big>^{(s)}$ and $\big< F_{0}\big>$\;
			}
			evaluate $V_0 = \mathbb{V}\big[\big< F_{0}\big> \big]$ and $C_0$ based on $N_0$ realizations \;
        \texttt{$\bullet$ Stage 2}\;
		\For{$\ell = 1, \ldots, L$}{
			\For{$n = 1, \ldots, \tilde N_{\ell}$}{
			execute MLHT  Algorithm  \ref{alg:mlht}  to compute $n^{th}$ realization at $\ell^{th}$ level\;
	      	}
			compute $\big<\Delta \phi_{\ell} \big>^{(s)}$ using  $\tilde N_\ell$ realizations \;	
			 compute $\big<\Delta F_{\ell} \big>^{(s)}$  \;
				\eIf{$s =1$}{
					$\big< \Delta \phi_{\ell} \big>$ $\leftarrow$ $\big< \Delta \phi_{\ell} \big>^{(s)}$, 
				$\big< \Delta F_{\ell} \big>$ $\leftarrow$ $\big< \Delta F_{\ell} \big>^{(s)}$\;
			}{
            	 $\big< \Delta \phi_{\ell}\big>$  $\leftarrow$ combine  $\big<\Delta  \phi_{\ell} \big>^{(s)}$ and $\big< \Delta  \phi_{\ell}\big>$ \;
				$\big< \Delta F_{\ell}\big>$  $\leftarrow$ combine  $\big<\Delta F_{\ell} \big>^{(s)}$ and $\big< \Delta F_{\ell}\big>$\;
			}
			evaluate $V_\ell = \mathbb{V}\big[\big< \Delta F_{\ell}\big> \big]$ and $C_\ell$  based on $N_{\ell}$ realizations\;
		}
		   compute	$\big< \boldsymbol{\phi}_L\big> =  \mathcal{I}_{0}^L \big< \boldsymbol{\phi}_0 \big>  +  \sum_{\ell=1}^L  \mathcal{I}_{\ell}^L \big< \Delta \boldsymbol{\phi}_\ell \big> $ \;
			compute  $\big< F_{L} \big>=  \big< F_0 \big>  +  \sum_{\ell=1}^L  \big< \Delta F_\ell \big> $\;
			$\Delta N_\ell  \gets \max\big( \big \lfloor  2\varepsilon^{-2} \sqrt{V_\ell / C_\ell} \sum_{\ell=0}^{L} \sqrt{V_\ell C_\ell} \big\rfloor  -   N_\ell, 0 \big)\, ,  \ell=0,\ldots,L$\;
			\eIf{$\Delta N_{\ell}  = 0 \ \forall \ell$}{
			stop}{ 
            $\tilde N_{\ell}  \leftarrow  \Delta N_{\ell}$\;
             $N_{\ell} \leftarrow N_{\ell}  +   \Delta N_{\ell}$\;} 
		}
\end{algorithm}

In this study, we consider the functionals $F[\phi]$  defined  as  an integral of the scalar flux  over a spatial region $A$:
\begin{equation} \label{F-A}
    \mathcal{F}_A \stackrel{\Delta}{=} \int_{A} \phi(x) dx \, ,
\end{equation}
where $A$ is either the whole spatial domain ($A=D$) or a cell  $\tau_{i,0}$ on the coarsest grid $G_0$.
We  also perform optimization across the whole problem domain by considering a vector of functionals  $\{\mathcal{F}_{\tau_{i,0} }\}_{i=1}^{I_0}$ defined on the set  of cells of the grid $G_0$
\begin{equation} \label{F-tau-i-0}
	\mathcal{F}_{\tau_{i,0} } \stackrel{\Delta}{=}  \int_{\tau_{i,0}} \phi(x)dx \, .
\end{equation}
In this case, we calculate the variances $\{ V_{i, \ell} \}_{i=1}^{I_0}$ of  $\{ \big< \Delta \mathcal{F}_{\tau_{i,0}, \ell} \big> \}_{i=1}^{I_0}$ and apply $V_{i,\ell}$ to calculate the optimal number of realizations for all cell, $N_{i,\ell}$.
The algorithm uses $N_{\ell} = \max_i N_{i,\ell}$.
Another option is to  use $\max_i V_{i,\ell}$ to define $N_\ell$.
In either case, the cell with the highest variance estimate determines the algorithm's parameters.  
The second version requires the same or more realizations to converge since
\begin{equation}
\max_{i}\bigg(\sqrt{\frac{V_{i,\ell}}{C_{\ell}}} \sum_{\ell=0}^{L} \sqrt{V_{i,\ell} C_\ell}\bigg) \le 
	\sqrt{\frac{ \max_{i} V_{i,\ell} }{ C_{\ell} } }\sum_{\ell=0}^{L} \sqrt{C_{\ell}\max_{i}V_{i,\ell} }  \, .
\end{equation}
The benefit of using the second version is its simplicity in computing the optimal $N_\ell$, since only the highest variance needs to be stored at each level.
However, the first method may require fewer MLMC samples, thereby reducing overall computational cost.

\FloatBarrier
\section{Numerical Results \label{sec:res}} 

We consider a group of 1D problems for a slab $x\in [0, 1]$ cm.
($X = 1$ cm) with the constant  total cross section $\Sigma_t(x) = 1$ cm$^{-1}$,  external source $q(x)=1, \ x\in [0,1]$ cm, and vacuum BCs $(\psi_{in}^{\pm}=0$).
The tests differ by the number of material regions and their parameters.
\begin{itemize}
	\item {\bf Test 1.} It is a one-region problem with the scattering ratio $c=0.9$ ($c=\frac{\Sigma_s(x)}{\Sigma_t(x)}$).
	\item {\bf Test 2.}   This is a two-region slab with  subdomains  given by
	\begin{itemize}
		\item Region 1:   $x\in  [0,0.5]$ cm,  $c_1=0.9$,
		\item Region 2:   $x\in  [0.5,1]$ cm,  $c_2= \{0.1,0.5\}$.
	\end{itemize}	
 \end{itemize}
MC calculations are performed using the implicit-capture method and Russian roulette to improve the efficiency of the Monte Carlo simulation.
Implicit capture was used to demonstrate that traditional variance-reduction techniques can also be applied to HMCD methods. 
In general, variance reduction might be necessary to encourage particles to visit all phase-space elements during the simulation.
We use a minimum weight of $10^{-4}$ for this set of results.

The results were generated by a code written in Julia, which uses the $\texttt{Xoshiro256++}$ PRNG sequence to provide the random numbers used in the MC simulations \citep{Julia-2017, Blackman_Vigna_2021}.
The MC simulation is parallelized using the $\texttt{Threads.@threads}$ macro.
Each thread has its own tally bin to collect results in, and then the results across threads are summed together into a flattened array at the end of the simulation, similar to an $\texttt{MPI\_SUM}$ operation.
Julia's PRNG uses a per-Task state, allowing reproducible results in multi-threaded simulations, provided the same version of Julia is used for repeated simulations \citep{Julia-2017}.
We verify this by comparing the MC results for the different simulations, which should match exactly.

\subsection{Comparison of HQD and HSM Schemes}
To analyze the accuracy of the HQD and HSM schemes (i.e., single-level hybrid transport methods), we use Test 1 and evaluate the relative error of numerical solutions in the $L_2$ norm given by
\vspace{-0.1cm}
\begin{equation} \label{eqn:rel_L2_Norm_Error}
 	\text{RE}_{L_2}(\phi) = \frac{|| \phi - \phi^{ex}||_{L_2}}{||\phi^{ex}||_{L_2}} =\sqrt{ \frac{\sum_{i=1}^I (\phi_i - \phi_i^{ex})^2\Delta x_i}
 		{\sum_{i=1}^I (\phi_i^{ex})^2 \Delta x_i} } \, .
 \end{equation}
The reference numerical solution $\phi^{ex}$  is computed by means of a deterministic transport method on a sequence of refined phase-space grids and Aitken extrapolation \citep{Sidi_2003,dahlquist2008numerical, Ganapol2013}. 
Appendix \ref{sec:appendix_a} provides details on the generation and the numerically exact values of the reference numerical solution $\phi^{ex}$.
 
The test is solved on  uniform spatial grids with $\Delta x = 2^{-m}$ cm, $m=2,\ldots,6$ and different numbers of particle histories $K = \{10^3, 10^4, 10^5 \}$.
The results demonstrate the effects  of a decrease in  discretization error in the approximation of low-order equations
with refinement of spatial grids and reduction in statistical error as the number of histories $K$ increases.
 
Table \ref{tab:mc_hybrid_l2_errors} presents the mean relative error $\big< RE_{L_2}(\phi) \big>$ and the standard deviation of the mean relative error $\sigma_{\big< RE_{L_2}\big>}$ calculated for numerical solutions of 100 simulations.
100 simulations were chosen to ensure an approximately normal distribution of $RE_{L_2}$ errors, which was verified by examining the histogram of error norms, quantile-quantile (Q-Q) plots, and performing statistical tests for normality, i.e., Shapiro-Wilk, Anderson-Darling, and Kolmogorov-Smirnov. 
Results indicated that the error norms are normally distributed, validating the use of the mean and standard deviation in this analysis.
The HQD and HSM  use different collections of particle histories in these calculations.
The results show that the HSM method has a lower relative $L_2$ error than the HQD method for $\Delta x \leq 2^{-3}$ cm.
We note that for $K=10^5$ and $\Delta x = 2^{-4}$ cm, the HQD solution has a statistically significantly lower error than the HSM solution.
There is no statistically significant difference in mean error for (a) $\Delta x\le 2^{-3}$ cm and $K=10^3, 10^4$ and (b)  $\Delta x\le 2^{-5}$ cm and $K=10^5$.
No statistically significant difference means one cannot distinguish between the two methods for more-refined grids.
This indicates that the effects of discretization error are small relative to statistical noise on refined grids.
  \begin{table}[htb]
  \centering
  \caption{Test 1. Mean relative error in  $L_2$-norm    $\big< RE_{L_2}(\phi) \big>$ and $\sigma_{\big< RE_{L_2}\big>}$ of  
  	HQD and HSM solutions  based on 100 simulations
  	 versus $\Delta x$ and $K$}
    \begin{tabular}{|c|c|c|}
        \hline
        $K = 10^3$         & HQD                                            & HSM                                              \\
        \hline
        $\Delta x =2^{-2}$ & $4.22 \times 10^{-2} \pm 1.57 \times 10^{-3} $ & $3.65 \times  10^{-2} \pm 1.44 \times  10^{-3}$  \\ 
        $\Delta x =2^{-3}$ & $4.31 \times 10^{-2} \pm 1.22 \times 10^{-3} $ & $4.65 \times  10^{-2} \pm 1.37 \times  10^{-3}$  \\
        $\Delta x =2^{-4}$ & $5.18 \times 10^{-2} \pm 1.19 \times 10^{-3} $ & $5.42 \times  10^{-2} \pm 1.01 \times  10^{-3}$  \\
        $\Delta x =2^{-5}$ & $6.08 \times 10^{-2} \pm 1.13 \times 10^{-3} $ & $6.05 \times  10^{-2} \pm 1.22 \times  10^{-3}$  \\
        $\Delta x =2^{-6}$ & $6.48 \times 10^{-2} \pm 1.13 \times 10^{-3} $ & $6.58 \times  10^{-2} \pm 1.11 \times  10^{-3}$  \\
        \hline  \hline
        $K = 10^4$         & HQD                                            & HSM                                              \\
        \hline
        $\Delta x =2^{-2}$ & $2.44 \times 10^{-2} \pm 6.29 \times 10^{-4} $ & $2.07 \times  10^{-2} \pm 6.38 \times  10^{-4}$  \\
        $\Delta x =2^{-3}$ & $1.45 \times 10^{-2} \pm 4.22 \times 10^{-4} $ & $1.48 \times  10^{-2} \pm 3.94 \times  10^{-4}$  \\
        $\Delta x =2^{-4}$ & $1.63 \times 10^{-2} \pm 3.38 \times 10^{-4} $ & $1.67 \times  10^{-2} \pm 4.09 \times  10^{-4}$  \\ 
        $\Delta x =2^{-5}$ & $1.89 \times 10^{-2} \pm 3.76 \times 10^{-4} $ & $1.93 \times  10^{-2} \pm 3.57 \times  10^{-4}$  \\ 
        $\Delta x =2^{-6}$ & $2.08 \times 10^{-2} \pm 3.00 \times 10^{-4} $ & $2.07 \times  10^{-2} \pm 3.13 \times  10^{-4}$  \\ 
        \hline  \hline
        $K = 10^5$         & HQD                                            & HSM                                              \\
        \hline
        $\Delta x =2^{-2}$ & $2.32 \times 10^{-2} \pm 2.21 \times 10^{-4} $ & $1.81 \times  10^{-2} \pm 2.23 \times  10^{-4}$  \\
        $\Delta x =2^{-3}$ & $7.24 \times 10^{-3} \pm 1.79 \times 10^{-4} $ & $6.20 \times  10^{-3} \pm 1.70 \times  10^{-4}$  \\
        $\Delta x =2^{-4}$ & $4.41 \times 10^{-3} \pm 1.15 \times 10^{-4} $ & $5.18 \times  10^{-3} \pm 1.20 \times  10^{-4}$  \\
        $\Delta x =2^{-5}$ & $6.05 \times 10^{-3} \pm 9.58 \times 10^{-5} $ & $5.93 \times  10^{-3} \pm 1.02 \times  10^{-4}$  \\
        $\Delta x =2^{-6}$ & $6.65 \times 10^{-3} \pm 9.28 \times 10^{-5} $ & $6.63 \times  10^{-3} \pm 1.05 \times  10^{-4}$  \\
        \hline
    \end{tabular}
  \label{tab:mc_hybrid_l2_errors}
\end{table}
\begin{table}[ht]
	\centering
	\caption{Test 1. Relative $L_2$ norm of $\big< \phi_{\ell}\big>$ (Eq. \eqref{phi_ell}) computed by the MLHQD algorithm}
	\begin{tabular}{|c|c||c|c|c|c|}
		\hline
		$\ell$&$N_{\ell}$ &$K_{\ell}=10^2$ &$K_{\ell}= 10^3$ &$K_{\ell}=10^4$ & $K_{\ell}=10^5$\\
		\hline
		0 &100 &$2.33 \times 10^{-2}$  & $2.08 \times 10^{-2}$   & $2.02 \times 10^{-2}$ & $2.01 \times 10^{-2}$  \\
		1 &50  &$1.91 \times 10^{-2}$  & $1.16 \times 10^{-2}$   & $1.01 \times 10^{-2}$ & $9.99 \times 10^{-3}$  \\
		2 &25  &$2.21 \times 10^{-2}$  & $8.94 \times 10^{-3}$   & $5.02 \times 10^{-3}$ & $4.55 \times 10^{-3}$  \\
		3 &10  &$3.49 \times 10^{-2}$  & $1.18 \times 10^{-2}$   & $3.78 \times 10^{-3}$ & $1.07 \times 10^{-3}$  \\
		\hline
	\end{tabular}
	\label{tab:test1-mlhqd-phi-ell}
\end{table}
\begin{table}[ht]
	\centering
	\caption{Test 1. Relative $L_2$ norm of $\big< \phi_{\ell}\big>$ (Eq. \eqref{phi_ell}) computed by the MLHSM algorithm}
	\begin{tabular}{|c|c||c|c|c|c|}
		\hline
		$\ell$&$N_{\ell}$ &$K_{\ell}=10^2$ &$K_{\ell}= 10^3$ &$K_{\ell}=10^4$ & $K_{\ell}=10^5$\\
		\hline
		0 &100 &$2.36 \times 10^{-2}$  & $2.08 \times 10^{-2}$   & $2.02 \times 10^{-2}$ & $2.01 \times 10^{-2}$  \\
		1 &50  &$1.94 \times 10^{-2}$  & $1.17 \times 10^{-2}$   & $1.01 \times 10^{-2}$ & $9.98 \times 10^{-3}$  \\
		2 &25  &$2.17 \times 10^{-2}$  & $9.06 \times 10^{-3}$   & $5.03 \times 10^{-3}$ & $4.55 \times 10^{-3}$  \\
		3 &10  &$3.45 \times 10^{-2}$  & $1.19 \times 10^{-2}$   & $3.79 \times 10^{-3}$ & $1.08 \times 10^{-3}$  \\
		\hline
	\end{tabular}
	\label{tab:test1-mlhsm-phi-ell}
\end{table}

\subsection{Global Numerical Solution  of MLHT Algorithms}
 
In this section, we analyze the performance of the multilevel algorithms in computing numerical solutions of the transport equation over the whole domain on the target grid $G_L$.
Test 1 is  solved on the grid with $\Delta x=2^{-7}$ cm by MLHT methods with $L=3$ using  a sequence uniform grids $G_\ell$ with  $I_\ell = 2 I_{\ell -1}$, where $I_0 = 16$. 
The MLHT methods (Algorithm \ref{alg:mlht}) use a prescribed set of realizations at each level $\bm{N} = \{ N_{\ell}\}_{\ell=0}^L$, namely, $\bm{N} =\{100, 50, 25, 10\}$.
The same collections of particle histories were used by both MLHT algorithms.
Figures \ref{fig:test1-mlhqd} and \ref{fig:test1-mlhsm} show plots of $\mathcal{I}_{0}^L\big< \phi_0 \big>$ and  $\mathcal{I}_{\ell}^L\big< \Delta \phi_\ell  \big>$ on each $G_\ell$ for Test 1  in the case of $K_{\ell}=10^4$, $\ell=0,\ldots, L$ presenting elements of global hybrid solutions prolongated on the target grid  $G_L$.
Figures \ref{fig:test2-01-mlhqd}  and \ref{fig:test2-01-mlhsm} show  $\mathcal{I}_{0}^L\big< \phi_0 \big>$ and  
$\mathcal{I}_{\ell}^L\big< \Delta \phi_\ell  \big>$ for the two-region Test 2 with $c_2=0.1$  solved by the MLHT algorithm with the same parameters as Test 1.
To illustrate intermediate components of the numerical solution of MLHT algorithms while moving through the grid levels,  we compute the solution obtained by reaching each of the levels and given by
\vspace{-0.1cm}
\begin{equation} \label{phi_ell}
    \big< \boldsymbol{\phi}_{\ell}\big> =  \mathcal{I}_{0}^L \big< \Delta \boldsymbol{\phi}_0 \big> +  \sum_{\ell^{\prime}=1}^{\ell}\mathcal{I}_{\ell^{\prime}}^L \big< \Delta \boldsymbol{\phi}_{\ell^{\prime}} \big> \, , \quad \ell=0,\ldots, L \, .
    \vspace{-0.1cm}
\end{equation}
Tables \ref{tab:test1-mlhqd-phi-ell} and \ref{tab:test1-mlhsm-phi-ell} present the relative error in $\big< \boldsymbol{\phi}_{\ell} \big>$ for the algorithms using different $K_{\ell}$.

\pagebreak
\begin{figure}[h]
    \begin{subfigure}{0.49\textwidth}
        \centering
        \includegraphics[width=\textwidth]{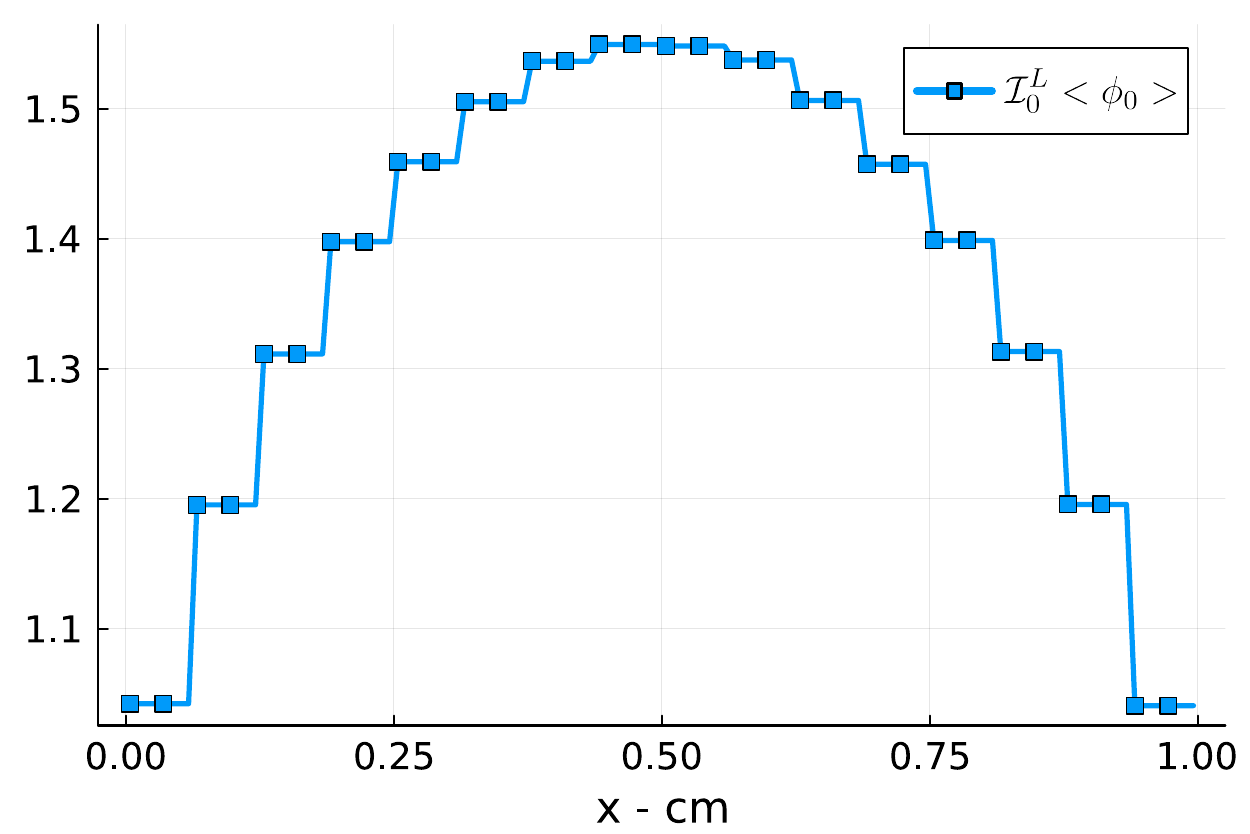}
        \caption{$\mathcal{I}_{0}^L\big<\phi_0 \big>$  \label{fig:test1-phi0-mlhqd}}
    \end{subfigure}
    \begin{subfigure}{0.49\textwidth}
		\centering
		\includegraphics[width=\textwidth]{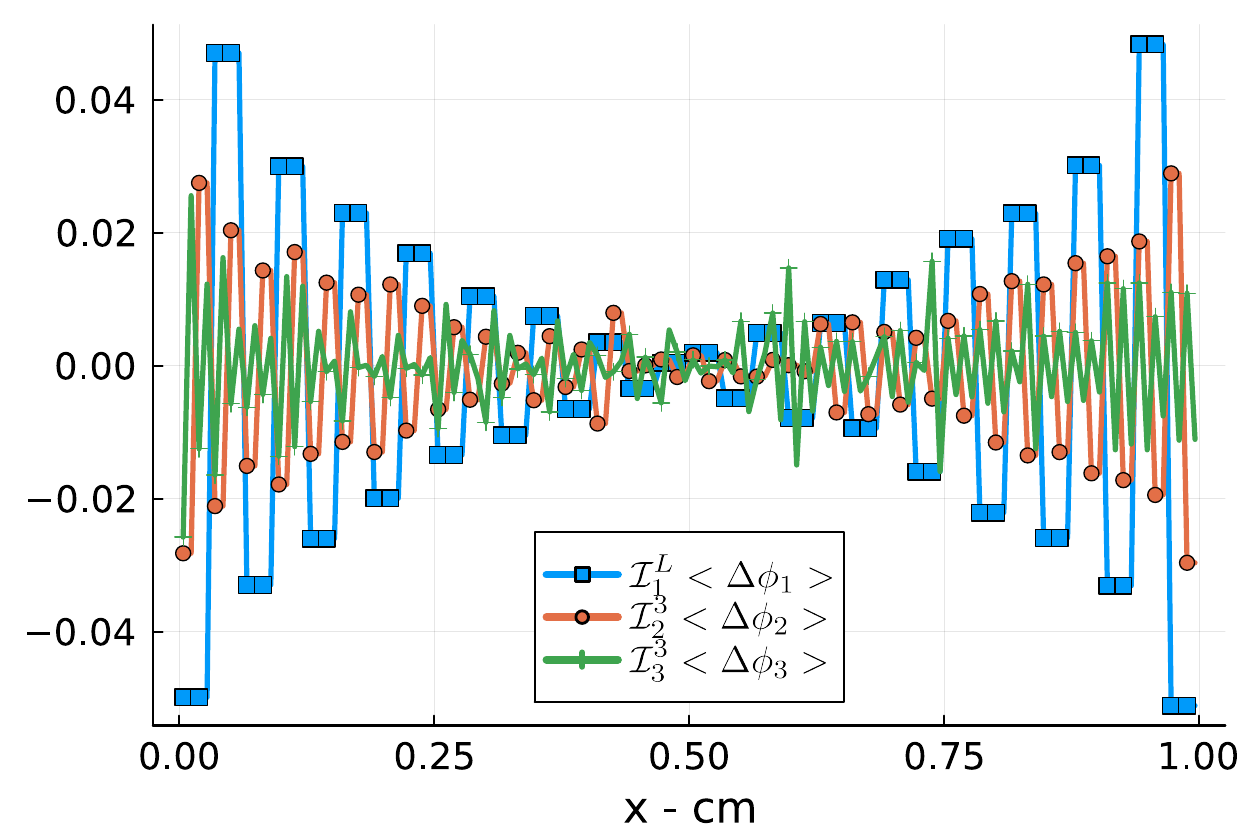}
		\caption{$\mathcal{I}_{\ell}^L\big<\Delta \phi_{\ell} \big>$ \label{fig:test1-delta-phi-mlhqd}}
	\end{subfigure}
	\caption{Test 1. $\mathcal{I}_{0}^L\big<\phi_0 \big>$ and $\mathcal{I}_{\ell}^L\big<  \Delta \phi_{\ell} \big>$ obtained by   the MLHQD algorithm with $L=3$ and $K_{\ell}=10^4$, $\ell=0,\ldots,L$. \label{fig:test1-mlhqd}}
\end{figure}

\begin{figure}[h!]
    \begin{subfigure}{0.49\textwidth}
		\centering
		\includegraphics[width=\textwidth]{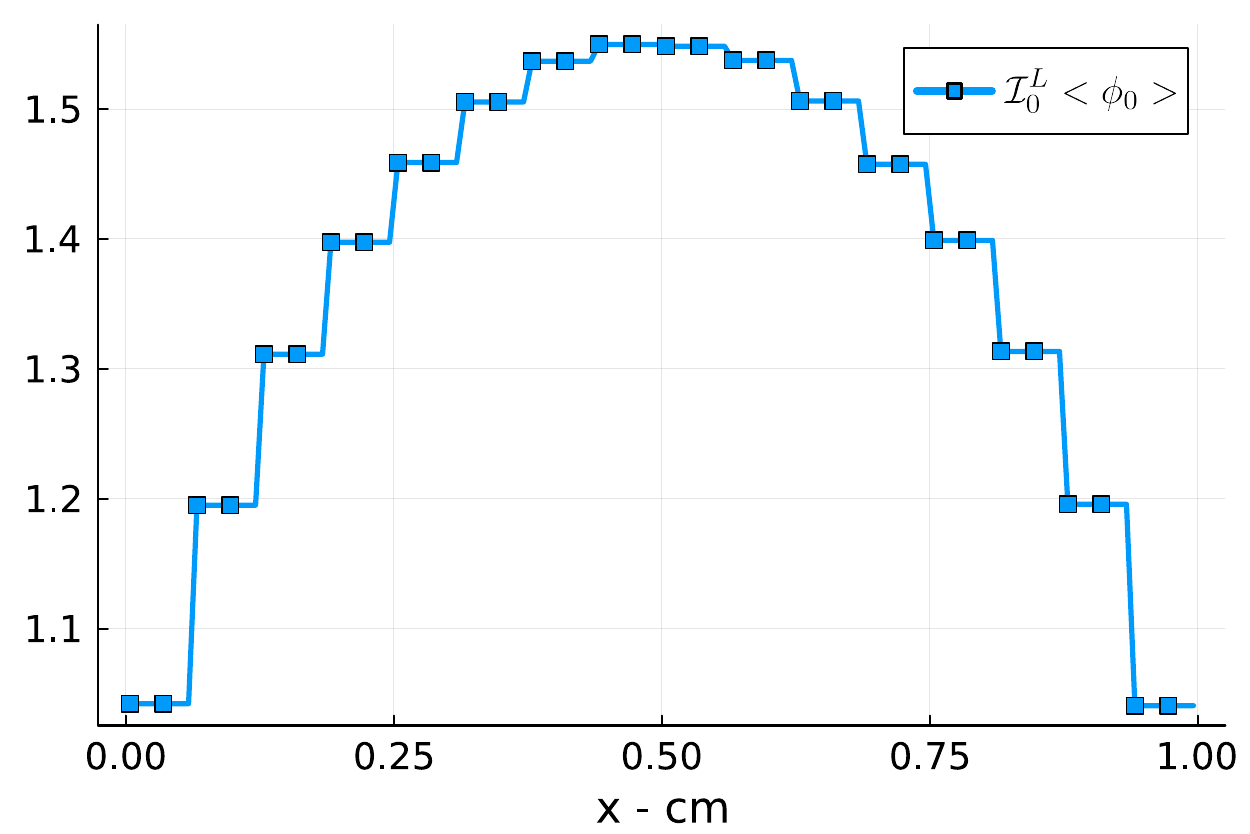}
		\caption{$\mathcal{I}_{0}^L\big<\phi_0 \big>$ \label{fig:test1-phi0-mlhsm}}
	\end{subfigure}
	\begin{subfigure}{0.49\textwidth}
		\centering
		\includegraphics[width=\textwidth]{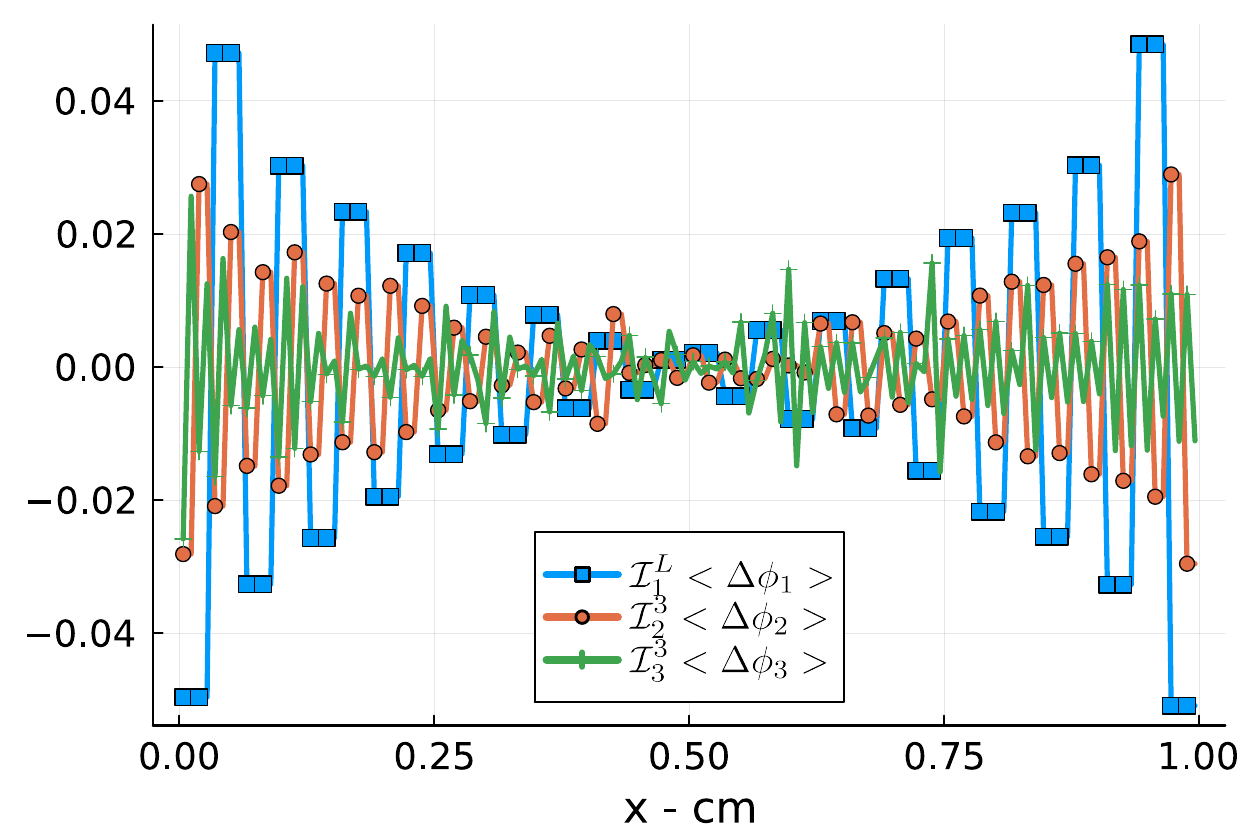}
		\caption{$\mathcal{I}_{\ell}^L\big<\Delta \phi_{\ell} \big>$ \label{fig:test1-delta-phi-mlhsm}}
	\end{subfigure}
	\caption{Test 1. $\mathcal{I}_{0}^L\big< \phi_0 \big>$ and $\mathcal{I}_{\ell}^L\big<  \Delta \phi_{\ell} \big>$ obtained by the MLHSM algorithm with $L=3$ and $K_{\ell}=10^4$, $\ell=0,\ldots,L$. \label{fig:test1-mlhsm}}
\end{figure}

\begin{figure}[h!]
	\begin{subfigure}{0.49\textwidth}
		\centering
		\includegraphics[width=\textwidth]{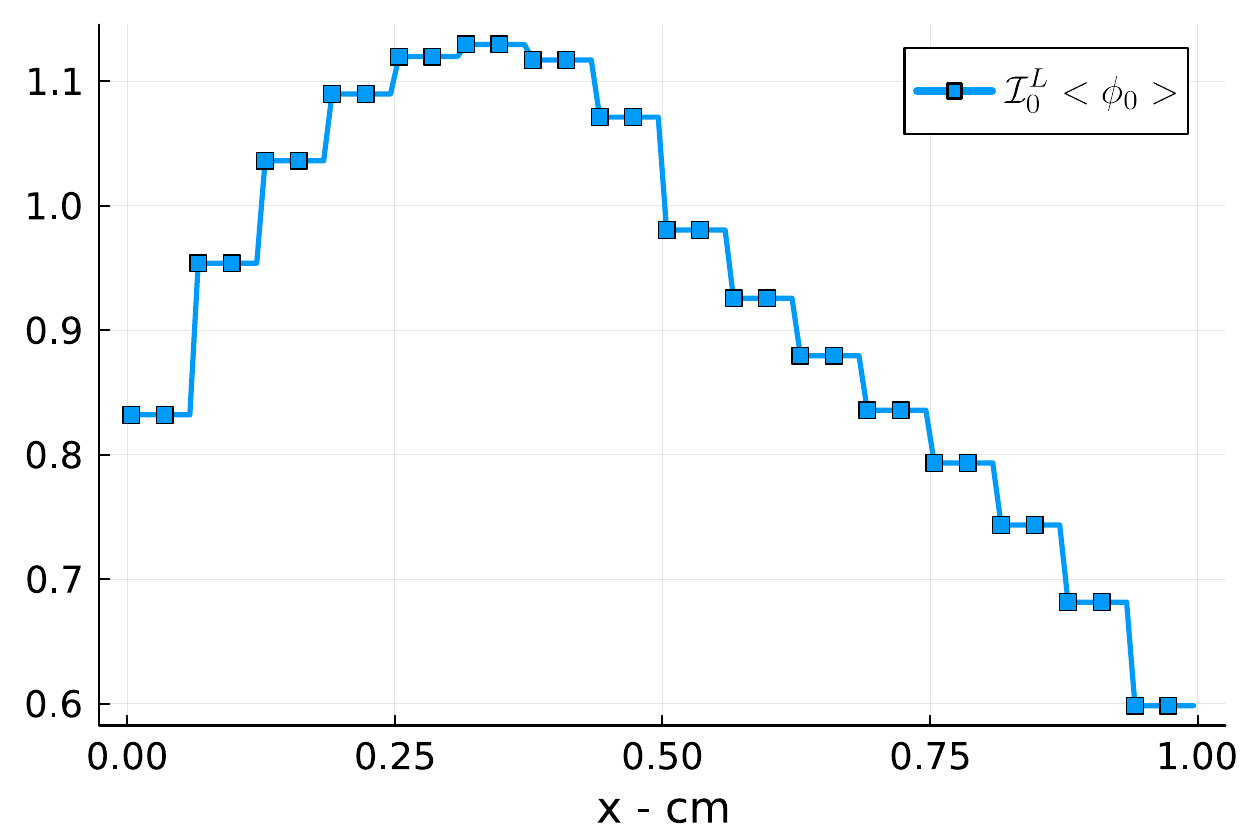}
		\caption{$\mathcal{I}_{0}^L\big<\phi_0 \big>$ \label{fig:test2-01-phi0-mlhqd}}
	\end{subfigure}
	\begin{subfigure}{0.49\textwidth}
		\centering
		\includegraphics[width=\textwidth]{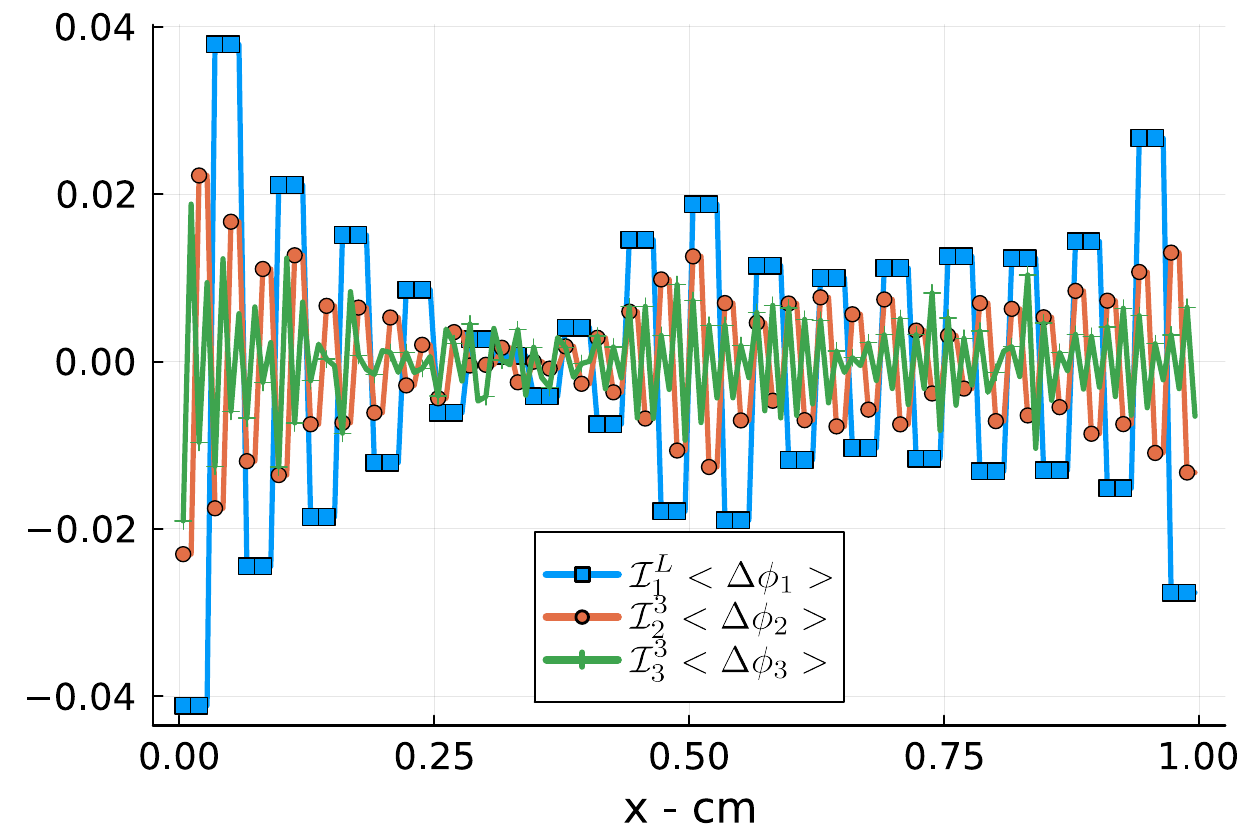}
		\caption{$\mathcal{I}_{\ell}^L\big<\Delta \phi_{\ell} \big>$ \label{fig:test2-01-delta-phi-mlhqd}}
	\end{subfigure}
	\caption{Test 2 with $c_2=0.1$. $\mathcal{I}_{0}^L\big<  \phi_0 \big>$ and $\mathcal{I}_{\ell}^L\big<  \Delta \phi_{\ell} \big>$  obtained by the MLHQD algorithm with $L=3$ and $K_{\ell}=10^4$, $\ell=0,\ldots,L$. \label{fig:test2-01-mlhqd} }
\end{figure}
 \FloatBarrier
\begin{figure}[h!]
	\begin{subfigure}{0.49\textwidth}
		\centering
		\includegraphics[width=\textwidth]{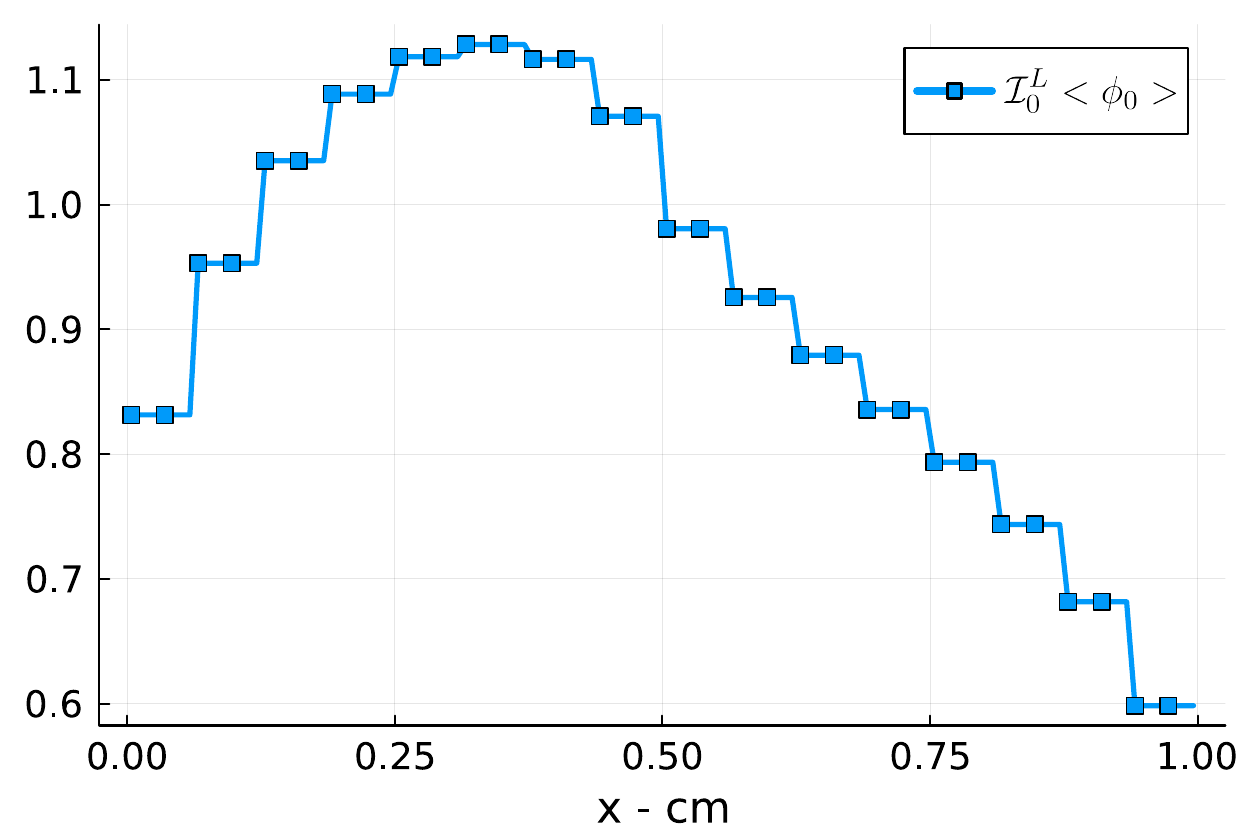}
		\caption{$\mathcal{I}_{0}^L\big<\phi_0 \big>$\label{fig:test2-01-phi0-mlhsm}}
	\end{subfigure}
	\begin{subfigure}{0.49\textwidth}
		\centering
		\includegraphics[width=\textwidth]{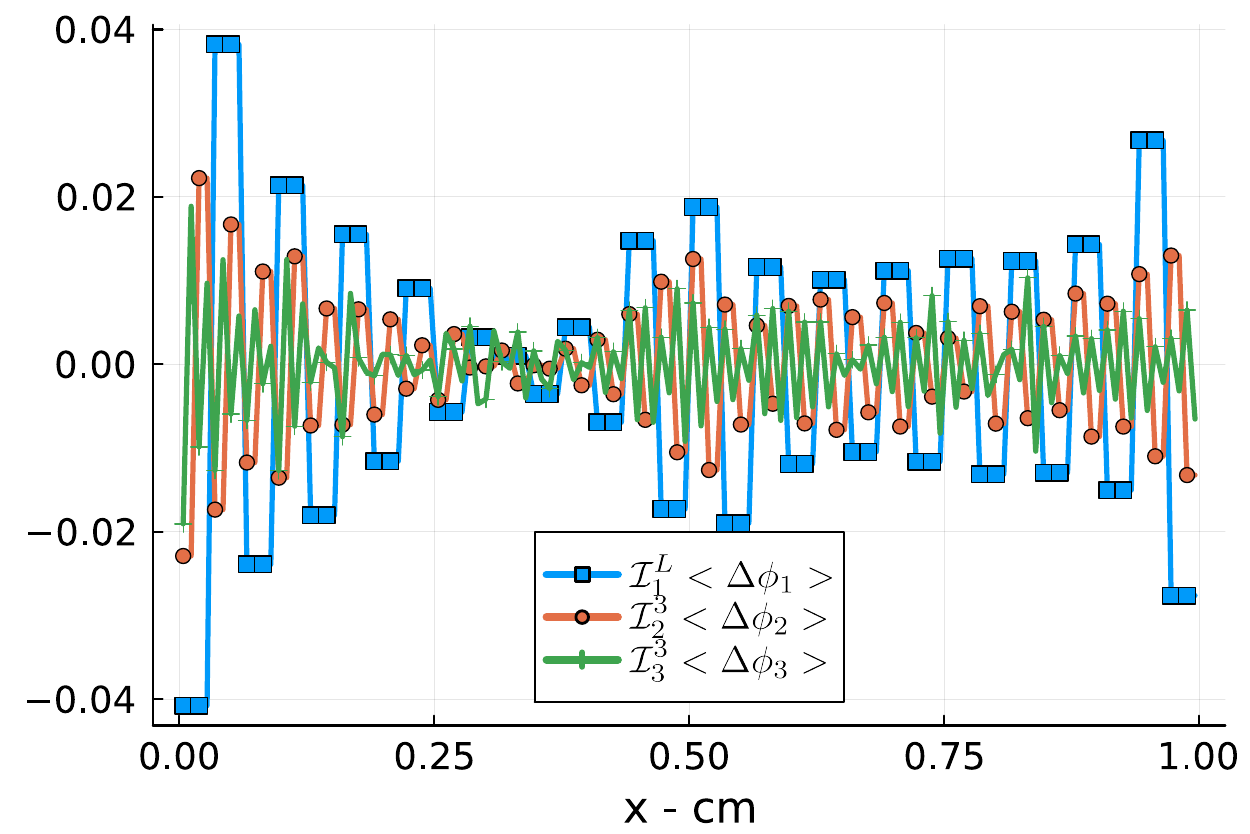}
		\caption{$\mathcal{I}_{\ell}^L\big<\Delta \phi_{\ell} \big>$ \label{fig:test2-01-delta-phi-mlhsm}}
	\end{subfigure}
	\caption{Test 2 with $c_2=0.1$. $\mathcal{I}_{0}^L\big<  \phi_0 \big>$ and $\mathcal{I}_{\ell}^L\big<  \Delta \phi_{\ell} \big>$ obtained by the MLHSM algorithm with $L=3$ and $K_{\ell}=10^4$, $\ell=0,\ldots,L$. \label{fig:test2-01-mlhsm}} 
\end{figure}

\subsection{Convergence of MLHT algorithms with MLMC Optimization for Functionals \label{sec:MLHT-convergence} }

This section presents numerical results of  MLMC-HQD and MLMC-HSM algorithms optimizing calculations of some given functional $F[\phi]$ (Algorithm \ref{alg:hybrid_mlmc}).
We apply these algorithms with $L=3$ to solve Test 2 on the grid having $\Delta x=2^{-7}$  and use uniform grids $G_\ell$ with  $I_\ell = 2 I_{\ell -1}$, where $I_0 = 16$. 
The initial stage of  MLMC algorithm ($s=1$) is executed with $\bm{N}^{ini} = \{ 10,10,10,10\}$.
The calculations are performed for different values of error $\varepsilon$ and the scattering ratio in Region 2, $c_2$. 
We note that the case of $c_2=0.9$ corresponds to a one-region problem of Test 1.
    
First, we solve this test to compute $ F = \mathcal{F}_D$ (Eq. \eqref{F-A}).     
Tables \ref{tab:mlmc_mixed_1e3} and \ref{tab:mlmc_mixed_1e3_sm} present $\alpha$, $\beta$, and $\gamma$, which characterize convergence of the functional and its variance, and change in computational costs of the algorithms, respectively (see Theorem \ref{theorem}).
These values are estimated by performing a linear fit to the $\log_2$ of the mean, variance, and computational costs, respectively.
The tables also show the total number of realizations $N_{\ell}$ performed by the MLMC algorithm, along with data on the evaluation of weak convergence (Eq. \eqref{weak-conv}).
The results of calculations with more particle histories, namely, $K=10^4$, are presented in Tables \ref{tab:mlmc_mixed_1e4} and \ref{tab:mlmc_mixed_1e4_sm}.
      
The data shows $\alpha>0$, $\beta>0$, $\gamma>0$, $\alpha> \frac{1}{2} \min(\beta,\gamma)$, and hence performance of the algorithms meets the condition of Theorem \ref{theorem}.
The convergence rate $\alpha \approx 2$. 
This is expected due to the second-order accuracy of spatial discretization schemes for the low-order equations.  
We notice that $\beta > \gamma$. 
Thus, the variance decreases faster than the cost of calculations increases.
The results also indicate that as $\varepsilon$ decreases, the required number of samples requested increases, since the variance convergence criteria are more stringent.
Increasing the number of particle histories ($K_{\ell}$) per realization also decreased the number of requested samples due to lower variance in the scalar flux sample and the functional.
Most of the computational work is placed on the coarsest level.
In addition, we note that the weak convergence criteria for the number of levels is met, meaning the MSE is bounded by $\varepsilon^2$.

       \begin{table}[htb]
      	\centering
      	\caption{Test 2. $F=\mathcal{F}_D$,  MLMC-HQD,   $K_{\ell} = 10^3$}
      	\begin{tabular}{|c|c||c|c|c|c|c|c|c|c|}
      		\hline
      		$c_2$ & $\varepsilon$   & $\alpha$ & $\beta$ & $\gamma$ & $N_0$ & $N_1$ & $N_2$ & $N_3$ & $\max_{\hat \ell}{ W_{\hat{\ell}}}$\\
      		\hline
      		  & $1 \times 10^{-2}$ & 2.10  & 2.50 & 0.61 & 10   & 10 & 10 & 10 & $3.1 \times 10^{-4}$ \\
      		0.1 & $5 \times 10^{-3}$ & 2.05  & 2.08 & 0.66 & 59   & 10 & 10 & 10 & $3.1 \times 10^{-4}$ \\
      		  & $1 \times 10^{-3}$ & 2.03  & 1.58 & 0.63 & 840  & 10 & 10 & 10 & $3.2 \times 10^{-4}$ \\
            \hline
      		  & $1 \times 10^{-2}$ & 1.99  & 2.06 & 0.61 & 12   & 10 & 10 & 10 & $3.6 \times 10^{-4}$ \\
      		0.5 & $5 \times 10^{-3}$ & 2.02  & 1.51 & 0.62 & 19   & 10 & 10 & 10 & $3.4 \times 10^{-4}$ \\
      		  & $1 \times 10^{-3}$ & 2.03  & 2.76 & 0.64 & 1182 & 10 & 10 & 10 & $3.4 \times 10^{-4}$ \\
            \hline
      		  & $1 \times 10^{-2}$ & 1.96  & 2.58 & 0.56 & 21   & 10 & 10 & 10 & $5.0 \times 10^{-4}$ \\
      		0.9 & $5 \times 10^{-3}$ & 2.01  & 1.61 & 0.62 & 57   & 10 & 10 & 10 & $4.9 \times 10^{-4}$ \\
      		  & $1 \times 10^{-3}$ & 2.01  & 1.85 & 0.61 & 3388 & 10 & 10 & 10 & $4.8 \times 10^{-4}$ \\
      		\hline
      	\end{tabular}
      	\label{tab:mlmc_mixed_1e3}
\vspace{0.5cm}
      	\centering
      	\caption{Test 2. $F=\mathcal{F}_D$,  MLMC-HSM,   $K_{\ell} = 10^3$}
      	\begin{tabular}{|c|c||c|c|c|c|c|c|c|c|}
      		\hline
      		$c_2$ & $\varepsilon$   & $\alpha$ & $\beta$ & $\gamma$ & $N_0$ & $N_1$ & $N_2$ & $N_3$ & $\max_{\hat \ell} {W_{\hat{\ell}}}$\\
      		\hline
      		  & $1 \times 10^{-2}$ & 1.98  & 3.75 & 0.64 & 10   & 10 & 10 & 10 & $2.5 \times 10^{-4}$ \\
      		0.1 & $5 \times 10^{-3}$ & 1.97  & 2.79 & 0.62 & 21   & 10 & 10 & 10 & $2.5 \times 10^{-4}$ \\
      		  & $1 \times 10^{-3}$ & 2.00  & 2.51 & 0.63 & 724  & 10 & 10 & 10 & $2.5 \times 10^{-4}$ \\
            \hline
      		  & $1 \times 10^{-2}$ & 2.01  & 3.45 & 0.65 & 10   & 10 & 10 & 10 & $2.8 \times 10^{-4}$ \\
      		0.5 & $5 \times 10^{-3}$ & 2.00  & 2.00 & 0.63 & 30   & 10 & 10 & 10 & $2.8 \times 10^{-4}$ \\
      		  & $1 \times 10^{-3}$ & 2.00  & 2.77 & 0.61 & 888  & 10 & 10 & 10 & $2.8 \times 10^{-4}$ \\
            \hline
      		  & $1 \times 10^{-2}$ & 2.00  & 3.54 & 0.61 & 12   & 10 & 10 & 10 & $3.8 \times 10^{-4}$ \\
      		0.9 & $5 \times 10^{-3}$ & 2.00  & 3.05 & 0.64 & 60   & 10 & 10 & 10 & $3.8 \times 10^{-4}$ \\
      		  & $1 \times 10^{-3}$ & 2.00  & 2.41 & 0.62 & 1536 & 10 & 10 & 10 & $3.8 \times 10^{-4}$ \\
      		\hline
      	\end{tabular}
      	\label{tab:mlmc_mixed_1e3_sm}
\vspace{0.5cm}
  	\centering
      	\caption{Test 2. $F=\mathcal{F}_D$,  MLMC-HQD, $K_{\ell} = 10^4$}
  	\begin{tabular}{|c|c||c|c|c|c|c|c|c|c|}
  		\hline
  		$c_2$ & $\varepsilon$   & $\alpha$ & $\beta$ & $\gamma$ & $N_0$ & $N_1$ & $N_2$ & $N_3$ & $\max_{\hat \ell} {W_{\hat{\ell}}}$ \\
  		\hline
  		  & $1 \times 10^{-2}$ & 2.01  & 2.92 & 0.66 & 10   & 10 & 10 & 10 & $3.2 \times 10^{-4}$ \\
  		0.1 & $5 \times 10^{-3}$ & 2.00  & 3.03 & 0.64 & 10   & 10 & 10 & 10 & $3.2 \times 10^{-4}$ \\
  		  & $1 \times 10^{-3}$ & 2.01  & 3.21 & 0.65 & 62   & 10 & 10 & 10 & $3.2 \times 10^{-4}$ \\
        \hline
  		  & $1 \times 10^{-2}$ & 1.99  & 2.91 & 0.68 & 10   & 10 & 10 & 10 & $3.6 \times 10^{-4}$ \\
  		0.5 & $5 \times 10^{-3}$ & 2.00  & 2.79 & 0.67 & 10   & 10 & 10 & 10 & $3.5 \times 10^{-4}$ \\
  		  & $1 \times 10^{-3}$ & 2.00  & 3.43 & 0.66 & 101  & 10 & 10 & 10 & $3.5 \times 10^{-4}$ \\
        \hline
  		  & $1 \times 10^{-2}$ & 2.00  & 3.57 & 0.69 & 10   & 10 & 10 & 10 & $4.9 \times 10^{-4}$ \\
  		0.9 & $5 \times 10^{-3}$ & 2.00  & 3.12 & 0.66 & 10   & 10 & 10 & 10 & $4.9 \times 10^{-4}$ \\
  		  & $1 \times 10^{-3}$ & 2.00  & 1.98 & 0.69 & 182  & 10 & 10 & 10 & $4.9 \times 10^{-4}$ \\
  		\hline
  	\end{tabular}
  	\label{tab:mlmc_mixed_1e4}
 \vspace{0.5cm}
  	\centering
      	\caption{Test 2. $F=\mathcal{F}_D$, MLMC-HSM, $K_{\ell} = 10^4$}
  	\begin{tabular}{|c|c||c|c|c|c|c|c|c|c|}
  		\hline
  		$c_2$ & $\varepsilon$   & $\alpha$ & $\beta$ & $\gamma$ & $N_0$ & $N_1$ & $N_2$ & $N_3$ & $\max_{\hat \ell} {W_{\hat{\ell}}}$ \\
  		\hline
  		    & $1 \times 10^{-2}$ & 2.00  & 1.95 & 0.69 & 10   & 10 & 10 & 10 & $2.5 \times 10^{-4}$ \\
  		0.1 & $5 \times 10^{-3}$ & 1.99  & 2.84 & 0.67 & 10   & 10 & 10 & 10 & $2.5 \times 10^{-4}$ \\
  		  & $1 \times 10^{-3}$ & 2.00  & 3.62 & 0.65 & 54   & 10 & 10 & 10 & $2.5 \times 10^{-4}$ \\
        \hline
  		  & $1 \times 10^{-2}$ & 2.00  & 3.04 & 0.68 & 10   & 10 & 10 & 10 & $2.8 \times 10^{-4}$ \\
  		0.5 & $5 \times 10^{-3}$ & 2.00  & 3.30 & 0.65 & 10   & 10 & 10 & 10 & $2.8 \times 10^{-4}$ \\
  		  & $1 \times 10^{-3}$ & 2.00  & 3.98 & 0.65 & 164  & 10 & 10 & 10 & $2.8 \times 10^{-4}$ \\
        \hline
  		  & $1 \times 10^{-2}$ & 1.99  & 3.34 & 0.70 & 10   & 10 & 10 & 10 & $3.8 \times 10^{-4}$ \\
  		0.9 & $5 \times 10^{-3}$ & 2.00  & 2.44 & 0.71 & 10   & 10 & 10 & 10 & $3.8 \times 10^{-4}$ \\
  		  & $1 \times 10^{-3}$ & 2.00  & 3.18 & 0.65 & 200  & 10 & 10 & 10 & $3.8 \times 10^{-4}$ \\
  		\hline
  	\end{tabular}
  	\label{tab:mlmc_mixed_1e4_sm}
  \end{table}
 \FloatBarrier

Figures \ref{fig:mlmc_results_0.5} and \ref{fig:mlmc_results_0.5_sm} present extra data on the performance of the MLMC-HQD and MLMC-HSM algorithms calculating $F=\mathcal{F}_D$ in Test 2 with $c_2=0.5$ in case of $K_{\ell}=10^4$ and $\varepsilon =10^{-3}$.
Figure \ref{mean-F-ell-test2-1e4-loqd} demonstrates plots of $\big< F_{\ell} \big>=\big< F_{0} \big> + \sum_{\ell^{\prime}=1}^{\ell} \big< \Delta F_{\ell^{\prime} } \big>$ and $\big< \Delta F_{\ell} \big>$.
The plots of variances $\mathbb{V}\big[\big< F_{\ell} \big>\big]$ and $\mathbb{V}\big[\big< \Delta F_{\ell} \big>\big]$ are shown in Fig. \ref{var-F-ell-test2-1e4-loqd}.
We observe a monotonic decrease in the mean value of the correction $\big< \Delta F_{\ell} \big>$ and its variance as the algorithm moves through levels. 
The estimate of costs, $C_{\ell}$,  for a sample increases due to refinement of the computational grid at each level (see Fig. \ref{cost-test2-1e4-loqd}).
 To analyze the convergence of the variance, we compute the kurtosis defined by
\begin{equation}
    \kappa_{\ell} = \mathbb{E} \left[ \left(\frac{F -\mathbb{E} \big[F_{\ell}\big]}{\sigma\big[ F_{\ell}\big] } \right)^4 \right] 
\end{equation}
for the initial number of realizations $\bm{N}^{ini}$ (see Fig. \ref{kappa-test2-1e4-loqd}). 
Kurtosis demonstrates the convergence of the variance by giving an order of the number of samples needed for convergence, i.e., $N_\ell = \mathcal{O}(\kappa_{\ell})$.
The results we obtained demonstrate that we have run a sufficient number of simulations to estimate the variance of our functional.
To ensure the validity of the telescoping summation, we perform a consistency check by computing \citep{giles-AN-2015}
\begin{equation}
    \eta_{\ell} = \frac{ \big<F_{\ell-1}\big> - \big<F_{\ell}\big> + \big<\Delta F_{\ell}\big>}{3\Big(\sqrt{\mathbb{V}[F_{\ell-1}]} + \sqrt{\mathbb{V} [F_{\ell}]}+ \sqrt{\mathbb{V}[\Delta F_\ell]}\Big)} \, .
\end{equation}
It is shown in Fig \ref{eta-test2-1e4-loqd}. 
$\eta_{\ell}$ should be less than $1$ otherwise estimates of functional $\Delta F_\ell$ are not being calculated correctly.
In our case, we demonstrate the validity of our implementation by showing $\eta_{\ell} < 1.0$ for all levels.
Figures \ref{fig:mlmc_results_0.9} and \ref{fig:mlmc_results_0.9_sm} shows results for one-region Test 1 with  $K_{\ell}=10^4$ and $\varepsilon = 1 \times 10^{-3}$.
The behavior in this test is largely the same as in the previous case, except that more samples were requested at the coarsest level.

The next problem we solve is Test 2 to calculate the functional $F~=~\mathcal{F}_{\tau_{8,0}}$ (Eq. \eqref{F-tau-i-0}), which is the integral of the scalar flux over the cell $\tau_{8,0}=[0.4375, 0.5]$ on the coarsest grid $G_0$.
This cell is adjacent to the domain center and lies at the interface between two regions.
The results are listed in  Tables \ref{tab:mlmc_1e4_cell_8} and \ref{tab:mlmc_1e4_cell_8_losm}.
For this study, we decreased $\varepsilon$ since the functional  $F~=~\mathcal{F}_{\tau_{8,0}}$ has a smaller value.
The results obtained are similar to those for the functional $\mathcal{F}_D$.

We now solve Test 2 to compute the vector of functionals $\{\mathcal{F}_{\tau_{i,0} }\}_{i=1}^{I_0}$ which are integrals of the scalar flux over each cell of the coarsest grid $G_0$.
To optimize these computations, the algorithm determines the number of realizations based on the cell with the maximum functional variance at each computational level.
Tables \ref{tab:mlmc_1e4_cell_max_var} and \ref{tab:mlmc_1e4_cell_max_var_losm} present the obtained results.
In this case, the algorithms require more realizations compared to calculations with optimization of a single functional $\mathcal{F}_{\tau_{8,0} }$ (see Tables \ref{tab:mlmc_1e4_cell_8}  and \ref{tab:mlmc_1e4_cell_8_losm}).
This difference is expected, as the coarse cell in the test above was not the one with the highest variance.
The convergence parameters  $\alpha$, $\beta$, and the weak convergence are calculated for every  functional $\mathcal{F}_{\tau_{i,0}}$.
The minimum values of $\alpha$ and $\beta$ are listed.
Values of $\alpha$ are slightly smaller than in the previous tests, and $\beta$ of the MLMC-HQD method is significantly smaller.
We note that $\beta > \gamma$ remains true for every computational cell.
In addition, the maximum weak-convergence check across all cells is analyzed to ensure that the error converges throughout the problem domain.
The results show the weak convergence check fails for MLMC-HQD method in the case $c_2=0.1$ and $\varepsilon=5\times 10^{-5}$ ( $\max W_{\hat \ell, i} > \frac{\varepsilon}{\sqrt{2}} = 3.5 \times 10^{-5}$) meaning additional computational level should be added to guarantee convergence under the MLMC theorem. 
     
\begin{figure*}[h!]
    \centering
    \begin{subfigure}[b]{0.49\textwidth}
        \centering
        \includegraphics[width=\textwidth]{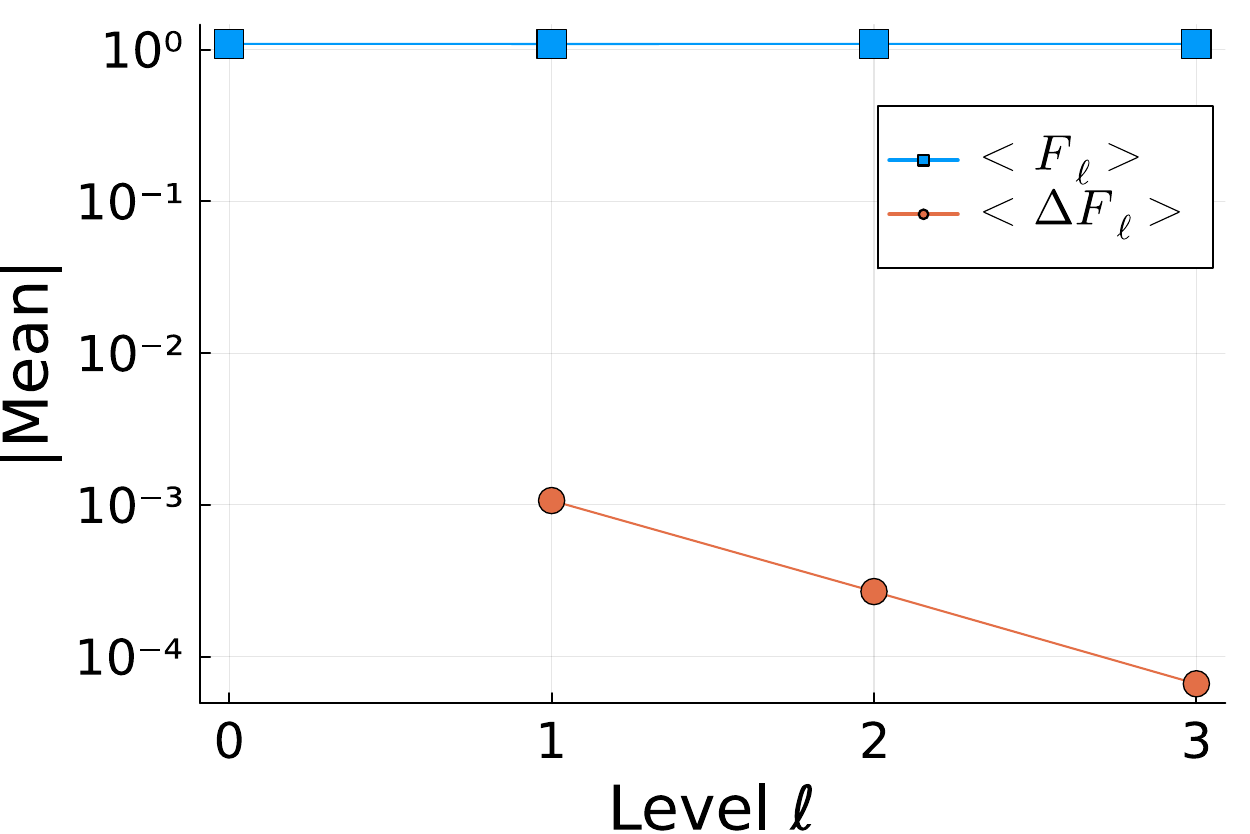}
        \caption{ $\big< F_{\ell} \big>$ and $\big< \Delta F_{\ell} \big>$ \label{mean-F-ell-test2-1e4-loqd}}
    \end{subfigure}
    \begin{subfigure}[b]{0.49\textwidth}
        \centering
        \includegraphics[width=\textwidth]{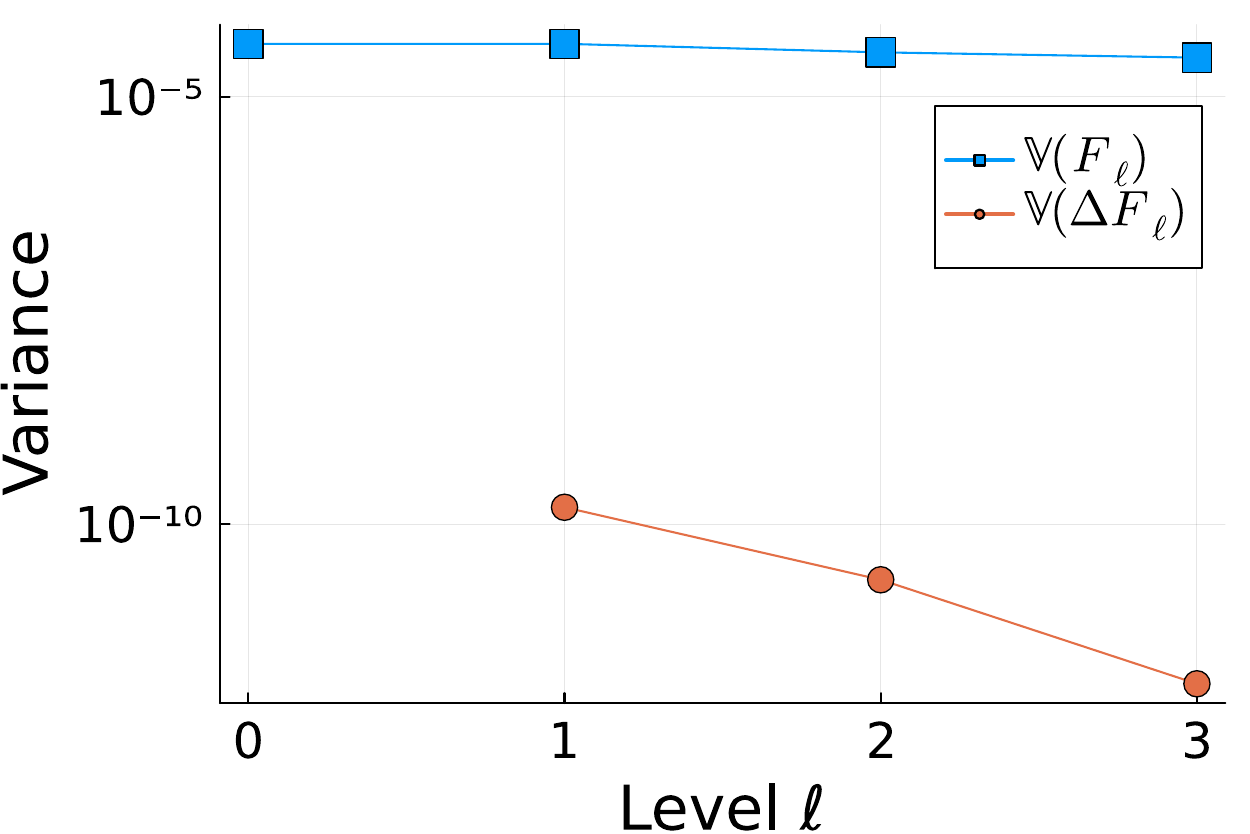}
        \caption{ $\mathbb{V}\big[\big< F_{\ell} \big>\big]$ and $\mathbb{V}\big[\big< \Delta F_{\ell} \big>\big]$ 
        \label{var-F-ell-test2-1e4-loqd}}
    \end{subfigure}
    \begin{subfigure}[b]{0.49\textwidth}
        \centering
        \includegraphics[width=\textwidth]{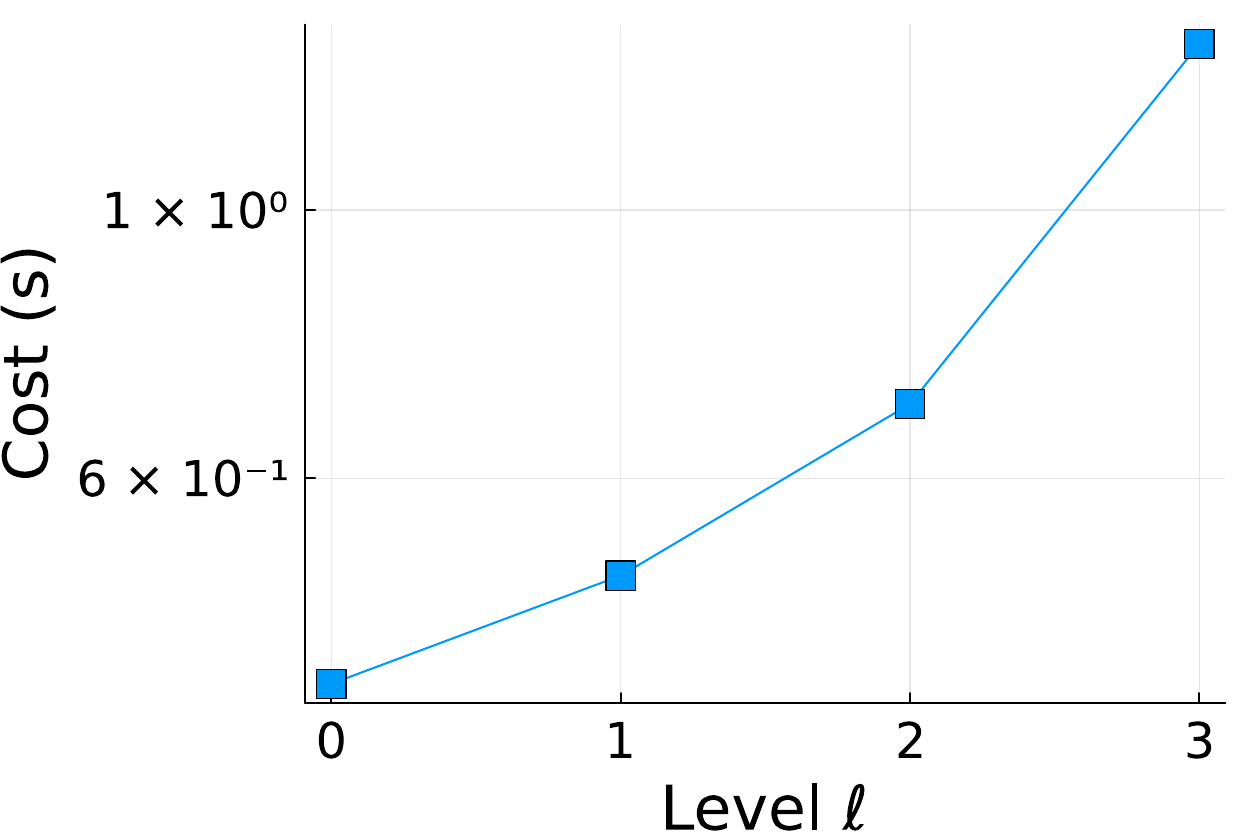}
        \caption{Computational costs, $C_{\ell}$  \label{cost-test2-1e4-loqd}}
    \end{subfigure}
    \begin{subfigure}[b]{0.49\textwidth}
        \centering
        \includegraphics[width=\textwidth]{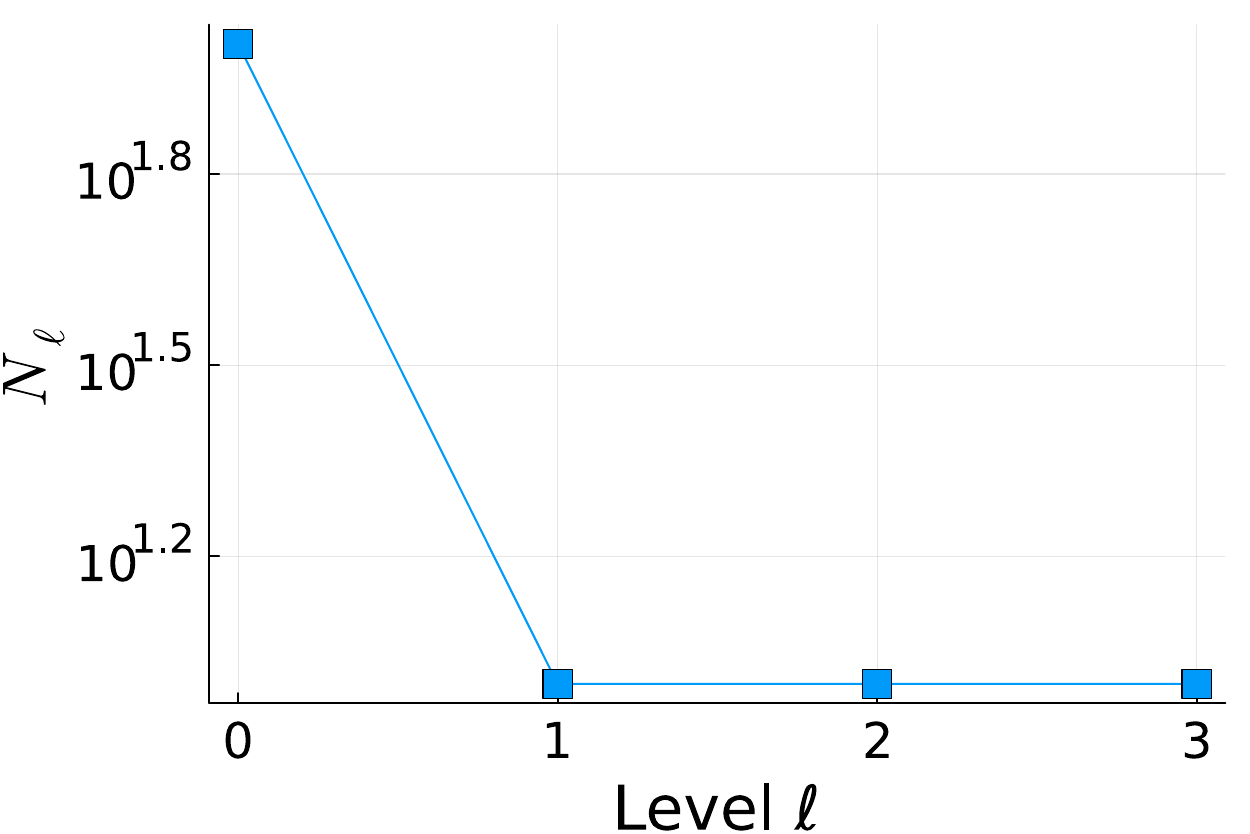}
        \caption{Estimates of number of realizations, $N_{\ell}$ \label{n-ell-test2-1e4-loqd} }
    \end{subfigure}
    \begin{subfigure}[b]{0.49\textwidth}
        \centering
        \includegraphics[width=\textwidth]{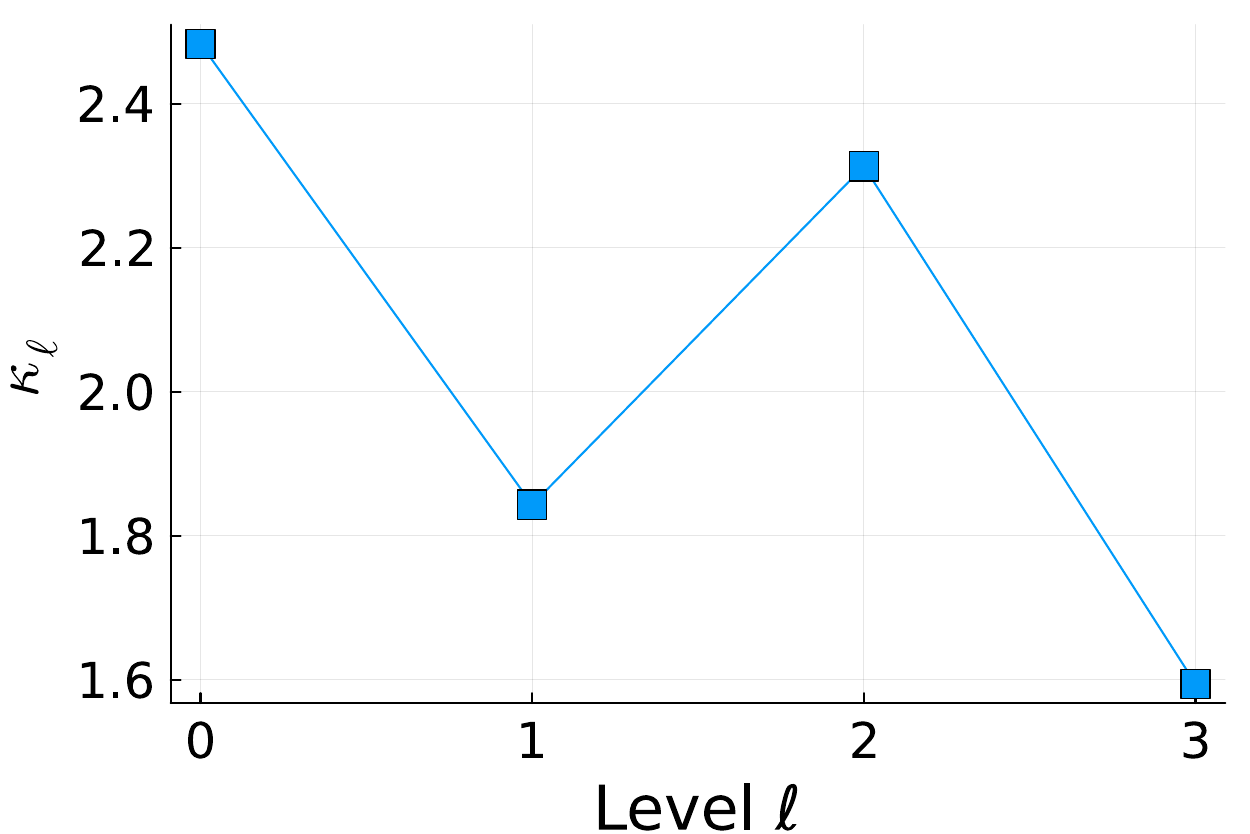}
        \caption{Kurtosis, $\kappa_{\ell}$ 
        \label{kappa-test2-1e4-loqd} }
    \end{subfigure}
    \begin{subfigure}[b]{0.49\textwidth}
        \centering
        \includegraphics[width=\textwidth]{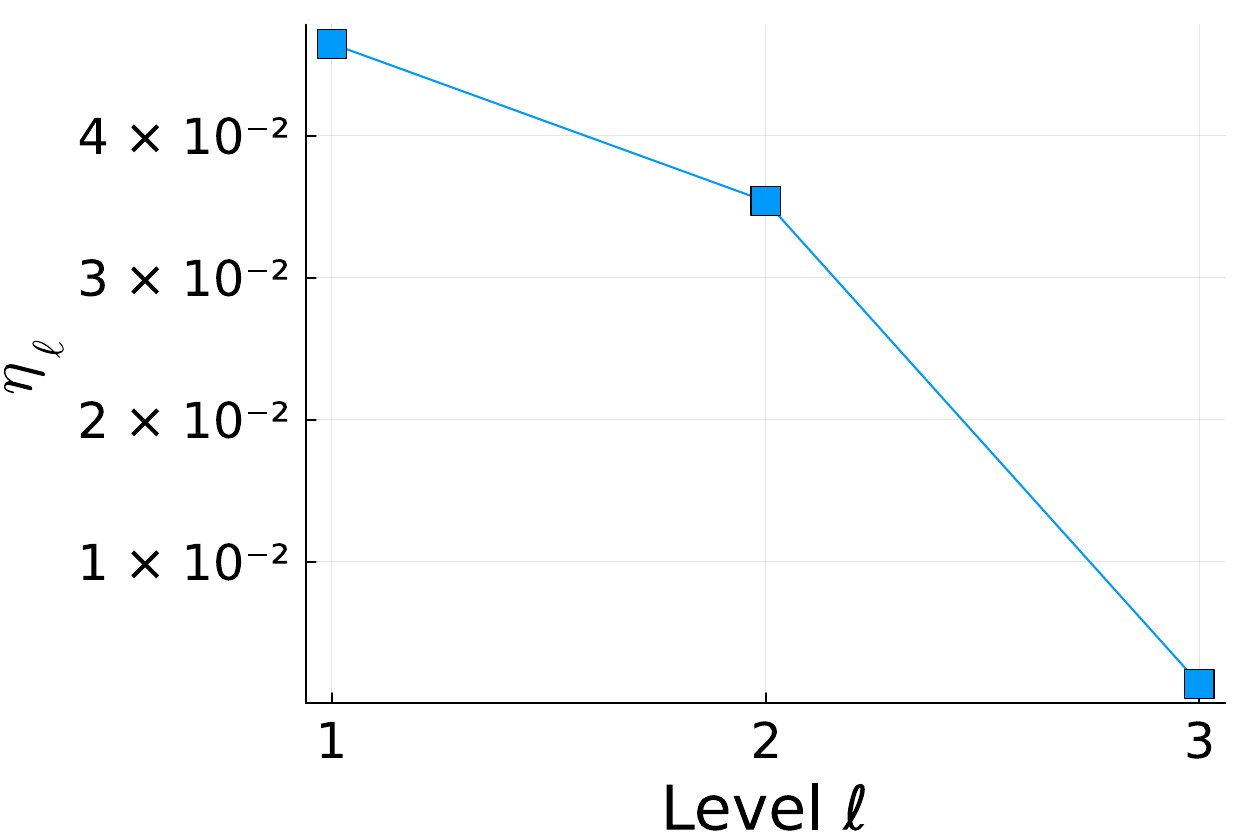}
        \caption{Consistency check, $\eta_{\ell}$ 
        \label{eta-test2-1e4-loqd} }
    \end{subfigure}
    \caption{Test 2, $c_2=0.5$,  $F=\mathcal{F}_D$.  Data on convergence of $\big< F\big>$ computed  by the MLMC-HQD algorithm with $K_{\ell}=10^4$  for $\varepsilon = 1 \times 10^{-3}$.}
    \label{fig:mlmc_results_0.5}
\end{figure*}

\begin{figure*}[h!]
    \centering
    \begin{subfigure}[b]{0.49\textwidth}
        \centering
        \includegraphics[width=\textwidth]{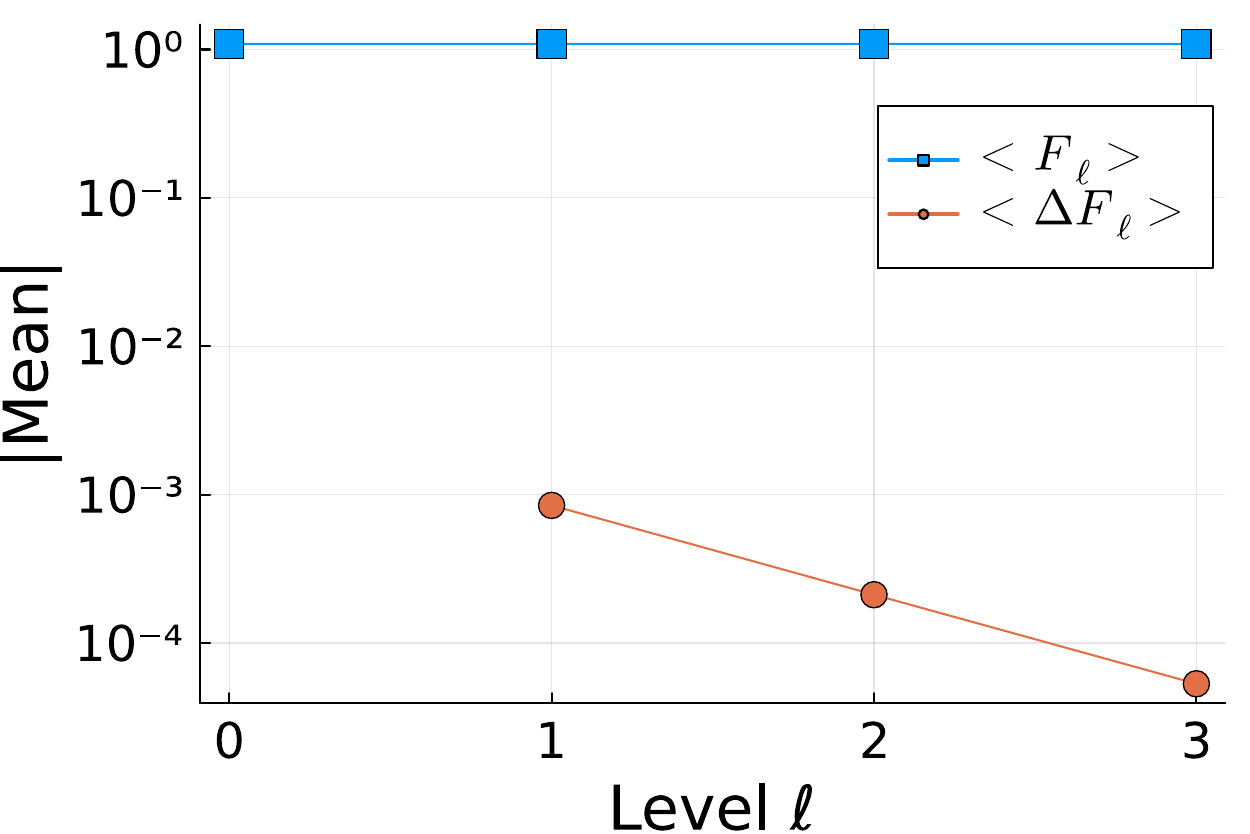}
        \caption{ $\big< F_{\ell} \big>$ and $\big< \Delta F_{\ell} \big>$ \label{mean-F-ell-test2-1e4-losm}}
    \end{subfigure}
    \begin{subfigure}[b]{0.49\textwidth}
        \centering
        \includegraphics[width=\textwidth]{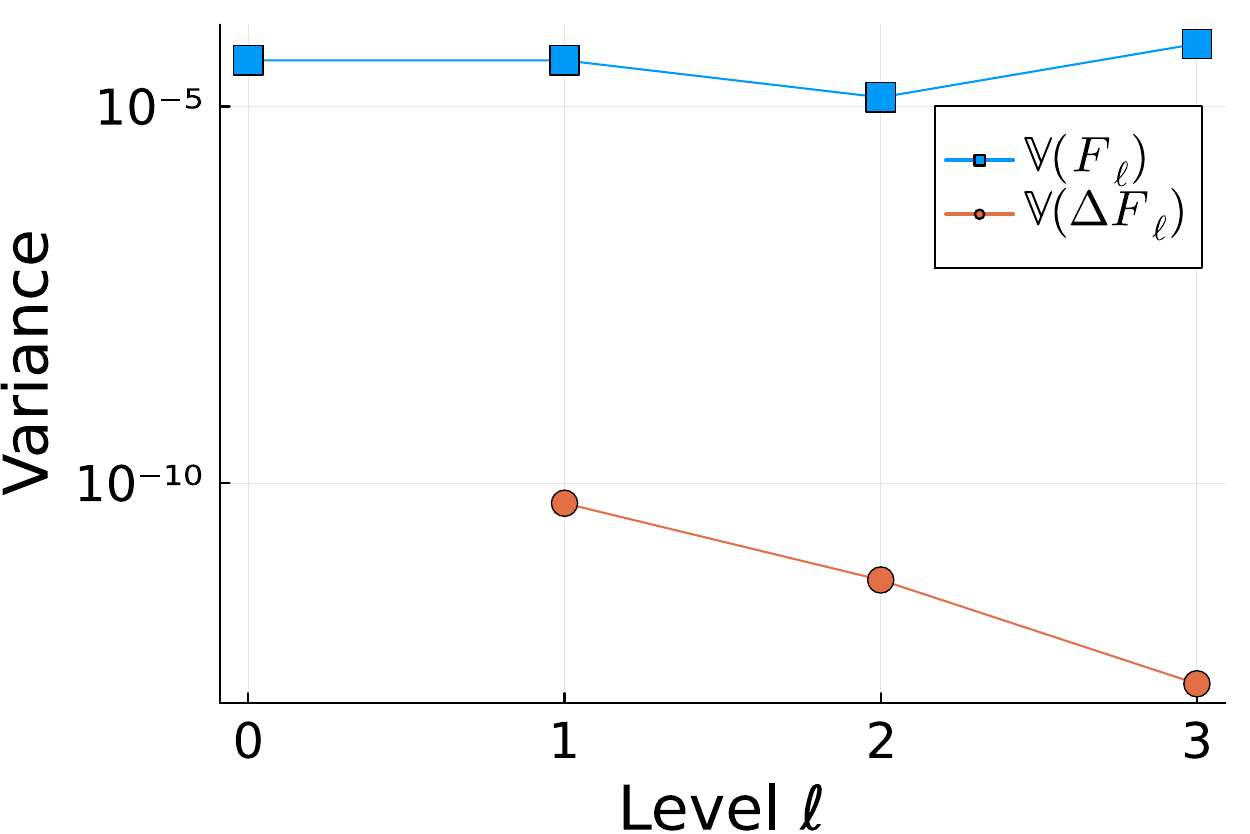}
        \caption{ $\mathbb{V}\big[\big< F_{\ell} \big>\big]$ and $\mathbb{V}\big[\big< \Delta F_{\ell} \big>\big]$
        \label{var-F-ell-test2-1e4-losm}}
    \end{subfigure}
    \begin{subfigure}[b]{0.49\textwidth}
        \centering
        \includegraphics[width=\textwidth]{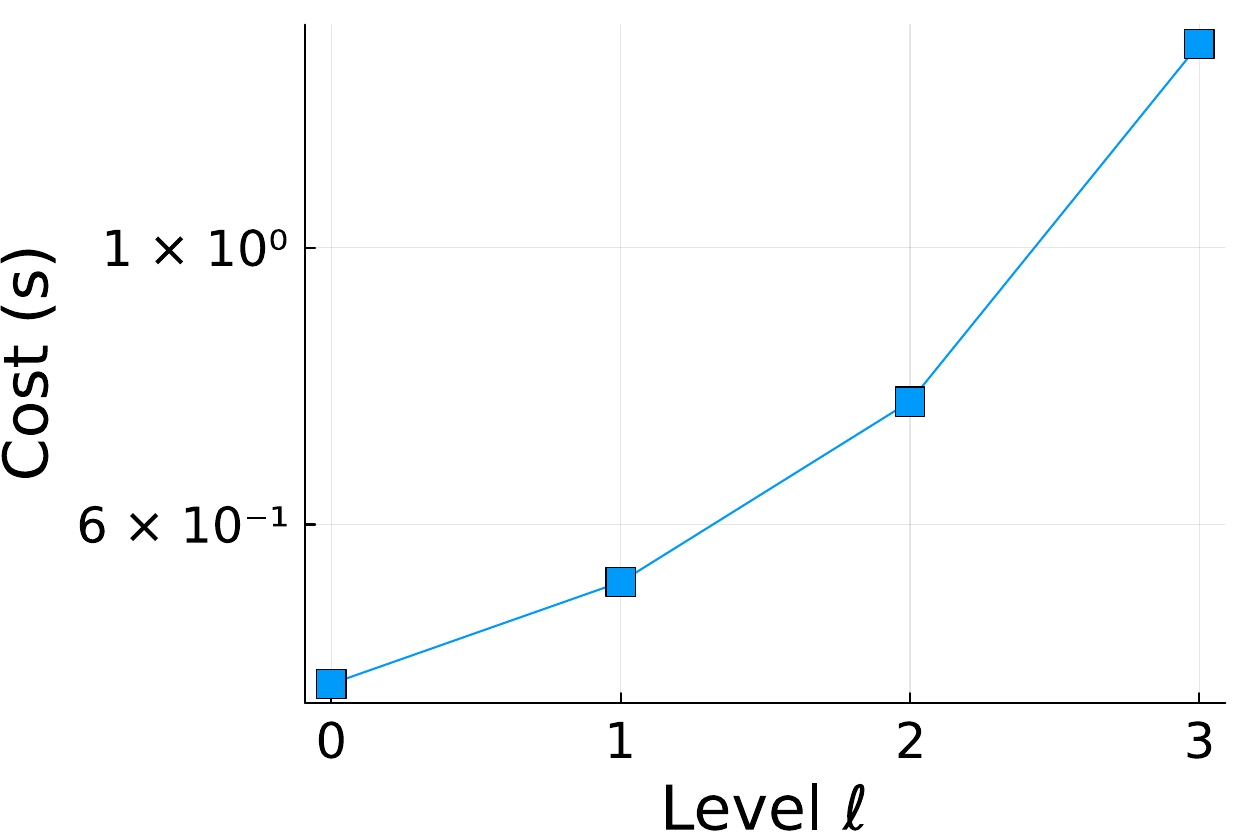}
        \caption{Computational costs, $C_{\ell}$  \label{cost-test2-1e4-losm}}
    \end{subfigure}
    \begin{subfigure}[b]{0.49\textwidth}
        \centering
        \includegraphics[width=\textwidth]{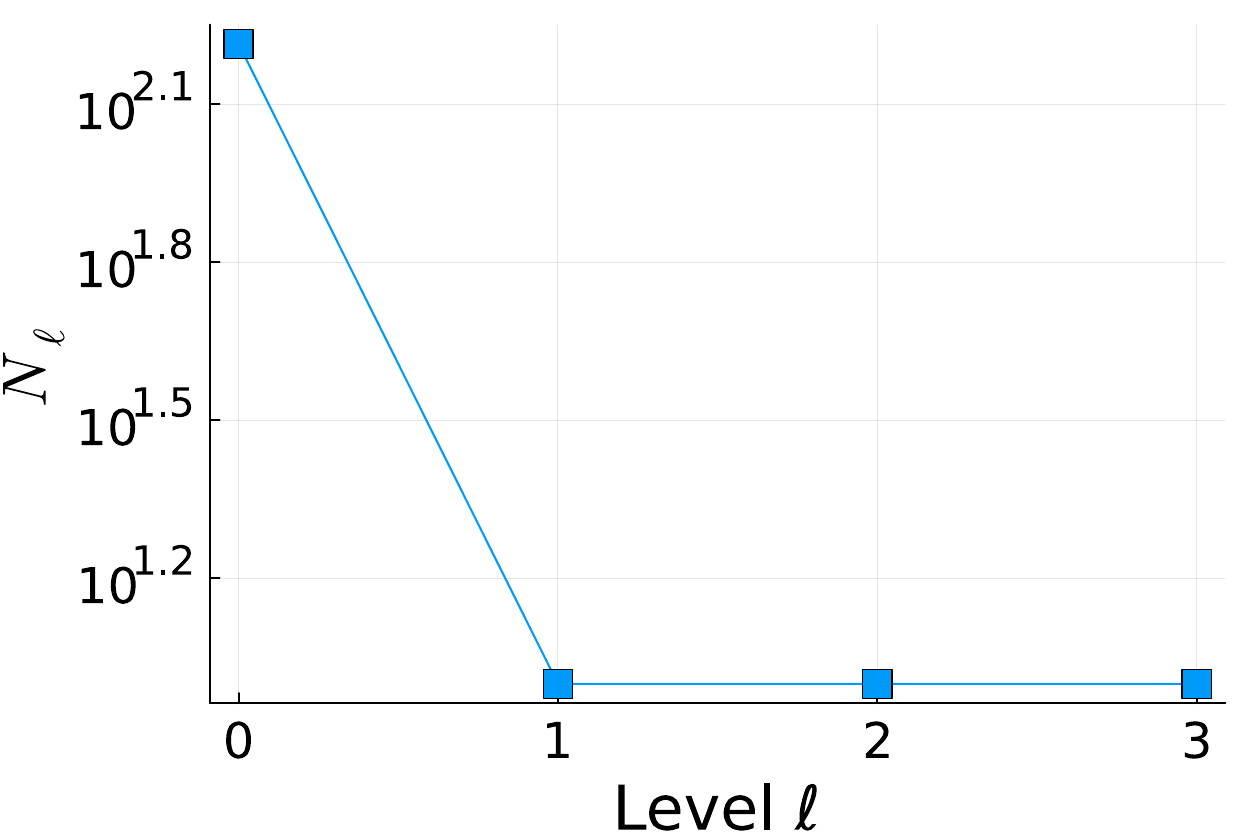}
        \caption{Estimates of number of realizations, $N_{\ell}$ \label{n-ell-test2-1e4-losm} }
    \end{subfigure}
    \begin{subfigure}[b]{0.49\textwidth}
        \centering
        \includegraphics[width=\textwidth]{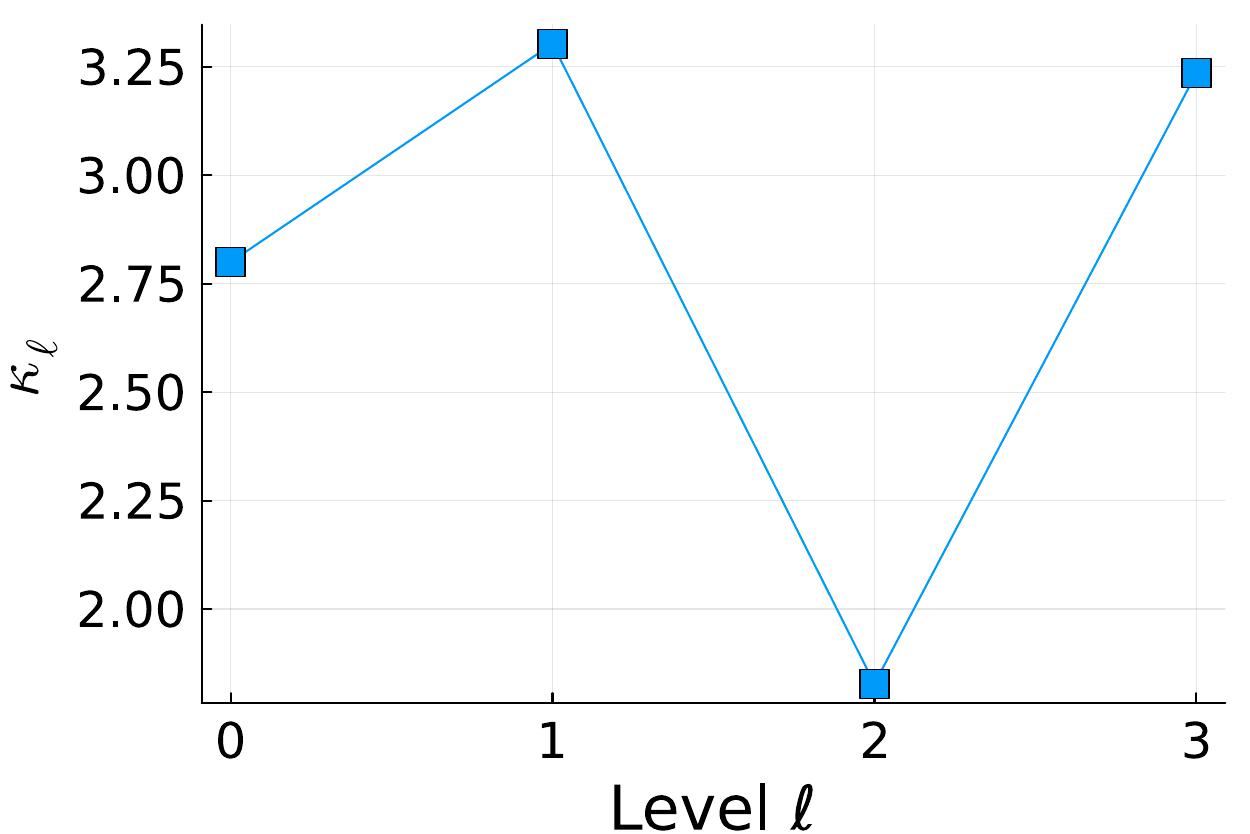}
        \caption{Kurtosis, $\kappa_{\ell}$ \
        \label{kappa-test2-1e4-losm} }
    \end{subfigure}
    \begin{subfigure}[b]{0.49\textwidth}
        \centering
        \includegraphics[width=\textwidth]{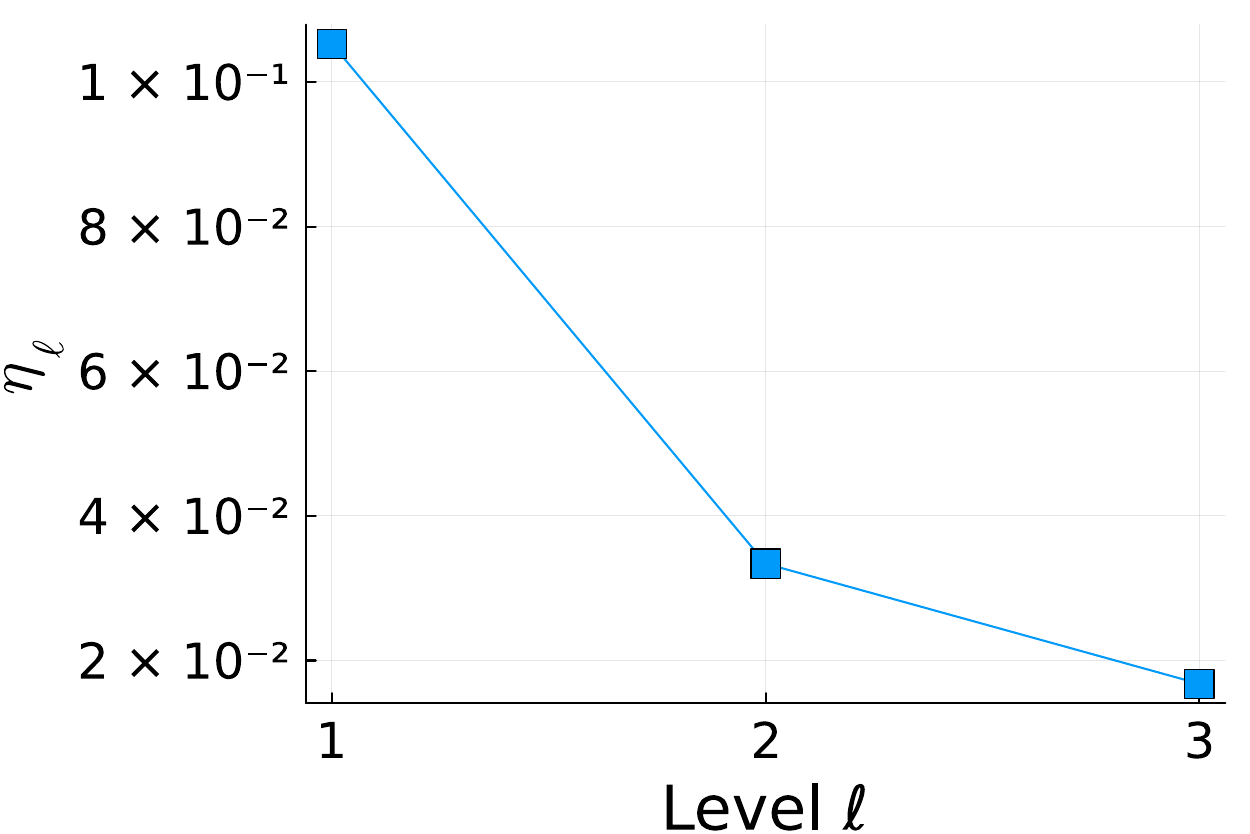}
        \caption{Consistency check, $\eta_{\ell}$ 
        \label{eta-test2-1e4-losm} }
    \end{subfigure}
    \caption{Test 2, $c_2=0.5$,  $F=\mathcal{F}_D$.  Data on convergence of $\big< F\big>$ computed  by the MLMC-HSM algorithm with $K_{\ell}=10^4$  for $\varepsilon = 1 \times 10^{-3}$.}
    \label{fig:mlmc_results_0.5_sm}
\end{figure*}
\begin{figure*}[h!]
	\centering
	\begin{subfigure}[b]{0.49\textwidth}
		\centering
		\includegraphics[width=\textwidth]{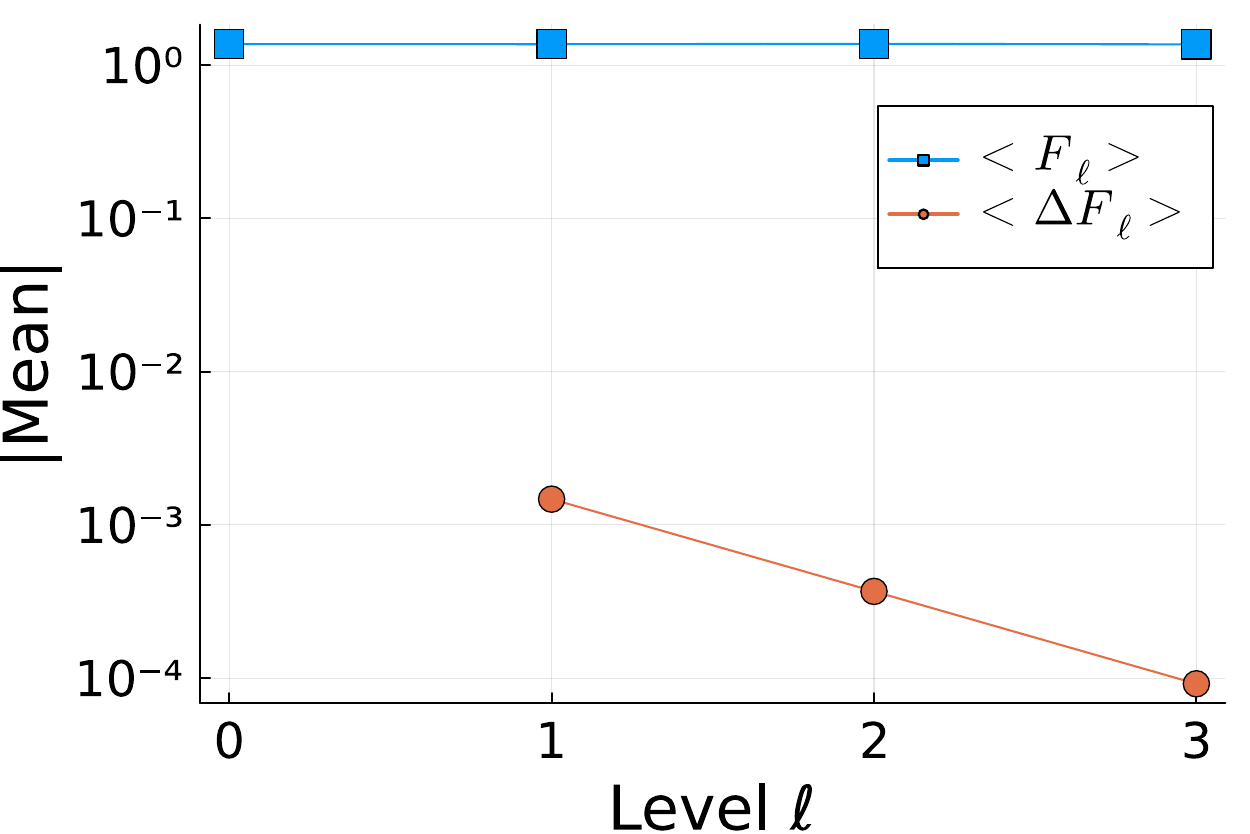}
		\caption{ $\big< F_{\ell} \big>$ and $\big< \Delta F_{\ell} \big>$ \label{mean-F-ell-test1-1e4-loqd}}
	\end{subfigure}
	\begin{subfigure}[b]{0.49\textwidth}
		\centering
		\includegraphics[width=\textwidth]{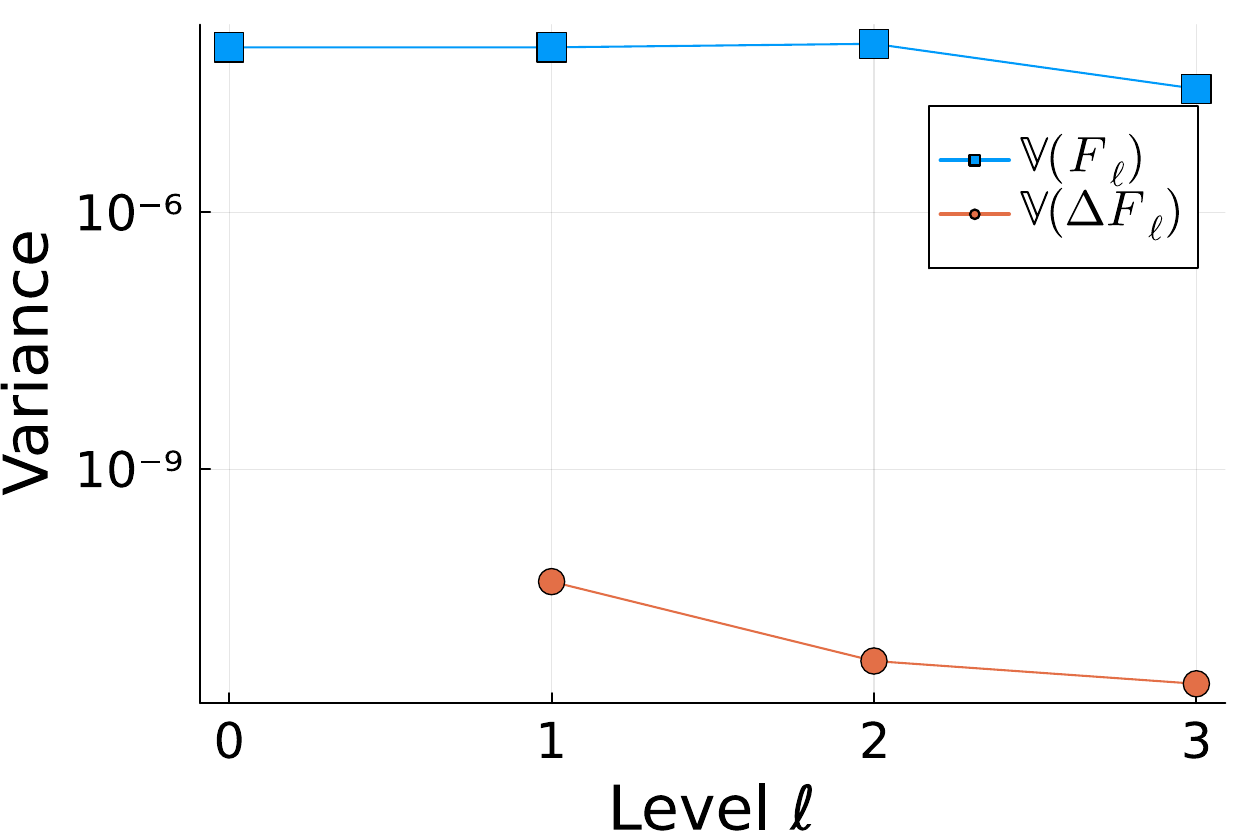}
			\caption{ $\mathbb{V}\big[\big< F_{\ell} \big>\big]$ and $\mathbb{V}\big[\big< \Delta F_{\ell} \big>\big]$ 
			\label{var-F-ell-test1-1e4-loqd}}
	\end{subfigure}
	\begin{subfigure}[b]{0.49\textwidth}
		\centering
		\includegraphics[width=\textwidth]{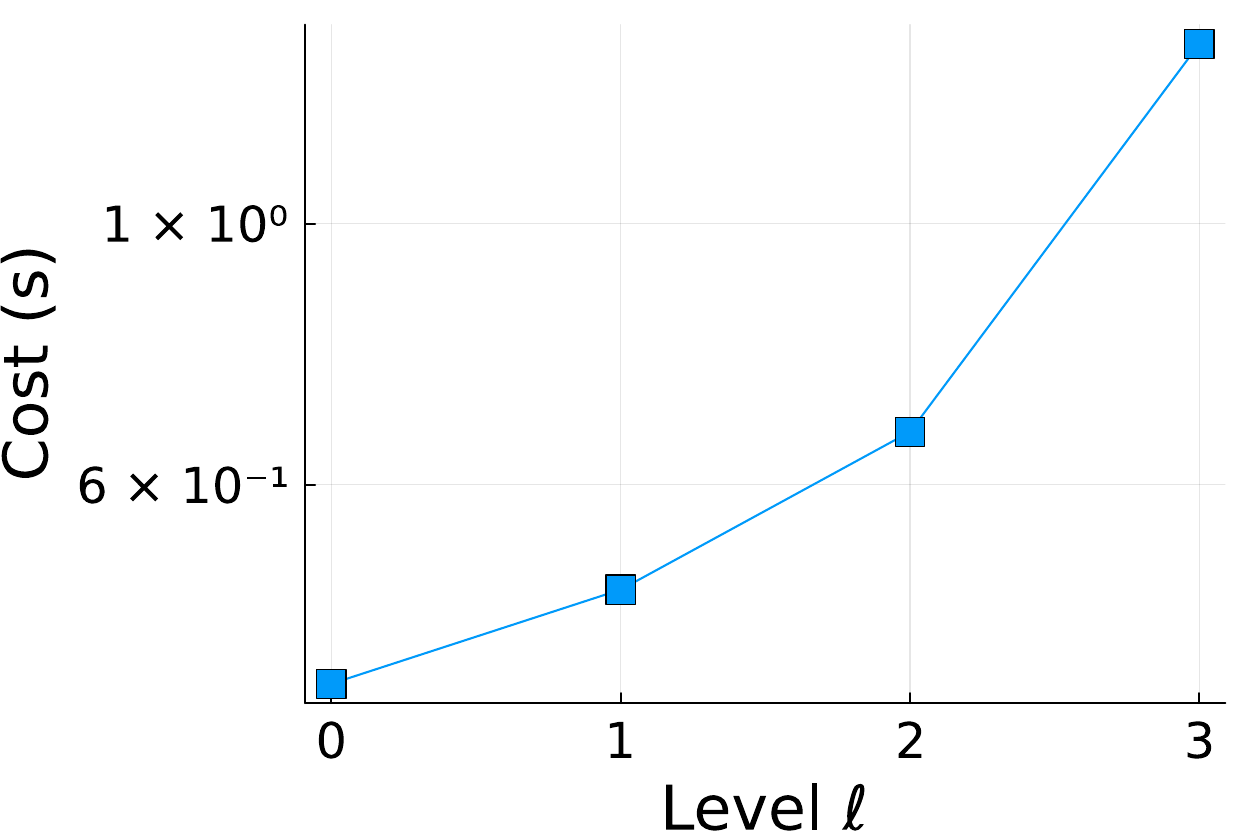}
		\caption{Computational costs, $C_{\ell}$  \label{cost-test1-1e4-loqd}}
	\end{subfigure}
	\begin{subfigure}[b]{0.49\textwidth}
		\centering
		\includegraphics[width=\textwidth]{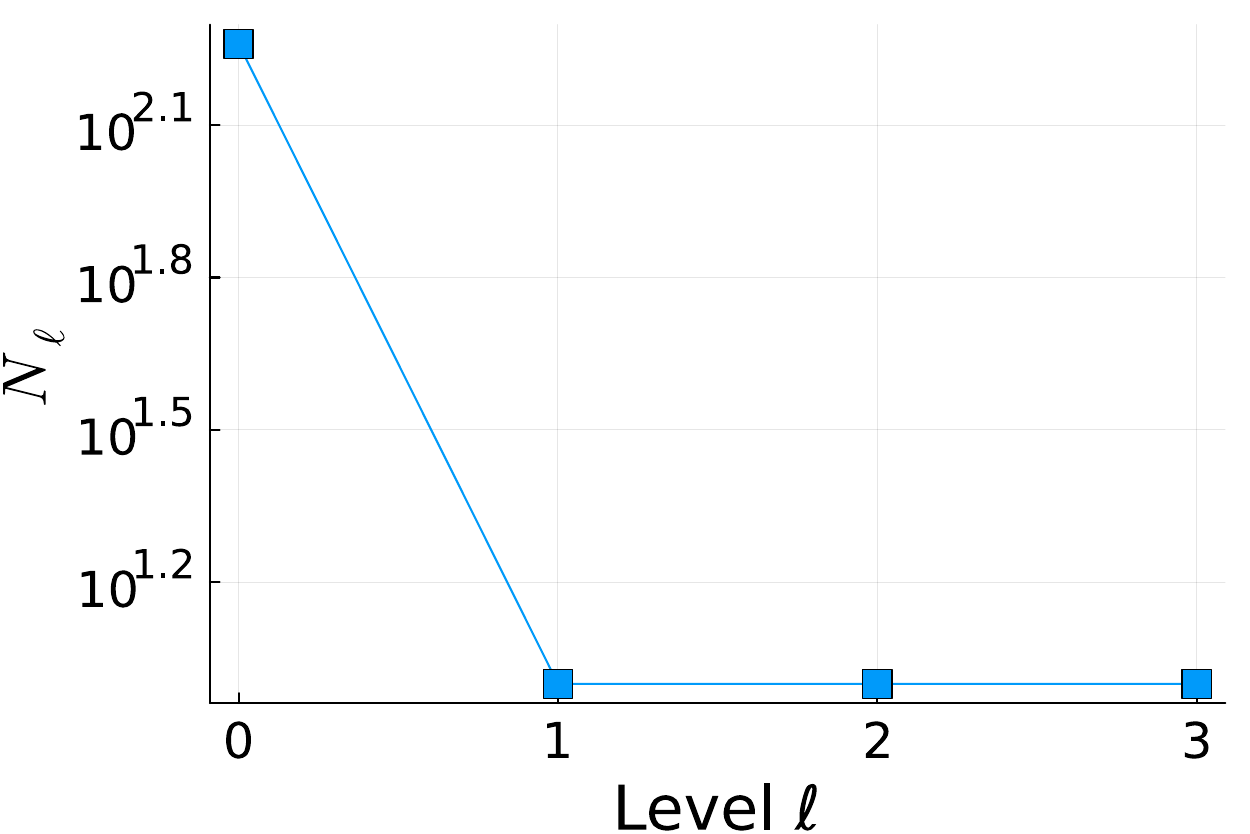}
		\caption{Estimates of number of realizations, $N_{\ell}$ \label{n-ell-test1-1e4-loqd} }
	\end{subfigure}
	\begin{subfigure}[b]{0.49\textwidth}
		\centering
		\includegraphics[width=\textwidth]{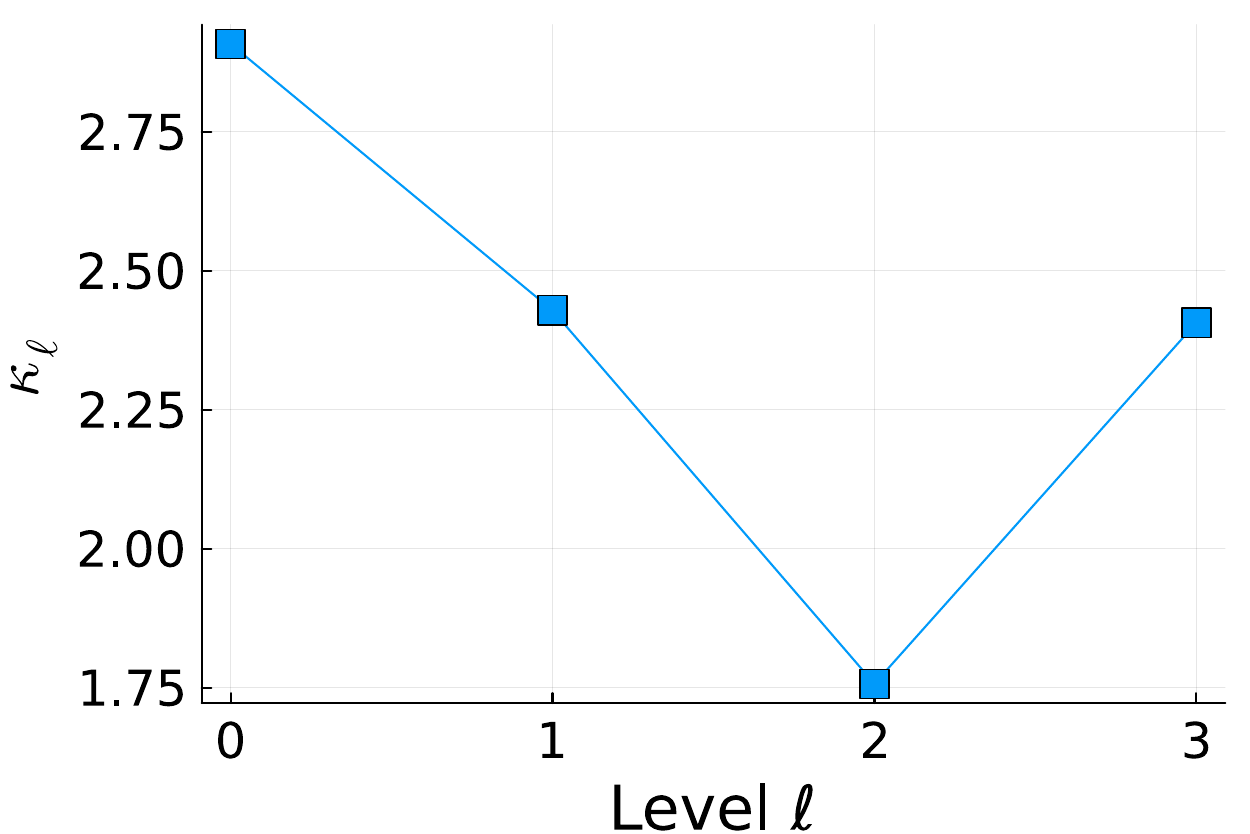}
		\caption{Kurtosis, $\kappa_{\ell}$ 
		\label{kappa-test1-1e4-loqd} }
	\end{subfigure}
	\begin{subfigure}[b]{0.49\textwidth}
		\centering
		\includegraphics[width=\textwidth]{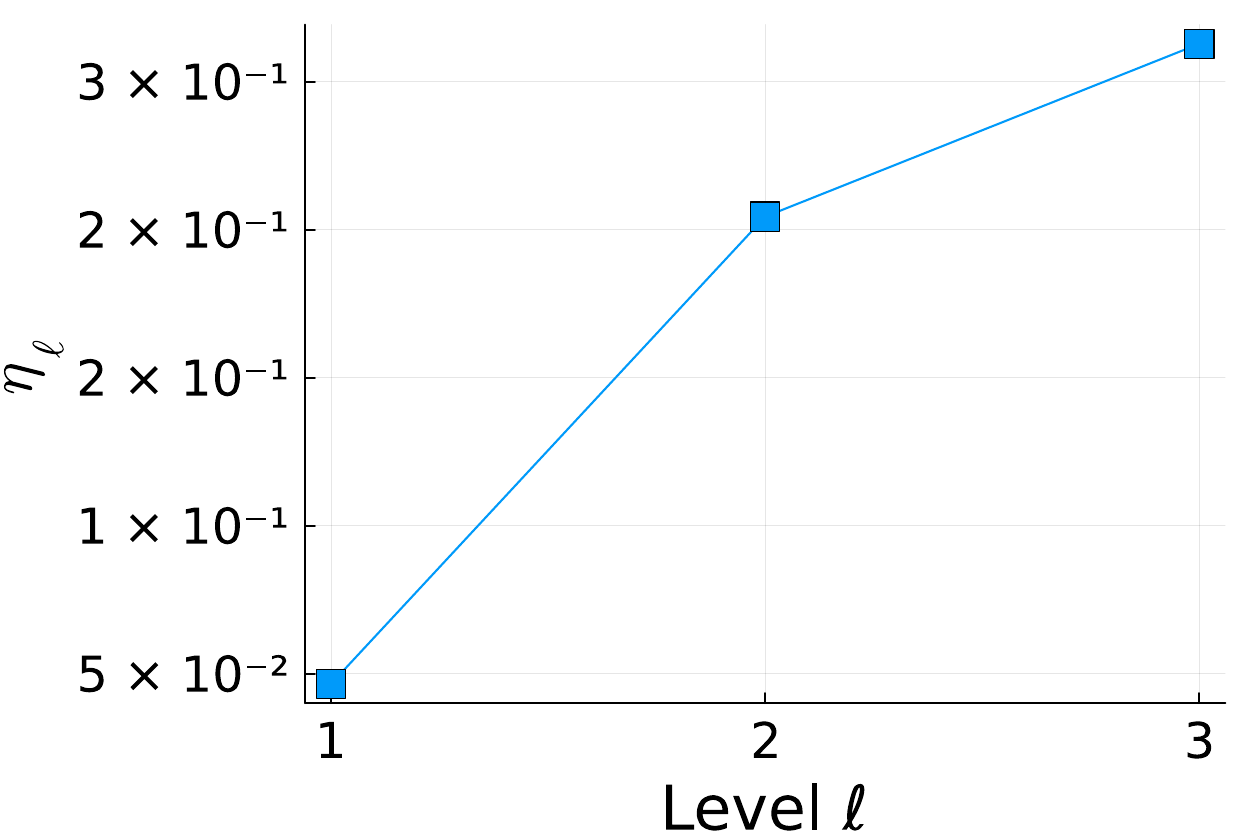}
		\caption{Consistency check, $\eta_{\ell}$  
		\label{eta-test1-1e4-loqd} }
	\end{subfigure}
	\caption{Test 1 $F=\mathcal{F}_D$.  Data on convergence of $\big< F\big>$ computed  by the MLMC-HQD algorithm with $K_{\ell}=10^4$  for $\varepsilon = 1 \times 10^{-3}$.}
	\label{fig:mlmc_results_0.9}
\end{figure*}

\begin{figure*}[h!]
	\centering
	\begin{subfigure}[b]{0.49\textwidth}
		\centering
		\includegraphics[width=\textwidth]{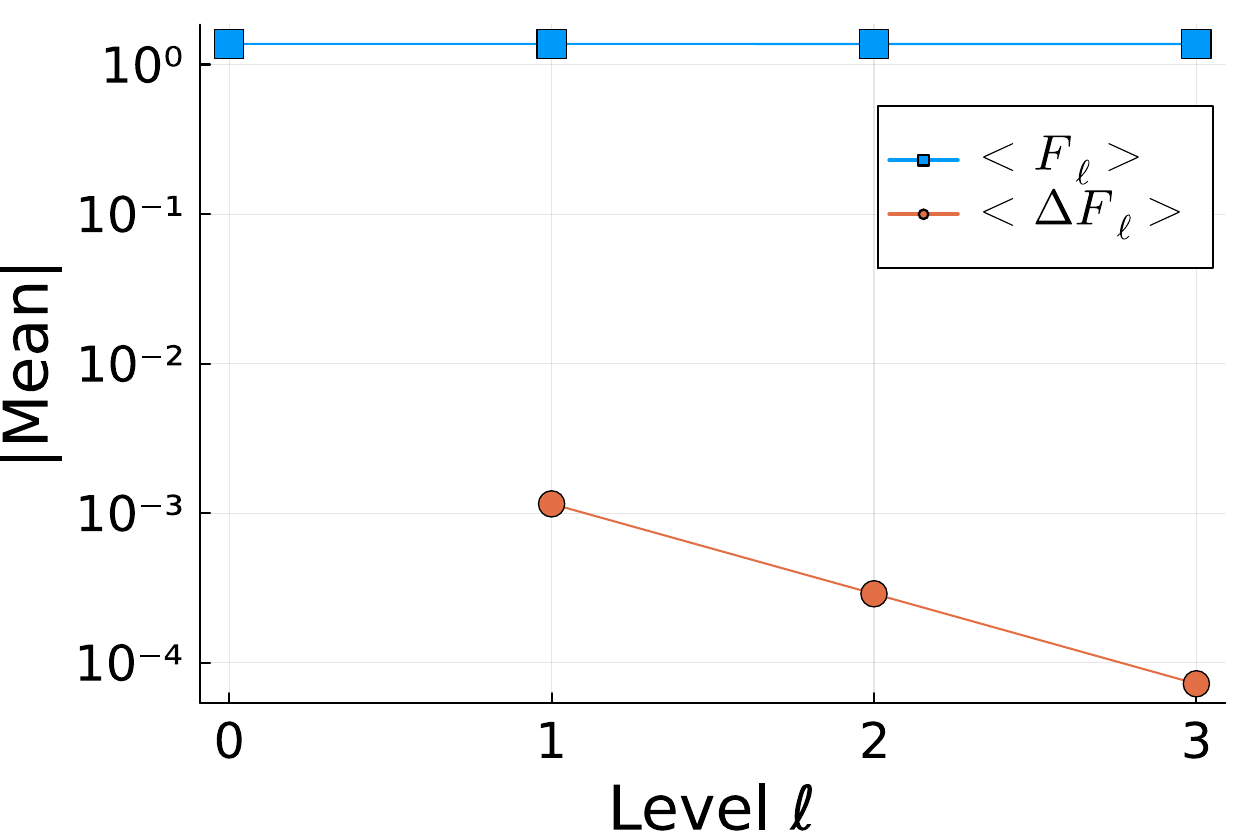}
		\caption{ $\big< F_{\ell} \big>$ and $\big< \Delta F_{\ell} \big>$ \label{mean-F-ell-test1-1e4-losm}}
	\end{subfigure}
	\begin{subfigure}[b]{0.49\textwidth}
		\centering
		\includegraphics[width=\textwidth]{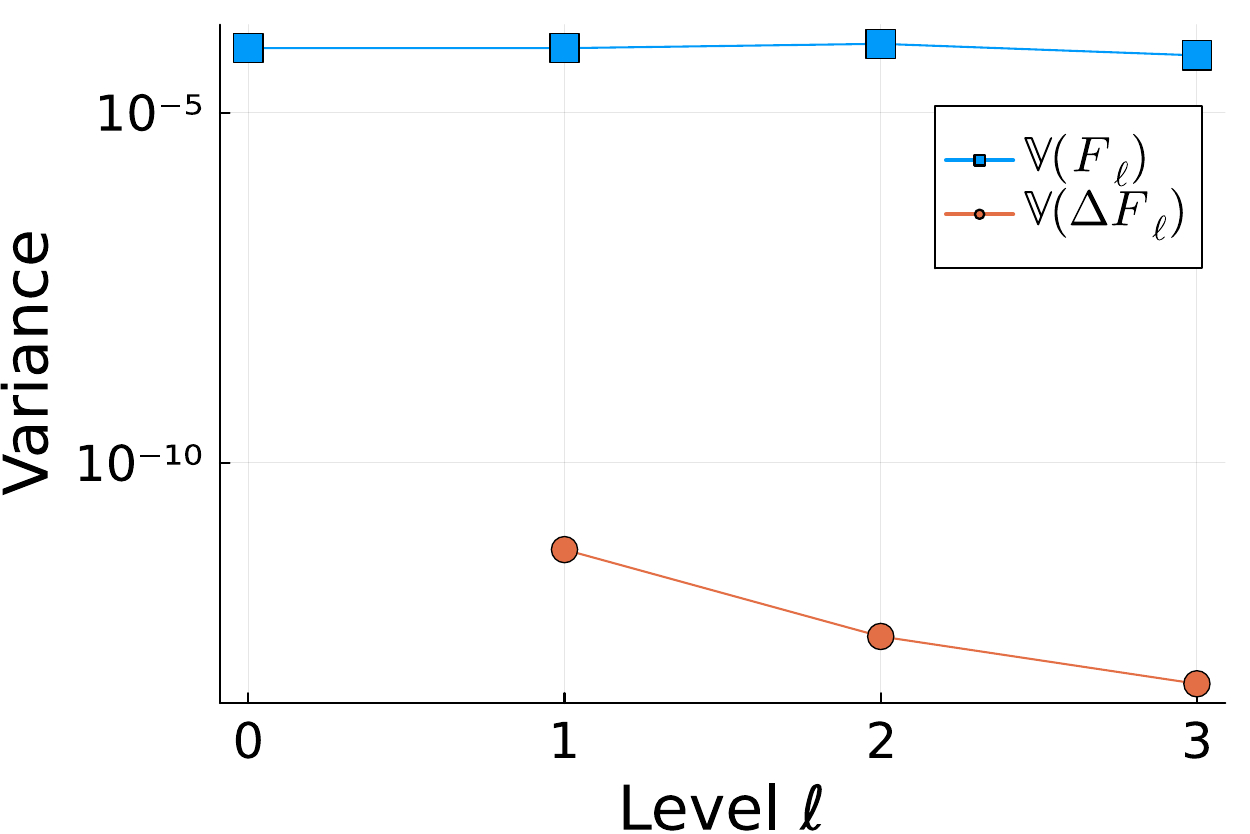}
			\caption{ $\mathbb{V}\big[\big< F_{\ell} \big>\big]$ and $\mathbb{V}\big[\big< \Delta F_{\ell} \big>\big]$
			\label{var-F-ell-test1-1e4-losm}}
	\end{subfigure}
	\begin{subfigure}[b]{0.49\textwidth}
		\centering
		\includegraphics[width=\textwidth]{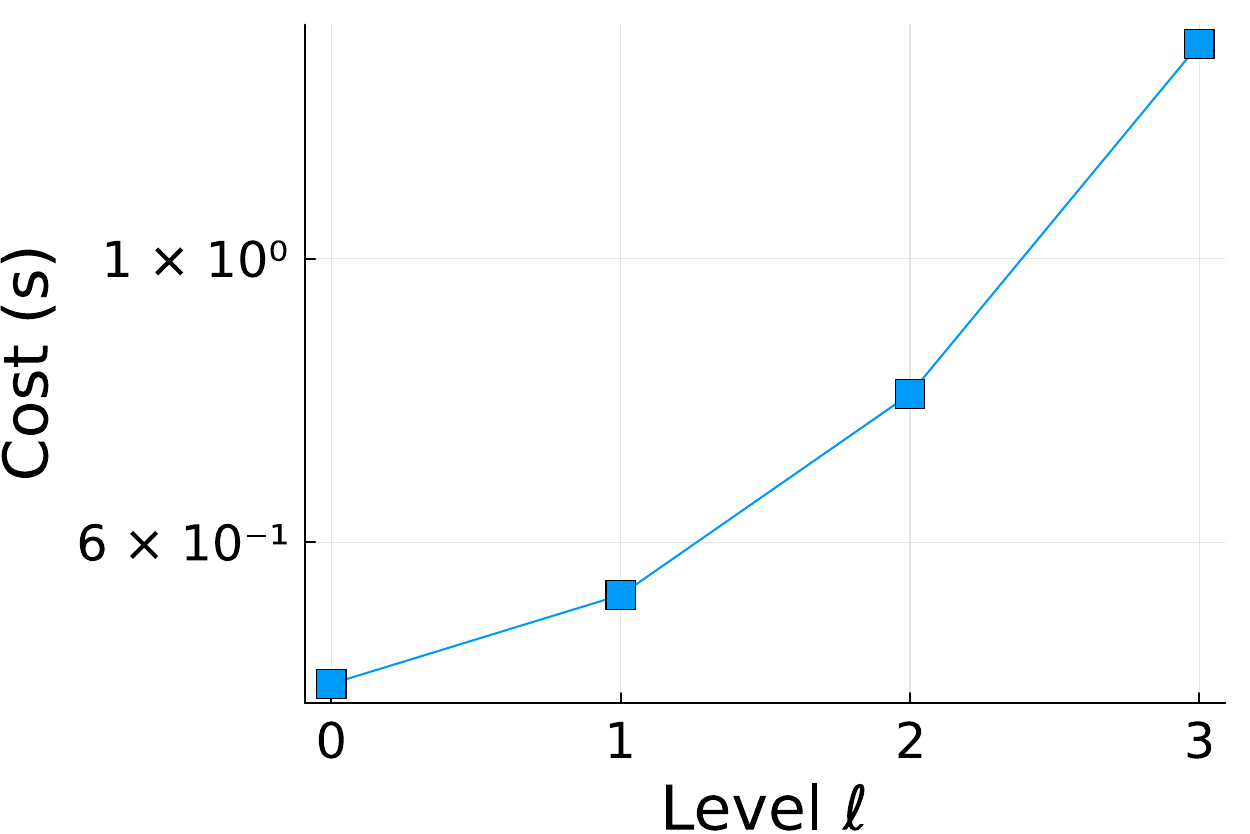}
		\caption{Computational costs, $C_{\ell}$  \label{cost-test1-1e4-losm}}
	\end{subfigure}
	\begin{subfigure}[b]{0.49\textwidth}
		\centering
		\includegraphics[width=\textwidth]{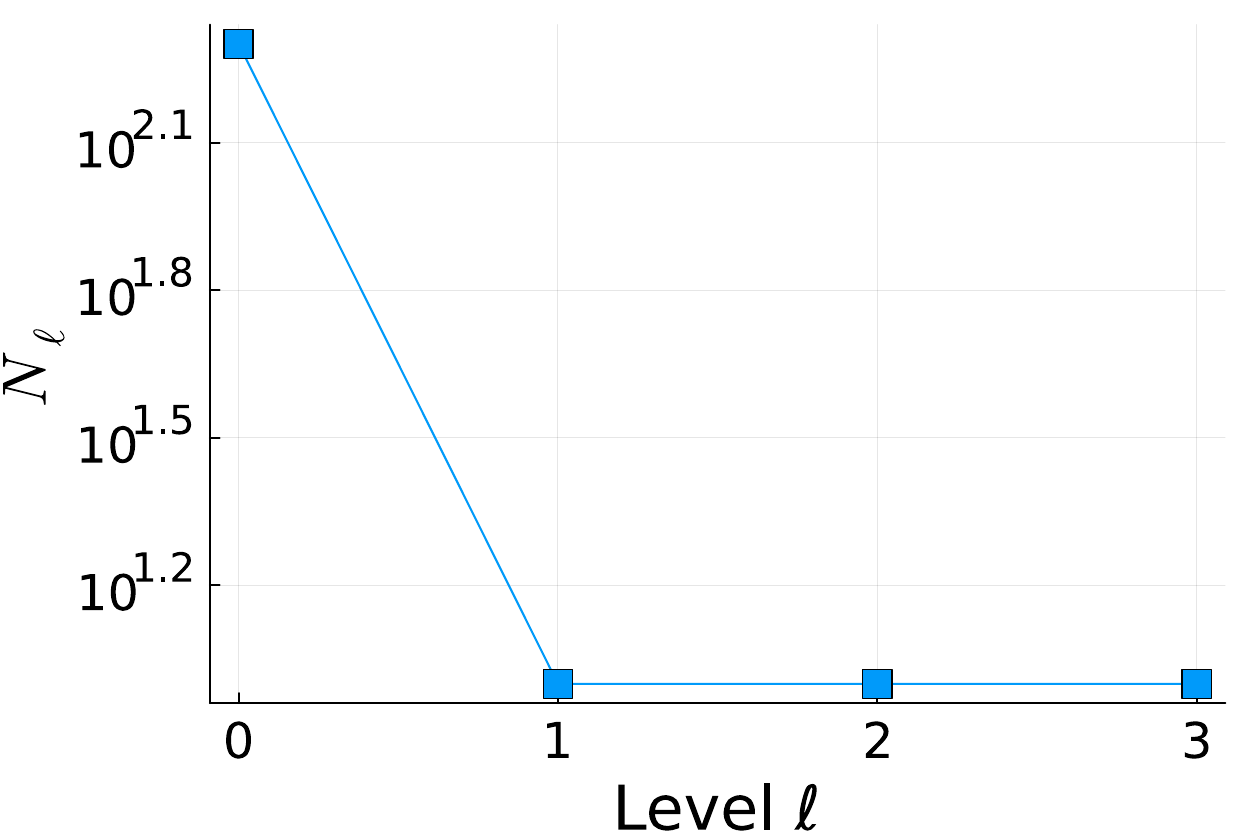}
		\caption{Estimates of number of realizations, $N_{\ell}$ \label{n-ell-test1-1e4-losm} }
	\end{subfigure}
	\begin{subfigure}[b]{0.49\textwidth}
		\centering
		\includegraphics[width=\textwidth]{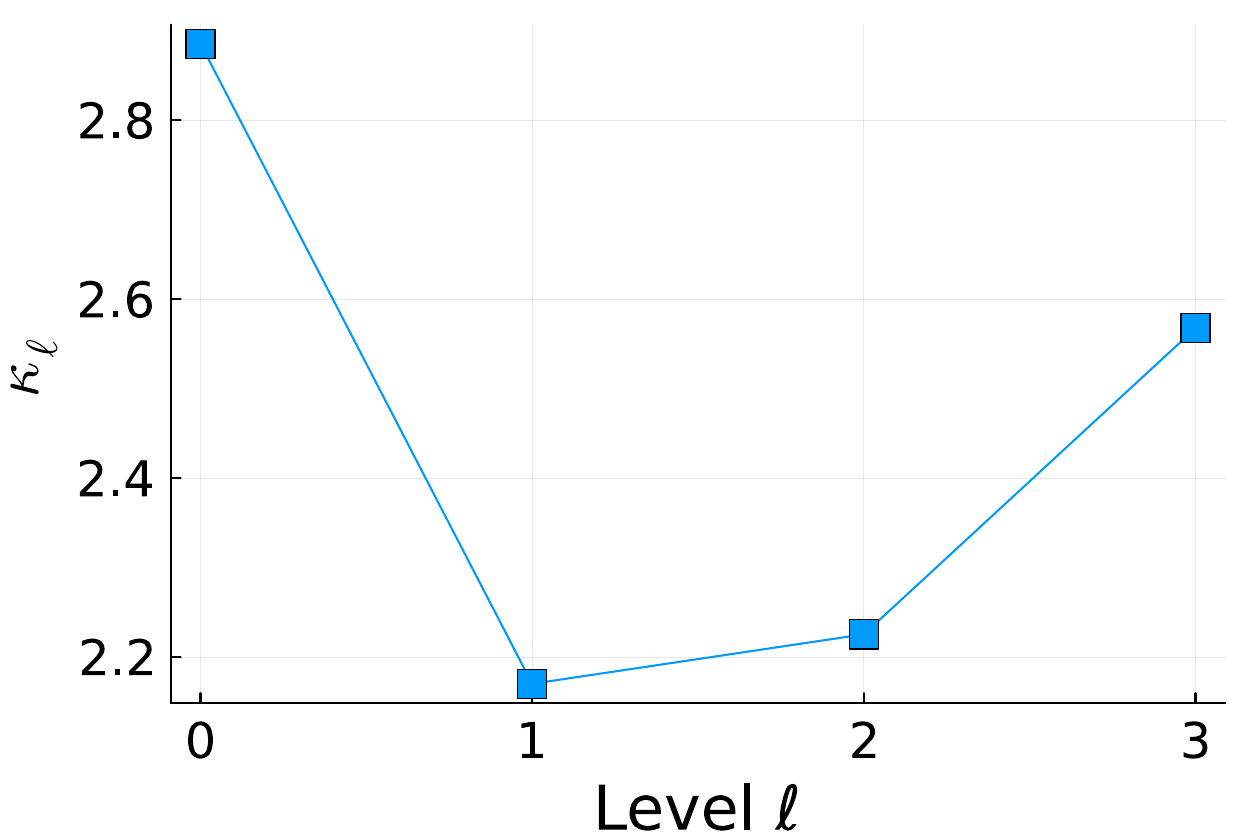}
		\caption{Kurtosis, $\kappa_{\ell}$
		\label{kappa-test1-1e4-losm} }
	\end{subfigure}
	\begin{subfigure}[b]{0.49\textwidth}
		\centering
		\includegraphics[width=\textwidth]{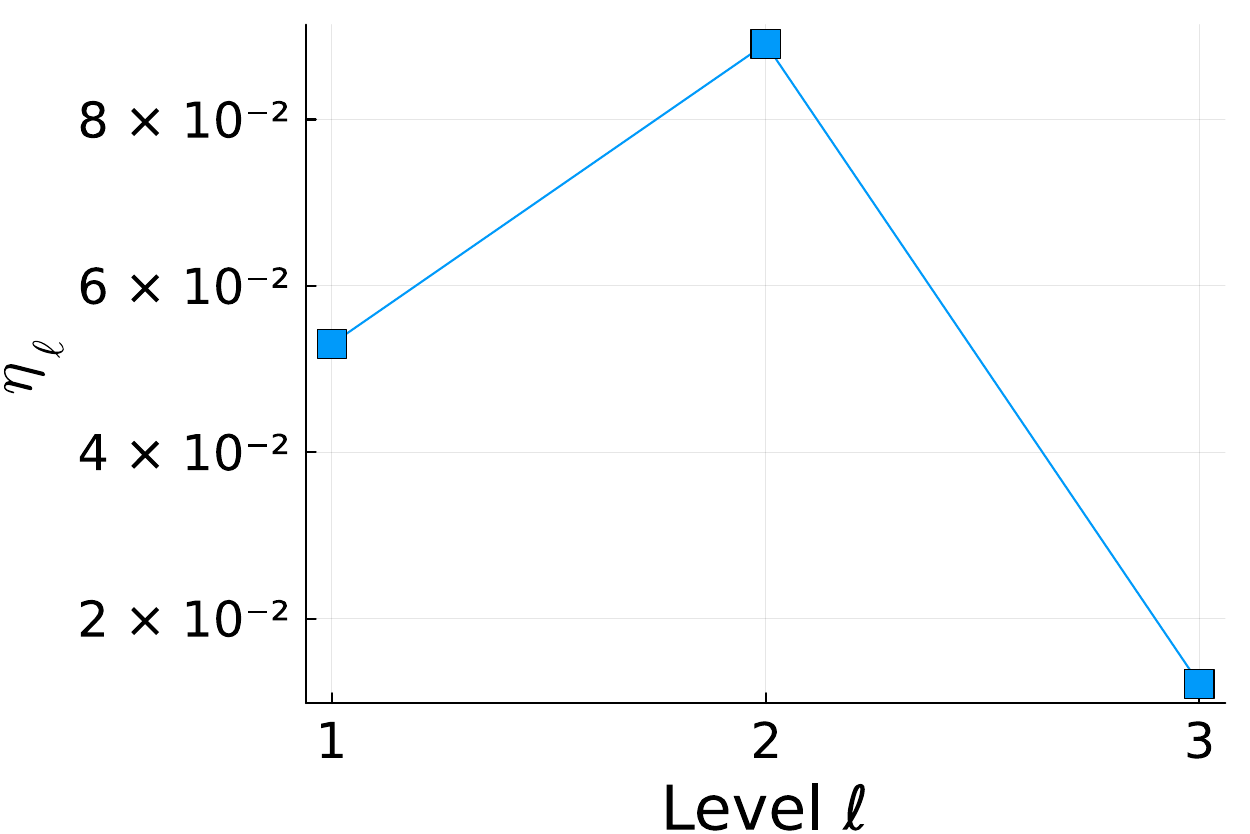}
		\caption{Consistency check, $\eta_{\ell}$
		\label{eta-test1-1e4-losm} }
	\end{subfigure}
    \caption{Test 1  $F=\mathcal{F}_D$.  Data on convergence of $\big< F\big>$ computed  by the MLMC-HSM algorithm with $K_{\ell}=10^4$  for $\varepsilon = 1 \times 10^{-3}$.}
\label{fig:mlmc_results_0.9_sm}
\end{figure*}

\FloatBarrier

\begin{table}[h!]
 	\centering
 	\caption{Test 2. $F=\mathcal{F}_{\tau_{8,0}}$,  MLMC-HQD,   $K_{\ell} = 10^4$}
 	\begin{tabular}{|c|c||c|c|c|c|c|c|c|c|}
 		\hline
 		$c_2$ & $\varepsilon$     & $\alpha$ & $\beta$ & $\gamma$ & $N_0$ & $N_1$ & $N_2$ & $N_3$ & $\max_{\hat \ell}{W_{\hat{\ell}}}$\\
 		\hline
 		    & $5 \times 10^{-4}$ & 2.00  & 3.25 & 0.66 & 11   & 10 & 10 & 10 & $3.5 \times 10^{-5}$ \\
 		0.1 & $1 \times 10^{-4}$ & 2.00  & 2.66 & 0.63 & 314  & 10 & 10 & 10 & $3.5 \times 10^{-5}$ \\
 	        & $5 \times 10^{-5}$ & 2.00  & 3.71 & 0.63 & 1257 & 10 & 10 & 10 & $3.5 \times 10^{-5}$ \\
        \hline
 		    & $5 \times 10^{-4}$ & 2.02  & 2.12 & 0.68 & 11   & 10 & 10 & 10 & $3.1 \times 10^{-5}$ \\
 		0.5 & $1 \times 10^{-4}$ & 2.01  & 2.39 & 0.66 & 464  & 10 & 10 & 10 & $3.1 \times 10^{-5}$ \\
 		    & $5 \times 10^{-5}$ & 1.98  & 3.13 & 0.62 & 1288 & 10 & 10 & 10 & $3.2 \times 10^{-5}$ \\
 		\hline
            & $5 \times 10^{-4}$ & 2.01  & 2.94 & 0.67 & 10   & 10 & 10 & 10 & $3.0 \times 10^{-5}$ \\
 		0.9 & $1 \times 10^{-4}$ & 2.01  & 2.46 & 0.69 & 447  & 10 & 10 & 10 & $3.0 \times 10^{-5}$ \\
 		    & $5 \times 10^{-5}$ & 2.05  & 2.20 & 0.67 & 1866 & 10 & 10 & 10 & $3.0 \times 10^{-5}$ \\
 		\hline
 	\end{tabular}
 	\label{tab:mlmc_1e4_cell_8}
 \end{table}
 \begin{table}[h!]
 	\centering
 	\caption{Test 2. $F=\mathcal{F}_{\tau_{8,0}}$,  MLMC-HSM,   $K_{\ell} = 10^4$}
 	\begin{tabular}{|c|c||c|c|c|c|c|c|c|c|}
 		\hline
 		$c_2$ & $\varepsilon$     & $\alpha$ & $\beta$ & $\gamma$ & $N_0$ & $N_1$ & $N_2$ & $N_3$ & $\max_{\hat \ell}{W_{\hat{\ell}}}$\\
 		\hline
 		    & $5 \times 10^{-4}$ & 1.99  & 3.99 & 0.66 & 10   & 10 & 10 & 10 & $2.6 \times 10^{-5}$ \\
 		0.1 & $1 \times 10^{-4}$ & 2.00  & 2.66 & 0.63 & 430  & 10 & 10 & 10 & $2.6 \times 10^{-5}$ \\
 		    & $5 \times 10^{-5}$ & 1.99  & 2.55 & 0.69 & 1185 & 10 & 10 & 10 & $2.6 \times 10^{-5}$ \\
 		\hline
            & $5 \times 10^{-4}$ & 2.00  & 3.44 & 0.63 & 16   & 10 & 10 & 10 & $2.4 \times 10^{-5}$ \\
 		0.5 & $1 \times 10^{-4}$ & 2.00  & 2.57 & 0.68 & 410  & 10 & 10 & 10 & $2.4 \times 10^{-5}$ \\
 		    & $5 \times 10^{-5}$ & 2.01  & 3.35 & 0.70 & 1637 & 10 & 10 & 10 & $2.4 \times 10^{-5}$ \\
 		\hline
            & $5 \times 10^{-4}$ & 2.00  & 2.71 & 0.66 & 15   & 10 & 10 & 10 & $2.3 \times 10^{-5}$ \\
 		0.9 & $1 \times 10^{-4}$ & 2.00  & 3.13 & 0.66 & 446  & 10 & 10 & 10 & $2.3 \times 10^{-5}$ \\
 		    & $5 \times 10^{-5}$ & 2.00  & 3.07 & 0.67 & 1847 & 10 & 10 & 10 & $2.3 \times 10^{-5}$ \\
 		\hline
 	\end{tabular}
 	\label{tab:mlmc_1e4_cell_8_losm}
 \end{table}

\begin{table}[h!]
	\centering
	\caption{Test 2. $\{\mathcal{F}_{\tau_{i,0} }\}_{i=1}^{I_0}$,  MLMC-HQD,   $K_{\ell} = 10^4$}
	\begin{tabular}{|c|c||c|c|c|c|c|c|c|c|}
		\hline
		$c_2$ & $\varepsilon$     & $\min\alpha_i$ & $\min\beta_i$ & $\gamma$ & $N_0$ & $N_1$ & $N_2$ & $N_3$ & $\max{W_{\hat{\ell},i}}$\\
		\hline
		    & $5 \times 10^{-4}$ & 1.97  & 1.64 & 0.77 & 10   & 10 & 10 & 10 & $3.2 \times 10^{-5}$ \\
		0.1 & $1 \times 10^{-4}$ & 1.90  & 2.19 & 0.68 & 445  & 10 & 10 & 10 & $3.7 \times 10^{-5}$ \\
		    & $5 \times 10^{-5}$ & 1.87  & 1.58 & 0.70 & 1988 & 10 & 10 & 10 & $3.7 \times 10^{-5}$ \\
       \hline
		    & $5 \times 10^{-4}$ & 1.97  & 1.64 & 0.77 & 21   & 10 & 10 & 10 & $3.2 \times 10^{-5}$ \\
		0.5 & $1 \times 10^{-4}$ & 1.98  & 1.37 & 0.73 & 558  & 10 & 10 & 10 & $3.2 \times 10^{-5}$ \\
		    & $5 \times 10^{-5}$ & 1.93  & 1.69 & 0.74 & 1764 & 10 & 10 & 10 & $3.3 \times 10^{-5}$ \\
		\hline
           & $5 \times 10^{-4}$ & 1.98  & 1.14 & 0.73 & 33   & 10 & 10 & 10 & $3.1 \times 10^{-5}$ \\
		0.9 & $1 \times 10^{-4}$ & 1.96  & 1.39 & 0.75 & 493  & 10 & 10 & 10 & $3.2 \times 10^{-5}$ \\
		    & $5 \times 10^{-5}$ & 1.96  & 1.46 & 0.75 & 3893 & 10 & 10 & 10 & $3.2 \times 10^{-5}$ \\
		\hline
	\end{tabular}
	\label{tab:mlmc_1e4_cell_max_var}
\end{table}
\begin{table}[h!]
	\centering
	\caption{Test 2. $\{\mathcal{F}_{\tau_{i,0} }\}_{i=1}^{I_0}$,  MLMC-HSM,   $K_{\ell} = 10^4$}
	\begin{tabular}{|c|c||c|c|c|c|c|c|c|c|}
		\hline
		$c_2$ & $\varepsilon$     & $\min\alpha_i$ & $\min\beta_i$ & $\gamma$ & $N_0$ & $N_1$ & $N_2$ & $N_3$ & $\max{W_{\hat{\ell},i}}$\\
		\hline
		    & $5 \times 10^{-4}$ & 1.88  & 2.56 & 0.67 & 31   & 10 & 10 & 10 & $2.7 \times 10^{-5}$ \\
		0.1 & $1 \times 10^{-4}$ & 1.96  & 1.85 & 0.70 & 410  & 10 & 10 & 10 & $2.6 \times 10^{-5}$ \\
		    & $5 \times 10^{-5}$ & 1.86  & 2.91 & 0.70 & 2507 & 10 & 10 & 10 & $2.7 \times 10^{-5}$ \\
		\hline
            & $5 \times 10^{-4}$ & 1.98  & 1.87 & 0.68 & 21   & 10 & 10 & 10 & $2.5 \times 10^{-5}$ \\
		0.5 & $1 \times 10^{-4}$ & 1.98  & 2.38 & 0.66 & 512  & 10 & 10 & 10 & $2.5 \times 10^{-5}$ \\
		    & $5 \times 10^{-5}$ & 1.99  & 2.22 & 0.67 & 1742 & 10 & 10 & 10 & $2.5 \times 10^{-5}$ \\
		\hline
            & $5 \times 10^{-4}$ & 2.00  & 2.48 & 0.66 & 37   & 10 & 10 & 10 & $2.5 \times 10^{-5}$ \\
		0.9 & $1 \times 10^{-4}$ & 2.00  & 2.25 & 0.69 & 743  & 10 & 10 & 10 & $2.5 \times 10^{-5}$ \\
		    & $5 \times 10^{-5}$ & 2.00  & 2.33 & 0.68 & 2821 & 10 & 10 & 10 & $2.5 \times 10^{-5}$ \\
		\hline
	\end{tabular}
	\label{tab:mlmc_1e4_cell_max_var_losm}
\end{table}  
 \FloatBarrier
 
\subsection{Accuracy of Functional Estimator}

In the previous section, we presented an analysis of the convergence of MLHT algorithms with MLMC optimization, verifying the conditions of Theorem \ref{theorem}.
We now evaluate the accuracy of the methods and analyze whether the obtained estimations of functional reach the expected MSE.
We performed 10 independent runs of the MLHT algorithms for each $\varepsilon$ value.
The functional of interest is estimated using MC and the MLHT algorithms with MLMC optimization.
We calculate the MSE of the solution 
\begin{equation}
	MSE\big( \big< F_L \big> \big) = \mathbb{E}\Big[ \big( \big< F_L \big> - F^{ex} \big)^2\Big]   \, 
\end{equation}
using the reference value $F^{ex}= F[\phi^{ex}]$.
We solve Test 1 to compute $F=\mathcal{F}_D$ using $K_\ell = 10^4$ particle histories on each level.
For this problem, the reference values of the functionals are $\mathcal{F}_D^{ex} = 1.37293$ and $F_{\tau_8}^{ex} = 9.67931 \times 10^{-2}$.
The MSE is presented in Figures \ref{fig:test1_mlmc_mse_full} and \ref{fig:test1_mlmc_mse_cell_8} for both MLMC-HQD and MLMC-HSM methods, along with an MC estimate of the functional, for each of the 10 simulations and $\varepsilon = 1.0 \times 10^{-3}$ for $F = F_D$ and $\varepsilon = 1.0 \times 10^{-4}$ for $F = F_{\tau_8}$.
In addition, the average MSE errors are presented in Tables \ref{tab:mlmc_accuracy_full_range_func}-\ref{tab:mlmc_accuracy_tau_8_func_sm}.
The results show that, on average, the MLMC algorithm converges the MSE to the expected accuracy for most values of $\varepsilon$.
Comparing to MC, we find the average error is lower for the MLMC optimized methods for $\varepsilon \leq 1 \times 10^{-3}$ for $F=F_D$ and $\varepsilon < 1 \times 10^{-4}$ for $F = F_{\tau_8}$.
Although a promising result, we caution readers that only 10 samples were used to compute the average error.
We note that the MLMC-HSM algorithm failed the weak convergence check for $\varepsilon = 5\times 10^{-4}$ and $\mathcal{F} = \mathcal{F}_D$, meaning an additional level is needed to converge our functional, and is likely why the results show an average MSE greater than $\varepsilon$ in the case of MLMC-HSM.
\begin{table}[htb]
   	\caption{Test 1. average $MSE\big( \big< F_L \big> \big)$, $F=\mathcal{F}_D$ for MLMC-HQD \label{tab:mlmc_accuracy_full_range_func} }
   	\centering
   	\begin{tabular}{|c|c|c|c|}
   		\hline
   		$\varepsilon$ & MC  & MLMC-HQD & $\varepsilon^2$ \\
   		\hline
   		$5 \times 10^{-3}$ & $1.89\times 10^{-6}$ & $6.54\times 10^{-6}$ & $2.5 \times 10^{-5} $\\
   		$1 \times 10^{-3}$ & $6.37\times 10^{-7}$ & $3.39\times 10^{-7}$ & $1.0 \times 10^{-6} $\\
   		$5 \times 10^{-4}$ & $3.28\times 10^{-7}$ & $1.70\times 10^{-7}$ & $2.5 \times 10^{-7} $\\
   		\hline
   	\end{tabular}
\end{table}
 
\begin{table}[htb]
   	   	\caption{Test 1. average $MSE\big( \big< F_L \big> \big)$, $F=\mathcal{F}_D$ for MLMC-HSM}\label{tab:mlmc_accuracy_full_range_func_sm}
   	\centering
   	\begin{tabular}{|c|c|c|c|}
   		\hline
   		$\varepsilon$ & MC  & MLMC-HSM & $\varepsilon^2$ \\
   		\hline
   		$5 \times 10^{-3}$ & $1.74\times 10^{-6}$ & $6.60\times 10^{-6}$ & $2.5 \times 10^{-5} $\\
   		$1 \times 10^{-3}$ & $7.93\times 10^{-7}$ & $6.17\times 10^{-7}$ & $1.0 \times 10^{-6} $\\
   		$5 \times 10^{-4}$ & $6.24\times 10^{-7}$ & $2.62\times 10^{-7}$ & $2.5 \times 10^{-7} $\\
   		\hline
   	\end{tabular}
\vspace{1.5cm}
   	\caption{Test 1. average $MSE\big( \big< F_L \big> \big)$, $F=\mathcal{F}_{\tau_8}$ \label{tab:mlmc_accuracy_tau_8_func} for MLMC-HQD}
   	\centering
   	\begin{tabular}{|c|c|c|c|}
   		\hline
   		$\varepsilon$ & MC  & MLMC-HQD & $\varepsilon^2$ \\
   		\hline
   		$5 \times 10^{-4}$ & $6.90\times 10^{-8}$ & $1.03\times 10^{-7}$ & $2.5 \times 10^{-7} $\\
   		$1 \times 10^{-4}$ & $2.60\times 10^{-9}$ & $2.83\times 10^{-9}$ & $1.0 \times 10^{-8} $\\
   		$5 \times 10^{-5}$ & $7.50\times 10^{-10}$& $6.31\times 10^{-10}$& $2.5 \times 10^{-9} $\\
   		\hline
   	\end{tabular}
\vspace{1.5cm}
   	   	\caption{Test 1. average $MSE\big( \big< F_L \big> \big)$, $F=\mathcal{F}_{\tau_8}$ \label{tab:mlmc_accuracy_tau_8_func_sm}for MLMC-HSM}
   	\centering
   	\begin{tabular}{|c|c|c|c|}
   		\hline
   		$\varepsilon$ & MC  & MLMC-HSM & $\varepsilon^2$ \\
   		\hline
   		$5 \times 10^{-4}$ & $3.59\times 10^{-8}$ & $1.19\times 10^{-7}$ & $2.5 \times 10^{-7} $\\
   		$1 \times 10^{-4}$ & $5.43\times 10^{-9}$ & $4.50\times 10^{-9}$ & $1.0 \times 10^{-8} $\\
   		$5 \times 10^{-5}$ & $2.65\times 10^{-9}$ & $1.29\times 10^{-9}$ & $2.5 \times 10^{-9} $\\
   		\hline
   	\end{tabular}
\end{table}  

\begin{figure}[h!]
	\vspace{1cm}
	\begin{subfigure}{0.49\textwidth}
		\centering
		\includegraphics[width=\textwidth]{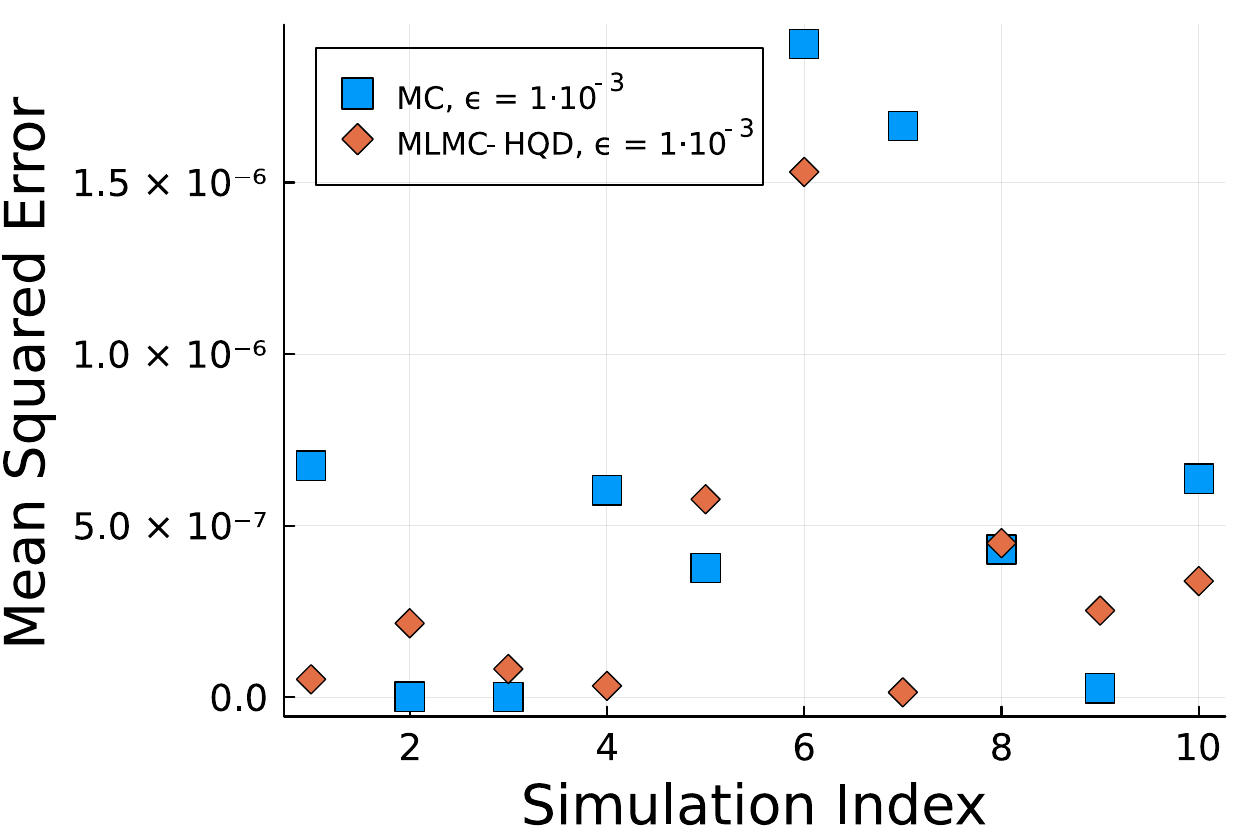}
		\caption{MLMC-HQD}
	\end{subfigure}
	\begin{subfigure}{0.49\textwidth}
		\centering
		\includegraphics[width=\textwidth]{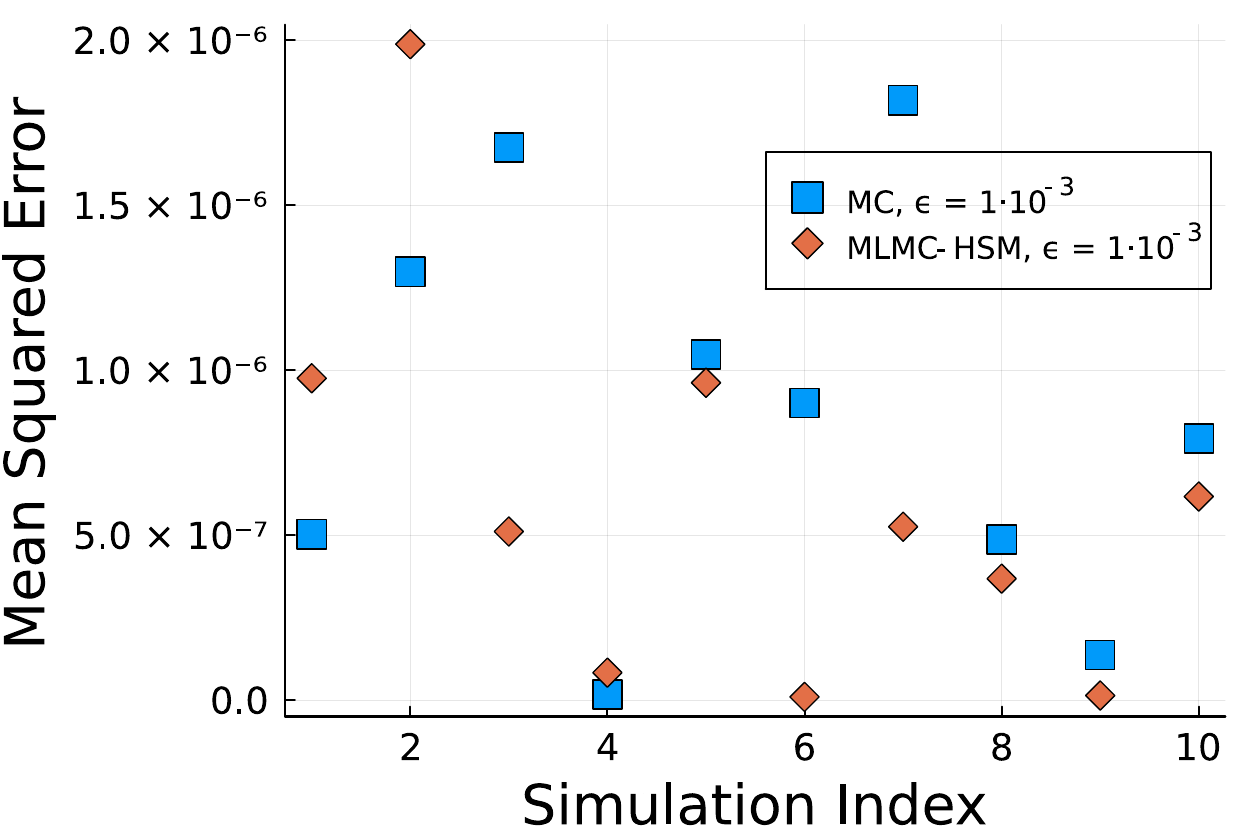}
		\caption{MLMC-HSM} 
	\end{subfigure}
	\caption{Test 1. MSE error in  the functional $F = F_D$ computed  by the MLMC-HQD and MLMC-HSM methods in each of 10 simulations with $\varepsilon=10^{-3}$} \label{fig:test1_mlmc_mse_full}
\end{figure}
\newpage
\begin{figure}[h!]
	\begin{subfigure}{0.49\textwidth}
		\centering
		\includegraphics[width=\textwidth]{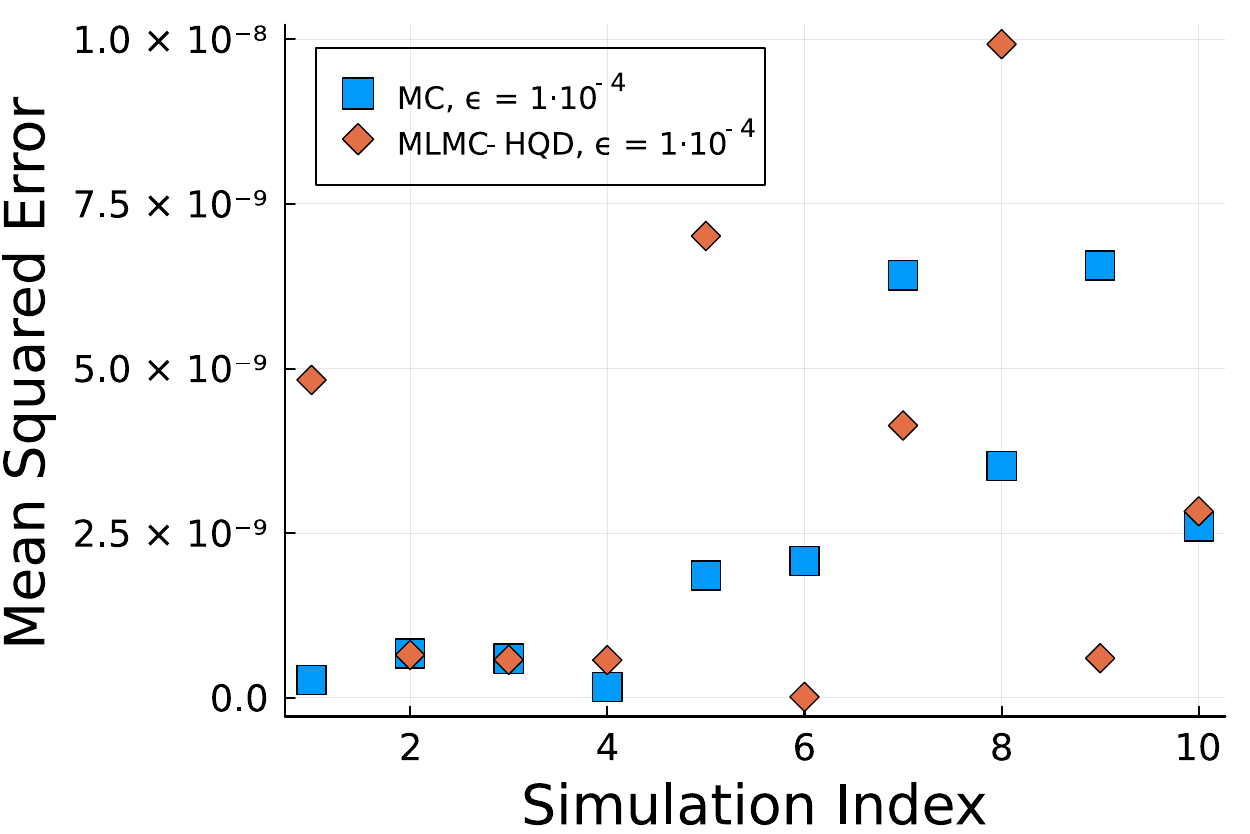}
		\caption{MLMC-HQD}
	\end{subfigure}
	\begin{subfigure}{0.49\textwidth}
		\centering
		\includegraphics[width=\textwidth]{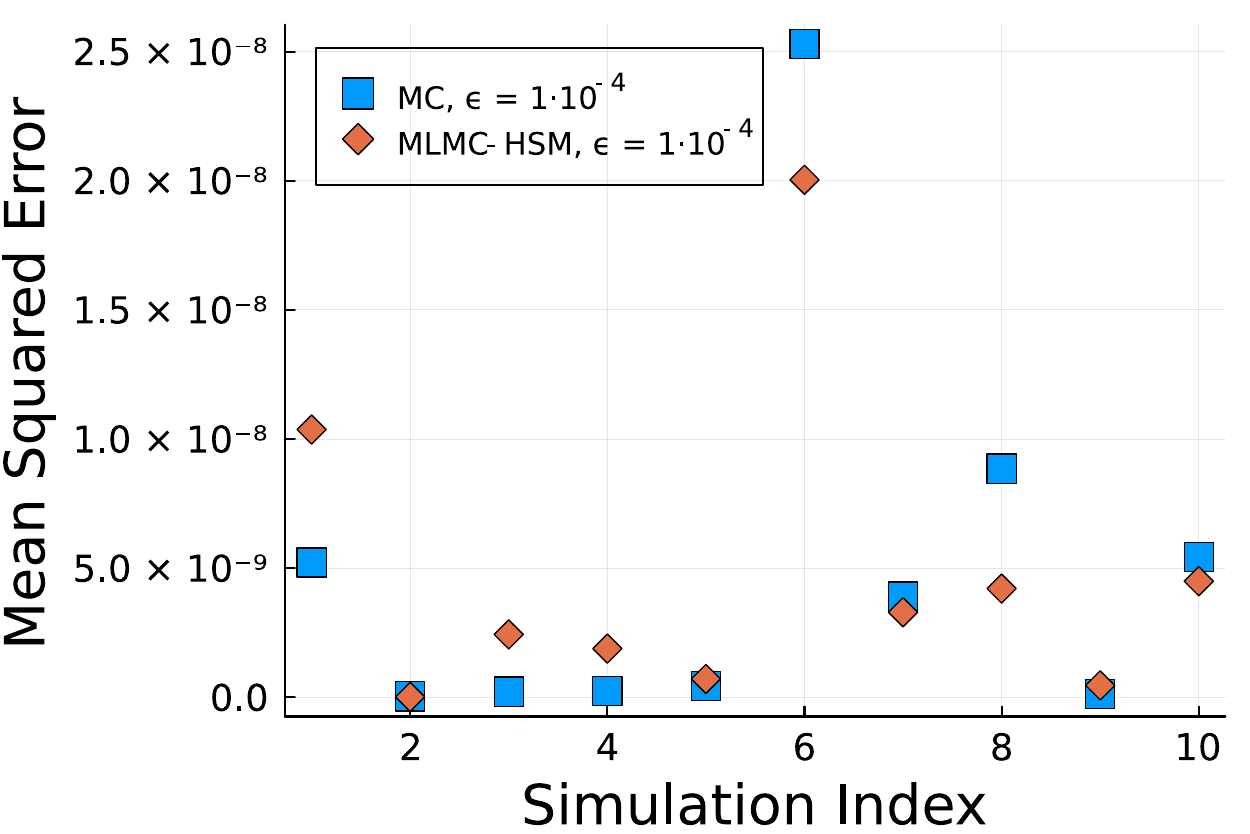}
		\caption{MLMC-HSM}
	\end{subfigure}
	\caption{Test 1. MSE error in  the functional $F_{\tau_8}$ computed  by the MLMC-HQD and MLMC-HSM methods in each of 10 simulations with $\varepsilon=10^{-3}$} \label{fig:test1_mlmc_mse_cell_8}
\end{figure}

\section{Conclusions } \label{sec:con}

In this paper, we presented novel MLHT methods for solving particle transport problems based on the fundamental ideas of MLMC and geometric multigrid in space.
The proposed MLTH is formulated using hybrid MC/deterministic numerical schemes, namely the HQD and HSM methods, which are defined by the LOQD and LOSM equations, respectively, and are discretized with $2^{nd}$-order FV schemes.

Analysis of the HQD and HSM methods showed that, with spatial refinement and increasing particle counts, the two methods exhibit similar convergence in the $L_2$ error norms of the scalar flux solution.
The developed MLHT algorithms reduce the magnitude of the correction functional as $\ell$ increases, indicating that the method converges with increasing solution fidelity.
The conditions of the MLMC theorem (Theorem \ref{theorem}) are met by both MLMC-HQD and MLMC-HSM algorithms.
To demonstrate convergence of the MLHT algorithm with MLMC optimization, the true MSE was evaluated for functionals and for a vector of functionals representing the solution over a spatial domain.
The average MSE was below the selected $\varepsilon^2$ for each set of results that met the criterion of Theorem \ref{theorem}.

The MSE was typically below $\varepsilon^2$ across most runs of the MLHT algorithm with MLMC optimization, indicating that the accuracy estimates provided by the MLHT algorithm are reasonable for the transport problems considered.
One possible reason the observed MSE is below $\varepsilon^2$ in some cases is the accuracy of the variance estimate for the functional.
The variance estimator can be noisy, leading to premature termination of MLHT sample generation when additional work should have been requested.
This is supported by the estimates of $\beta$ obtained, which showed significant variability across simulations.

The developed MLHT methodology led to the construction of an unbiased estimator for a functional of the scalar flux using a hybrid MC/deterministic technique that incurs discretization error.
Using an MLMC approach allows us to correct for the bias in the hybrid solution and potentially yields a more efficient solution, since most of the work is performed on the coarser computational grid, where the MC simulation and low-order solves require less compute time than on a finer grid. 
In addition, this algorithm does not preclude the use of variance-reduction techniques in Monte Carlo particle simulations; for example, we used implicit capture with Russian Roulette. 
The effects of variance reduction on the MLHT algorithm will be beneficial, since the number of request simulations at each computational level will be reduced by roughly the same amount as the variance is reduced.
We observed this effect when comparing the results for $10^3$ vs $10^4$ particle histories in Test 1: a $10$-fold reduction in the number of particles increases the variance by a factor of $10$, yielding roughly $10$ times as many samples requested for $K_\ell = 10^3$.

Future work will include examining MLHT algorithms with HMCD methods having a higher order of spatial convergence using high-order discretization schemes.
Such MLHT methods have the potential to converge to higher accuracy with fewer computational levels than the second-order schemes we examined here. 
MLHT methods with advanced prolongation operators will be considered.  
The proposed methodology can be extended to develop hybrid transport calculations based on Multi-Fidelity Monte Carlo and approximate control-variate approaches \citep{MFMC-Peherstorfer-2016, MFmethods-siam-review-2018, ACV-Bomarito-2022}.

\section*{ACKNOWLEDGEMENTS}

This work was supported by the Center for Exascale Monte Carlo Neutron Transport (CEMeNT), a PSAAP-III project funded by the Department of Energy, grant number DE-NA003967.
A part of the work of the first author (VNN) was supported under a University Nuclear Leadership Program Graduate Fellowship, grant number DE-FOA-0002265. Any opinions, findings, conclusions or recommendations expressed in this publication are those of the author(s) and do not necessarily reflect the views of the Department of Energy Office of Nuclear Energy. 

\section*{Declarations of Interest}
The authors report there are no competing interests to declare.

\section*{Author Contribution Statement}
VNN was involved with the formal analysis, funding acquisition, methodology, software, visualization, writing: original draft, validation, and investigation.
DYA was involved with conceptualization, formal analysis, funding acquisition, methodology, project administration, supervision, writing: review \& editing, resources, validation, and investigation.
Both VNN and DYA agree to be accountable for all aspects of the work.
\bibliographystyle{elsarticle-num}
\bibliography{novellino_anistratov_mlmc_arxiv_2}

\clearpage
\appendix
\section{Reference Solution} \label{sec:appendix_a}

To compute the reference solution for Test 1, we use the deterministic QD method for solving the transport equation \citep{gol'din-cmmp-1964}. 
The high-order transport equation is discretized in angle using the method of discrete ordinates and approximated in space using step characteristics. 
LOQD equations are discretized with the 2$^{nd}$ order finite volume scheme described in Section \ref{sec:LOQD} \citep{dya-vyag-vant-1986,dya-vyag-ttsp}. 
The numerical integration method for the angular variable is based on a second-order technique.
The phase-space grid is refined uniformly. 
Results on each of three successive grids are extrapolated by means of Aikten's method \citep{Sidi_2003, dahlquist2008numerical, Ganapol2013}. 
The refinement proceeds until the Aitken-extrapolated numerical solution converges in phase space for the given tolerance.  
Table \ref{tab:Reference solution for Test 1} shows 7 significant digits of the numerical reference solution of Test 1 in terms of the cell-average scalar fluxes and the value at the boundary on the spatial mesh with $I=2^7$. 
The solution is symmetric in space; therefore, we present the results for the left half of the domain, i.e., $x \in [0,0.5]$ cm.

{\footnotesize
\begin{longtable}{|c|c|c|c|}
    \caption{Reference solution for Test 1 on the spatial mesh with $\Delta x = 2^{-7}$ presented as the vector of cell-average scalar fluxes, including the value at the boundary.}
    \label{tab:Reference solution for Test 1}
     \\ \hline
    $x_i$ & $\phi^{ex}$ & $x_i$ & $\phi^{ex}$ \\
    \hline
	$0.000 $              & $0.923204$ & &              \\
	$3.906 \times 10^{-2}$& $0.945787$ &$2.539 \times 10^{-1}$& $1.432468$ \\
	$1.172 \times 10^{-2}$& $0.979368$ &$2.617 \times 10^{-1}$& $1.440153$ \\
	$1.953 \times 10^{-2}$& $1.007905$ &$2.695 \times 10^{-1}$& $1.447552$ \\
	$2.734 \times 10^{-2}$& $1.033603$ &$2.773 \times 10^{-1}$& $1.454670$ \\
	$3.516 \times 10^{-2}$& $1.057304$ &$2.852 \times 10^{-1}$& $1.461510$ \\
	$4.297 \times 10^{-2}$& $1.079455$ &$2.930 \times 10^{-1}$& $1.468078$ \\
	$5.078 \times 10^{-2}$& $1.100336$ &$3.008 \times 10^{-1}$& $1.474376$ \\
	$5.859 \times 10^{-2}$& $1.120137$ &$3.086 \times 10^{-1}$& $1.480408$ \\
	$6.641 \times 10^{-2}$& $1.138997$ &$3.164 \times 10^{-1}$& $1.486177$ \\
	$7.422 \times 10^{-2}$& $1.157020$ &$3.242 \times 10^{-1}$& $1.491687$ \\
	$8.203 \times 10^{-2}$& $1.174288$ &$3.320 \times 10^{-1}$& $1.496939$ \\
	$8.984 \times 10^{-2}$& $1.190866$ &$3.398 \times 10^{-1}$& $1.501937$ \\
	$9.766 \times 10^{-2}$& $1.206809$ &$3.477 \times 10^{-1}$& $1.506683$ \\
	$1.055 \times 10^{-1}$& $1.222161$ &$3.555 \times 10^{-1}$& $1.511179$ \\
	$1.133 \times 10^{-1}$& $1.236960$ &$3.633 \times 10^{-1}$& $1.515427$ \\
	$1.211 \times 10^{-1}$& $1.251239$ &$3.711 \times 10^{-1}$& $1.519429$ \\
	$1.289 \times 10^{-1}$& $1.265026$ &$3.789 \times 10^{-1}$& $1.523188$ \\
	$1.367 \times 10^{-1}$& $1.278345$ &$3.867 \times 10^{-1}$& $1.526703$ \\
	$1.445 \times 10^{-1}$& $1.291217$ &$3.945 \times 10^{-1}$& $1.529978$ \\
	$1.523 \times 10^{-1}$& $1.303663$ &$4.023 \times 10^{-1}$& $1.533014$ \\
	$1.602 \times 10^{-1}$& $1.315699$ &$4.102 \times 10^{-1}$& $1.535811$ \\
	$1.680 \times 10^{-1}$& $1.327339$ &$4.180 \times 10^{-1}$& $1.538371$ \\
	$1.758 \times 10^{-1}$& $1.338598$ &$4.258 \times 10^{-1}$& $1.540695$ \\
	$1.836 \times 10^{-1}$& $1.349489$ &$4.336 \times 10^{-1}$& $1.542785$ \\
	$1.914 \times 10^{-1}$& $1.360021$ &$4.414 \times 10^{-1}$& $1.544639$ \\
	$1.992 \times 10^{-1}$& $1.370207$ &$4.492 \times 10^{-1}$& $1.546261$ \\
	$2.070 \times 10^{-1}$& $1.380054$ &$4.570 \times 10^{-1}$& $1.547650$ \\
	$2.148 \times 10^{-1}$& $1.389571$ &$4.648 \times 10^{-1}$& $1.548806$ \\
	$2.227 \times 10^{-1}$& $1.398767$ &$4.727 \times 10^{-1}$& $1.549730$ \\
	$2.305 \times 10^{-1}$& $1.407648$ &$4.805 \times 10^{-1}$& $1.550423$ \\
	$2.383 \times 10^{-1}$& $1.416221$ &$4.883 \times 10^{-1}$& $1.550885$ \\
	$2.461 \times 10^{-1}$& $1.424493$ &$4.961 \times 10^{-1}$& $1.551116$ \\
\hline
\end{longtable}
}

\end{document}